\documentclass[12pt]{article}
\usepackage{amsthm,amsmath,epsf,amssymb,xy}
\xyoption{all} 
\newcommand{\cA}{{\mathcal A}}

\newcommand{\CC}{{\mathbb C}}
\newcommand{\Ca}{{\mathcal C}}
\newcommand{\Cpr}{{{\mathcal C}^{\prime \,\vee }}}
\newcommand{\Der}{{\bf {\rm D} }}
\newcommand{\Ea}{{\mathcal E}}
\newcommand{\EE}{{\mathbb E}}
\newcommand{\Fa}{{\mathcal F}}

\newcommand{\Ha}{{\mathcal H}}
\newcommand{\Hom}{{ {\rm Hom} }}
\newcommand{\Ia}{{\mathcal I}}

\newcommand{\La}{{\mathcal L}}
\newcommand{\LL}{{\mathbb L}}
\newcommand{\Ltensor}{\mathbin{\overset{\mathbf L}\otimes}}
\newcommand{\Mrm}{{\rm M}}
\newcommand{\Ma}{{\mathcal M}}
\newcommand{\Nrm}{{\rm N}}
\newcommand{\Na}{{\mathcal N}}
\newcommand{\ovl}[1]{{\overline{#1}}}
\newcommand{\Oa}{{\mathcal O}}
\newcommand{\ov}{{\overline{v}}}
\newcommand{\tv}{{\widetilde{v}}}
\newcommand{\PP}{{\mathbb P}}
\newcommand{\Pa}{{\mathcal P}}
\newcommand{\QQ}{{\mathbb Q}}
\newcommand{\rmu}{{\mu^*}}
\newcommand{\rnu}{{\nu^*}}
\newcommand{\Ra}{{\mathcal R }}
\newcommand{\RR}{{\mathbb R}}
\newcommand{\Rz}{{\rm Res}}
\newcommand{\tmu}{{\tilde{\mu}}}
\newcommand{\tnu}{{\tilde{\nu}}}
\newcommand{\Ta}{{\mathcal T}}
\newcommand{\Td}{{\rm Todd}}
\newcommand{\Sa}{{\mathcal S}}
\newcommand{\ZZ}{{\mathbb Z}}

\newtheorem*{conj}{Conjecture {\rm (M. Kontsevich)}}
\newtheorem{thm}{Theorem}[section]
\newtheorem{newconj}[thm]{Conjecture}
\newtheorem{corollary}[thm]{Corollary}
\newtheorem{definition}[thm]{Definition}
\newtheorem{example}[thm]{Example}
\newtheorem{proposition}[thm]{Proposition}
\newtheorem{lemma}[thm]{Lemma}
\newtheorem{condition}[thm]{Condition}

\theoremstyle{definition}
\newtheorem{remark}[thm]{Remark}
\newtheorem*{conv}{Convention}

\newif\iffigs\figstrue

\numberwithin{equation}{section}

\begin{document}


\begin{center}
{\Large Hypergeometric Functions and Mirror Symmetry in 
Toric Varieties
}

\end{center}
\bigskip
\centerline{Richard Paul Horja}
\bigskip
\centerline{\it Department of Mathematics}
\centerline{\it Duke University, Box 90320}
\centerline{\it Durham, NC 27708-0320, U.S.A.}
\centerline{\tt horja@math.duke.edu}

\begin{abstract}
We study aspects related to Kontsevich's homological mirror symmetry 
conjecture \cite{Kont1} in the case of Calabi--Yau complete intersections 
in toric varieties. In a 1996 lecture, Kontsevich \cite{Kont2}
indicated how his proposal implies that the groups of automorphisms of 
the two types of categories involved in the homological mirror 
symmetry conjecture should also be identified. Our results 
provide an explicit geometric construction of the correspondence between the 
automorphisms of the two types of categories in some cases, as well as a 
conjectural proposal for the general case.

We compare the monodromy calculations for the Picard--Fuchs system
associated with the periods of a Calabi--Yau manifold $M$ with the
algebro-geometric computations of the cohomology action
of Fourier--Mukai functors on the bounded derived category 
of coherent sheaves on the mirror Calabi--Yau manifold $W.$
We obtain the complete dictionary between the two sides 
for the one complex parameter case of Calabi--Yau complete 
intersections in weighted projective spaces, as well as 
for some two parameter cases. We also find the complex 
of sheaves on $W \times W$ that corresponds to a loop in the 
moduli space of complex structures on $M$ induced by 
a phase transition of $W.$
\end{abstract}

\section{Introduction}
\label{cha:introd}

Mirror symmetry was discovered some time 
ago by string theorists as an equivalence between two physical theories
associated to pairs of families of Calabi--Yau manifolds. The best known
manifestation of the phenomenon consisted in the 
predictions made by a group of physicists \cite{CDGP} about the number 
of rational curves on Calabi--Yau manifolds. In 1993, Witten 
\cite{Witten1} introduced the gauged linear sigma model and studied 
the different realizations of the physical theory, which he referred to as 
{\it phase transitions} in the theory. More detailed investigations 
of the mechanism through which mirror symmetry identifies the 
parameter spaces  corresponding to the two sides of the story  
(the K\"ahler moduli and the complex moduli, respectively)
were pursued by Aspinwall, Greene and Morrison 
\cite{AGM1},\cite{AGM2},\cite{AGM3}.  
In his ICM 1994 presentation \cite{Kont1}, Kontsevich proposed 
{\it the homological mirror symmetry conjecture} as
a new mathematical interpretation of mirror symmetry through
an equivalence between ``Fukaya's $A_\infty$ category'' of a Calabi--Yau 
manifold $M$ and the bounded derived category of coherent sheaves $\Der(W)$
of the mirror Calabi--Yau manifold $W.$ Although in the general 
Calabi--Yau case,
Fukaya's category remains quite elusive, in a 1996 lecture, Kontsevich
\cite{Kont2} indicated how his proposal implies that the 
groups of automorphisms of the two types of categories should also
be identified under mirror symmetry. Our
focus is not on the categories themselves, but 
rather on the identification of their automorphisms groups
in the case when the mirror Calabi--Yau
varieties $M$ and $W$ are complete intersections in toric
varieties. The braid group actions on derived categories of coherent
sheaves and their conjectural mirror symmetric relationship
with the symplectic generalized Dehn twists of Seidel \cite{seidel}
have been analyzed recently by Seidel and Thomas \cite{thomas}, 
\cite{ST}. 

The main results of the present work describe the correspondence 
between the automorphisms groups in the 1-parameter and some 
2-parameter cases, as well as a geometric construction of complexes in 
$\Der (W \times W)$ whose induced actions on $H^*(W,\CC)$ match the 
monodromy action on $H^*(M,\CC)$ provided by the actions of loops around 
the various components 
of the discriminant locus in the moduli space of complex structures
on $M$ (in the general case). The construction establishes an explicit 
realization in the 
toric case of the program presented by Kontsevich in his lecture 
\cite{Kont2}. Our method uses tools developed in the ``early'' days of mirror 
symmetry, such as the machinery of Gelfand--Kapranov--Zelevinsky 
hypergeometric systems \cite{GKZ1}, \cite{Bat2} and the description of the
phase transitions in type II string theory obtained 
by Witten \cite{Witten1}, and Aspinwall, Greene and Morrison 
\cite{AGM1}. Their investigation offers a detailed picture of 
how the variations of the complexified K\"ahler class on $W$ (inducing,
in particular, different birational models, or ``phases'', of $W$) are
identified with the variations of complex structures on $M.$ 
As a consequence of the relative simple structure of the (special) Lagrangian 
submanifolds in the elliptic curve and K$3$ surface situations,
several authors were able to prove special cases of 
Kontsevich's conjectures (see \cite{PolZas} 
for the case of elliptic curves, and \cite{BBS} for the case of K$3$
surfaces).

The string theoretic origins of the equivalence
require that the Calabi--Yau structure on a manifold $X$
should be parametrized by a choice
of a complex structure with trivial canonical bundle, a K\"ahler form
$J,$ and an element $B$ in $H^2(X;\RR)/H^2(X;\ZZ)$ (the so-called B-field).
The complexified K\"ahler form is then defined to be $K = B+ iJ$. 
According to Yau's solution of the Calabi conjecture, the complex structure 
and the K\"ahler form determine uniquely a Ricci-flat metric. It has
long been stated by string theorists that the mirror manifolds
$W$ and $M$ are just different manifestations of the same physical
theory. The automorphisms of the two types of categories
studied in this work provide a way of comparing the effects on the 
physical theory of the variations of the different parts of the parameter 
spaces of mirror Calabi--Yau structures.

Starting with the work of Mukai \cite{Mukai1}, the actions 
of functors on derived categories of coherent sheaves have
received much attention. Under some additional assumptions, Bondal and 
Orlov \cite{BO} and 
Orlov \cite{O} showed that for a smooth projective variety $W$  
any automorphism of $\Der (W)$ is given by an object $\Ea^\bullet$ in  
$\Der(W \times W)$. For every such object,  an automorphism
functor of $\Der (W)$ is defined by 
\begin{equation}
\label{eq:auto}
\Phi_{\Ea^\bullet}(\cdot):={\bf \rm R}p_{2_*}
(\Ea^\bullet \Ltensor p_1^*(\cdot)),
\end{equation}
where $p_i : W \times W \to W$, $i=1,2$ are projections on the 
first and second factor, respectively, and the functors involved
are the corresponding right and left derived functors. 
Due to the analogy with 
integral operators, $\Ea^\bullet$ is called a ``kernel'' and  the 
transformation is known as a ``Fourier--Mukai functor''. By 
considering the induced action on Chern characters of coherent sheaves, 
these equivalences define a group $G_1$ of automorphisms of the cohomology 
of $W.$ 

On the other hand, since the objects of Fukaya's category are 
Lagrangian submanifolds of $M$ endowed with unitary local systems,
some natural automorphisms of the category 
should come from symplectomorphisms of the manifold $M.$
In particular, any loop 
in the moduli space of complex structures on $M$ should define such
an equivalence. The induced action on $H^{*}(M,{\mathbb C})$
determines a group $G_2$ of automorphisms of the cohomology of $M.$

We now state a version of the conjecture presented
by Kontsevich in his lecture \cite{Kont2}.

\begin{conj} The actions of the groups $G_1$ and $G_2$
on the cohomologies of mirror Calabi--Yau manifolds $W$ and $M$ are 
interchanged under mirror symmetry. In particular, the loops in 
the  moduli space of complex structures of $M$ are identified with
kernels in $\Der(W \times W).$ 
\label{conj:kont}
\end{conj}
\noindent
Among other very interesting considerations, 
Kontsevich supported this conjecture with an explicit calculation for the 
1-parameter case, in particular for the 
quintic threefold and its mirror; 
we reproduce that calculation in section \ref{chap:1par}.  

Although the conjecture is not toric in nature, we need 
toric methods to be able to explicitly compute the action of the monodromy
representation. In fact, in chapter \ref{chap:4}
we state a refinement of Kontsevich's proposal (conjecture \ref{newconj}) 
specific to the case of Calabi-Yau complete intersections in toric varieties. 

We now review the plan of this work. The required 
concepts are reviewed in chapter \ref{chap:gengen}. The exposition is an 
adaptation for the present purposes of the vast literature available 
in the field. 

The analytic continuation computations are performed in chapter
\ref{chap:mon}. Theorem \ref{gensol} and corollary \ref{cor:contform}
present a recursive construction 
of the solutions to the Picard--Fuchs system associated with
the periods of the mirror Calabi--Yau manifold $M.$ The results provide an 
integral representation that captures the complete set of solutions to the 
Picard--Fuchs system 
around {\it any} toric boundary point in the moduli space of complex 
structures of $M$ (previous results \cite{yau}, \cite{Stienstra} contain 
explicit power series representations around a maximal unipotent point in 
the moduli space of complex structures).
The analytic continuation is achieved essentially  
by using the classical Mellin--Barnes integral representations method 
for the Gelfand--Kapranov--Zelevinsky hypergeometric series. 
The content of theorem \ref{thm:mon} is 
a monodromy formula for a particular type of loop $\gamma_F$
(previously studied by Aspinwall, Greene and Morrison in \cite{AGM3})
corresponding to an edge $F$ in the secondary polytope of Gelfand, Kapranov and 
Zelevinsky \cite{GKZ2} associated to the toric situation.

In chapter \ref{chap:4}, the explicit analytic computations of 
the monodromy representation are compared with the algebro-geometric 
calculations of the cohomology actions of kernels. 
The complete identification between kernels in $\Der(W \times W)$
and loops in the moduli space of complex structures on $M$
is established in theorem \ref{thm:1par}
for the 1-parameter case of smooth Calabi--Yau complete 
intersections in weighted projective spaces, and in theorem \ref{thm:2par}
for some 2-parameter cases
generalizing examples extensively studied in the physics literature
(study begun in \cite{2par}). 

A particularly interesting case is the 
automorphism of $\Der(W)$ corresponding to a phase transition from the 
smooth Calabi--Yau variety $W$ to an ``adjacent phase''. Such a phase 
transition corresponds to some edge $F$ in the secondary polytope associated 
to the toric context (an example of this situation is provided by a birational 
transformation of $W$). A comparison with the monodromy action 
of the loop $\gamma_F$ on $H^*(M,\CC)$ allows us to give 
a geometric construction  of the complex $\Ea^\bullet(F)$ in $\Der(W \times W)$
(definition \ref{def:bigdef})
whose cohomology action matches the monodromy action of the loop $\gamma_F$
(theorem \ref{thm:genthm}). We also show how the phase transitions induced by 
edges in the secondary polytope 
are naturally associated to elementary contractions of the ambient toric 
variety (described by Reid in \cite{reid}) in the sense of Mori theory. 
The toric transition from $X$ to the ``adjacent'' toric phase $X'$
is induced by the morphisms $f : X \to 
\ovl{X}$ and $f' : X' \to \ovl{X}.$ If $E_0$ is the exceptional 
locus of $f$ in the toric variety $X,$ and $E$ the intersection
of $E_0$ with $W,$ $E:=W \cap E_0,$ 
let $Z:=f(E)$ be the image of $E$ under the morphism $f.$
The definition \ref{def:bigdef} of the 
complex $\Ea^\bullet(F) \in \Der(W \times W)$ 
involves a Koszul resolution associated to the complete intersection
$E \subset W.$ The action on $H^*(W,\CC)$ induced by the associated 
Fourier--Mukai functor $\Phi_{\Ea^\bullet(F)}$ is 
essentially determined by a sheaf supported on $E \times_Z E.$
While this functor has the expected cohomological properties, we conjecture 
its invertibility, but we do not prove it here. A proof of the more general 
toric proposal (conjecture \ref{newconj}) identifying the kernel 
corresponding to any face of the secondary polytope requires 
an extension of the  monodromy 
calculation of theorem \ref{thm:mon}. It is tempting to speculate that a 
general formula of this type would involve the use of multiple Mellin--Barnes 
integral representations, but the analytical details of such an approach 
seem to be quite intricate.

After the advent of D-branes in string theory, Strominger, Yau and 
Zaslow \cite{SYZ} proposed in 1996 a geometrical mirror symmetry construction.
In their interpretation, mirror 
symmetry is described by means of compactifications of toroidal 
special Lagrangian submanifolds (see also \cite{M}, \cite{gross}
for further investigations of this conjecture). Although 
the interplay between the different points of view is not yet clarified, 
an appropriate interpretation of the results of this work
offers some clues about the Strominger, Yau and Zaslow 
geometrical mirror symmetry conjecture.
One prediction of the SYZ proposal is a geometric map between the different 
kinds of branes (cycles) in type IIA and type
IIB string theories (this issue has been the subject of recent 
investigations by string theorists \cite{doug}). In this spirit,
the explicit identification between kernels and loops explored here
can be equivalently viewed as a dictionary between complexes
of sheaves on $W \times W$ and vanishing cycles on
the mirror manifold $M$ corresponding to the various components of 
the discriminant locus in the moduli space of complex structures on $M.$
For example, our monodromy calculations in the one and 
two parameter cases confirm the prediction that 
the structure sheaf of the Calabi--Yau manifold $W$ has to be 
mirror symmetric to the vanishing cycle at the conifold locus
of $M.$

This dictionary should also have applications when used in conjunction 
with the recent approaches to the problem of counting special Lagrangian
submanifolds \cite{joyce} and may allow us to partially overcome 
the difficult existence problem for special Lagrangian 
submanifolds of general Calabi--Yau manifolds (although, see recent 
work of Bryant \cite{bryant}).
Another question arises in connection with the so-called 
``local mirror symmetry'' \cite{CKYZ} and 
its consequences for field theory \cite{KatzMayrVafa}. Is there a similar 
identification in this case? In a different realm, the monodromy calculations
of chapter \ref{chap:mon} have the flavor of the vertex operator calculations
of the type described for example in \cite{boris}. It would be interesting
to see if the vertex algebras formalism can be applied in our context. 

Another related issue is to try to understand in more detail 
the relationship between this work and the theory of exceptional 
collections of sheaves and braided automorphisms of derived categories 
\cite{Mukai1}, \cite{rudakov}, \cite{deligne}, \cite{zaslow}
and especially the work of Seidel and Thomas 
\cite{thomas}, \cite{ST}. Kontsevich's proposal 
implies that there should be an action on $\Der (W)$ of the fundamental 
group of the complement of the discriminant locus in the moduli space of 
complex structures on $M.$ Zariski \cite{Zariski} pioneered the study of
the fundamental group of the complement of a complex curve in $\CC^2,$
and the braid monodromy technique has been employed since then
to obtain detailed information about this fundamental group \cite{moi}.
On general grounds, in the cases where we proved the
complete equivalence between kernels and loops, the braid monodromy action
on the fundamental group induces a braid group action on $\Der(W).$ The braid
group is also a frequent occurrence in the study of hypergeometric functions
which, not surprisingly, play a distinguished r\^{o}le in this work.
It is an intriguing problem (at least for the author) to try to understand 
if our results could shed some light on how these concepts are tied 
together.

{\bf Acknowledgements:} This paper is 
part of the author's Ph.D. dissertation written 
at Duke University. I am grateful to my advisor Professor David Morrison
for his valuable guidance and support. I would like to thank 
Professors Paul Aspinwall, Ronen Plesser, Les Saper and Mark Stern for
many helpful discussions. 

\newpage

\section{Toric Varieties and GKZ Systems}
\label{chap:gengen}
\subsection{Generalities about Toric Varieties}
\label{chap:toricgen}
Our main reference for this section will be \cite{CK}. We will only
sketch the main facts about toric varieties, with an emphasis
on the notions that we will need most. For more
details on the definitions and proofs, see also \cite{Dani}, \cite{Fulton}. 

Let $\Nrm$ be the integral lattice ${\ZZ^n}$ in $\Nrm_{\RR}\cong {\RR^n},$ 
and let $\Mrm={\rm Hom}(\Nrm, \ZZ)$ be its dual. Denote by 
$\langle \cdot,\cdot \rangle : \Mrm \times \Nrm \to \ZZ$ the canonical 
non-degenerate pairing. A {\bf rational
polyhedral cone} $\sigma$ in $\Nrm_{\RR}$ is a subset of the form 
\begin{equation}
\sigma=\{\sum_{i=1}^{s} \lambda_i u_i, \lambda_i \geq 0 \},
\end{equation}
where $u_1, \ldots , u_s \in \Nrm$. We say that $\sigma$ is 
{\bf strongly convex} if it does not contain any line. Every cone 
$\sigma$ has a {\bf dual cone} $\check{\sigma}$ defined by 
\begin{equation}
\check{\sigma}=\{m \in \Mrm_{\RR} : \langle m,v \rangle  \geq 0, 
\hbox{for all} \ v \in \sigma\}.
\end{equation}

A {\bf fan} $\Sigma$ in $\Nrm_{\RR}$ is a finite
collection of strongly convex rational polyhedral cones such that (i) the 
faces of cones in $\Sigma$ are in $\Sigma$ and (ii) if $\sigma, \sigma' \in
\Sigma$, then $\sigma \cap \sigma'$ is face of each. For each $d$,  
$\Sigma(d)$ denotes the collection of $d$-dimensional cones in 
$\Sigma$. For $d=1$ we obtain the set $\Sigma(1)$ of {\bf rays} or 
{\bf primitive vectors}. The subset $|\Sigma| \subset \Nrm_{\RR}$ obtained as 
the union of all cones $\sigma \in \Sigma$ is called the {\bf support}
of $\Sigma.$

Given a fan $\Sigma$, each cone $\sigma \in \Sigma$ gives the affine
toric variety 
\begin{equation}
X_{\sigma}={\rm Spec}(\CC[\Mrm \cap \check{\sigma}]),
\end{equation}
where $\CC[M \cap \check{\sigma}]$ is the $\CC$-algebra with generators
$\chi ^m$ for each $m \in \Mrm \cap \check{\sigma}$ and relations $\chi ^m
\chi ^{m'}=\chi ^{m+m'}$. The toric variety is then obtained from 
these affine pieces by gluing together $X_{\sigma}$ and $X_{\sigma'}$
along $X_{\sigma \cap \sigma'}$.

The properties of the fan $\Sigma$ determine the geometry of the toric
variety $X_{\Sigma}$. For example, $X_{\Sigma}$ is  
compact if and only if the support of the fan $\vert \Sigma \vert
= \Nrm_{\RR}$. The smoothness of $X_{\Sigma}$ is equivalent to the fact that
each cone in $\Sigma$ is spanned by a part of a $\ZZ$-basis of the lattice 
$\Nrm$; if the fan is only simplicial the corresponding toric variety has
orbifold singularities. 

The toric variety $X_{\Sigma}$ is a disjoint union of its orbits
under the action of the torus $(\CC^*)^n$. Namely, for each 
$k$-dimensional cone $\sigma \in \Sigma$ there is a $k$-codimensional
open orbit $O_\sigma$. The closure of this orbit is denoted by
$V(\sigma)$ and it is a toric variety itself obtained as the disjoint 
union of all the orbits $O_\gamma$ for which $\gamma$ contains $\sigma$ 
as a subcone. Consequently the toric structure of $V(\sigma)$
is given by a fan $\Sigma(\tau)$ in $(\Nrm/\Nrm_{\tau})_\RR,$ where 
$\Nrm_{\tau}$ is the sublattice of $\Nrm$ generated (as a group)
by $\tau \cap \Nrm.$ The fan $\Sigma(\tau)$ consists of cones obtained as 
images under the projection $\Nrm \to \Nrm/\Nrm_{\tau}$ 
of all the cones in $\Sigma$ that have $\tau$ as a subcone. For proofs
of these facts, see, for example, section 3.1 in \cite{Fulton}. 

A lattice homomorphism $f : (\Nrm, \Sigma) \to (\Nrm', \Sigma')$ is
called a {\bf map of fans} if, for each cone $\sigma$ in $\Sigma,$ 
there is some cone $\sigma'$ in $\Sigma'$ such that 
$f(\sigma) \subset \sigma'.$ A map of fans determines a 
morphism $f_*$ of the corresponding toric varieties
$f_* : X_{\Sigma} \to X_{\Sigma'}.$ If 
$f^{-1}(\vert \Sigma' \vert) = (\vert \Sigma \vert),$ 
then $f_*$ maps any closed toric orbit  
$V(\sigma) \subset X_{\Sigma}$ 
onto some closed toric orbit $V(\sigma') \subset X_{\Sigma'},$ where
$\sigma'$ is the smallest cone of $\Sigma'$ that
contains $f(\sigma)$ (see, for example, page 56 in \cite{Fulton}
for an explicit statement of this property). 

Each ray $v_i \in \Sigma(1)$ determines a codimension 1
orbit. The corresponding closure of the orbit in $X_\Sigma$ is a 
toric divisor $D_i$. The divisors $D_i$ generate the group of Weil
divisors $A_{n-1}(X).$
A Weil divisor $D=\sum m_i D_i$ is Cartier
if and only if for each $\sigma \in \Sigma,$ there is $m_\sigma
\in M$ such that $\langle m_\sigma, v \rangle = -m_\sigma.$
A divisor $D=\sum m_i D_i$ is principal if and only
if $m_i=\langle m,v_i \rangle $ for some $m \in \Mrm,$ so we have the
following short exact sequence:
\begin{equation}\label{ses1}
0\rightarrow \Mrm \longrightarrow \ZZ^{\Sigma(1)} \stackrel{A}{\longrightarrow}
A_{n-1}(X) \rightarrow 0.
\end{equation}
where the map $A : \ZZ^{\Sigma(1)} \to \Nrm$ is given by 
$A(\lambda_1,\ldots,\lambda_r)=\lambda_1 [D_1] + \ldots +
\lambda_r [D_r],$ with $r$ being the cardinality of $\Sigma(1).$
This exact sequence can be dualized to obtain another short
exact sequence:
\begin{equation}\label{ses2}
0\rightarrow (A_{n-1}(X))^* \longrightarrow \ZZ^{\Sigma(1)}
\stackrel{B}{\longrightarrow} \Nrm \rightarrow 0,
\end{equation}
where the map $B : \ZZ^{\Sigma(1)} \to \Nrm$ is given by 
$B(\lambda_1,\ldots,\lambda_r)=\lambda_1 v_1 + \ldots +
\lambda_r v_r,$ $r$ being the cardinality of $\Sigma(1).$

In particular, for a smooth toric variety $X,$
the group $A_1(X) \cong H_2(X, \ZZ)$ is isomorphic to the lattice 
$\LL$ of relations between the elements of $\ZZ^{\Sigma(1)}:$  
\begin{equation}
\LL:= \{\lambda=(\lambda_i) \in \ZZ^{\Sigma(1)} : 
\sum_{v_i \in \Sigma(1)} \lambda_i v_i =0 \}.
\end{equation}
We now describe another well-known construction of the $n$-dimensional
toric variety $X_{\Sigma}$ as a quotient of on an open part of 
$\CC^{\Sigma(1)}$ by a free action of the torus $\LL \otimes \CC^*.$ 
Assume that 
that the set of rays $\Sigma(1)$ has $r$ elements, 
$\Sigma (1)= \{ v_1, \ldots, v_r \}.$
For a subset $I=\{i_1, \ldots, i_d\} \subset \{1, \ldots r\}$ 
we write in short $I \in \Sigma(d)$ if the subset $\{v_{i_1}, \ldots,
v_{i_d} \}$ generates a cone in the fan $\Sigma.$ 

Let $x_1, \ldots, x_r$ be coordinates on $\CC^{\;r}$ and define 
\begin{equation}\label{sympl}
\begin{split}
\CC_{I}^{\;d} &:= \{(x_1, \ldots, x_r) : x_i \not= 0 \ \hbox{if} \ i 
\not\in I \}
 \ \hbox{for} \ I \in \Sigma(d),\\
\CC_{\Sigma}^{\;r} &:= \bigcup_{I \in \Sigma(d), 0 \leq d \leq n}\CC_I^{\;d}.
\end{split}
\end{equation}
The diagonal action of the torus $\LL \otimes \CC^*$ on 
$\CC_{\Sigma}^{\;r}$ is free and one can show that
\begin{equation}
X_{\Sigma}=\CC_{\Sigma}^{\;r} / \LL \otimes \CC^*
\end{equation}
(see, for example, \cite{guillemin}).  

For each Cartier divisor $D= \sum m_i D_i$ there exists a {\bf support
function} $\Phi_D: |\Sigma| \to \RR$, which is linear on all cones and 
such that $\Phi_D(v_i)=-m_i.$ This implies that for each $\sigma \in \Sigma$,
there is $m_\sigma \in \Mrm$ such that $\Phi_D(v)=\langle m_\sigma,v \rangle$ 
for $v \in \sigma$. To each Cartier divisor $D= \sum m_i D_i$ we associate
the {\bf supporting polyhedron} $\Delta_D$ defined by 
\begin{equation}
\begin{split}
\Delta_D :&= \{ u \in \Mrm_\RR : \langle u, v_i \rangle \geq -m_i \}\\
&= \{ u \in \Mrm_\RR : \langle u, v_i \rangle \geq \Phi_D(v_i)\  
\hbox{on} \ |\Sigma| \}.
\end{split}
\end{equation}
This integral polyhedron supports the global sections of the line bundle
${\cal O} (D),$ and we have that: 
\begin{equation}
\displaystyle{\Gamma (X_{\Sigma}, {\cal O}_D)= \bigoplus_{u \in \Delta_D 
\cap \Mrm} \CC \cdot \chi^u}.
\end{equation}
If $X_\Sigma$ is complete and a divisor $D=\sum m_i D_i$ is Cartier,
then we have that 
\begin{itemize}\item $D$ is generated by global sections if and only if
\begin{equation}\label{cond:gen}
\langle m_\sigma, v_i \rangle \geq -m_i, \ \hbox{whenever $i \not\subset 
\sigma$}.
\end{equation} \item $D$ is ample if and only if
\begin{equation}\label{cond:ample}
\langle m_\sigma, v_i \rangle > -m_i, \ \hbox{whenever $i \not\subset 
\sigma$ and $\sigma$ is $n$-dimensional}.
\end{equation}
\end{itemize}
In the first case, the support function $\Phi_D$ is said to be convex. In
the second case $\Phi_D$ is said to be
strictly convex, i.e. 
$\Phi_D(v+w) \geq \Phi_D(v) + \Phi_D(w)$, with 
equality if and only if $v$ and $w$ lie in the same cone of $\Sigma$. 
In this case, we have that $\Delta_D$ is a an $n$-dimensional polytope 
with vertices $\{ m_{\sigma}: \sigma \in \Sigma (n) \}.$ For smooth toric 
varieties, ampleness is equivalent to very ampleness. 

Now we can determine when the toric variety $X_{\Sigma}$ 
is a toric Fano variety. In the smooth case, this is equivalent
to the fact that the anticanonical bundle is ample. 
We also need to consider singular toric Fano varieties (see page
46 in \cite{CK}). A toric variety is Gorenstein if and only if
the anticanonical divisor is Cartier.
We have the following result (lemma 3.5.2. in \cite{CK}).
\begin{lemma}\label{lemma:fano}
A complete toric variety $X$ is Fano if and only if the 
anticanonical divisor $-K_X= \sum_i D_i$ is Cartier and ample.
\end{lemma}
Consider a convex polyhedron $\Delta  \subset \Mrm_{\RR}$. The 
{\bf interior point fan} $\Sigma(\Delta)$ of the convex polyhedron
$\Delta$ consists of cones obtained as positive hulls of all faces of the 
polyhedron $\Delta.$ For each 
nonempty face $F \subset \Delta,$ consider the cone
\begin{equation}\label{eq:defncone}
\check{\sigma}_F :=\{\lambda (m-m') : m \in \Delta, m' \in F, \lambda 
\geq 0 \} \subset \Mrm_{\RR}.
\end{equation}
Its dual is a cone $\sigma_F \subset \Nrm_{\RR}$. Note that the codimension
of the cone $\sigma_F$ equals the dimension of the face $F.$
The collection of cones $\{\sigma_F : F$ is a nonempty face of
 $\Delta \}$ gives the {\bf normal fan} $\Na(\Delta)$ of the polytope 
$\Delta$. This 
is a complete fan $\Sigma$, and $X(\Sigma)$ is by definition the toric 
variety associated to the polytope $\Delta.$ We will follow tradition
and denote this toric variety by $\PP_{\Delta}$. It is well-known  
\cite{GKZ2} that a fan $\Sigma \subset \Nrm_{\RR}$ is the normal 
fan of some convex polyhedron $\Delta  \subset M_{\RR}$ if 
and only if the fan $\Sigma$ is the interior point fan of the {\bf
polar polyhedron} $\Delta^*,$ where the polar polyhedron is defined by
\begin{equation}
\Delta^*:= \{ v \in \Nrm_{\RR} : \langle m,v \rangle \ \geq -1 \ 
\hbox{for all} \ m \in \Delta \} \subset \Mrm_{\RR}.
\end{equation}

\begin{definition}\label{reflexive}\cite{Bat1}
A convex $n$-dimensional polyhedron $\Delta \subset \Mrm_{\RR}$ with
all the vertices in the lattice $\Mrm \subset \Mrm_{\RR}$ and containing the 
origin in its interior is called {\bf reflexive} 
if all the vertices of the polar polyhedron $\Delta^*$ belong
to the dual lattice $\Nrm \subset \Nrm_{\RR}.$ 
\end{definition}
If $\Delta \subset \Mrm_{\RR}$ is a 
reflexive polyhedron, then  $\Delta^* \subset \Nrm_{\RR}$ is again a
reflexive polyhedron and $(\Delta^*)^* =\Delta$. This involution
on the set of reflexive polyhedra plays a crucial r{\^ o}le 
in Batyrev's construction of mirror symmetric Calabi--Yau manifolds.
The toric variety $\PP_{\Delta}$ (resp. $\PP_{\Delta^*}$) is Fano if and
only if the polyhedron $\Delta$ (resp. $\Delta^*$) is reflexive. In 
particular, $\PP_\Delta$ and $\PP_{\Delta^*}$ have at worst canonical 
singularities. Special cases of reflexive polyhedra are the
unimodular polyhedra (all the $(n-1)$-dimensional faces are 
simplicial and of unit volume), which are customarily called 
Fano polyhedra in the literature (see for example \cite{BatFano}).
In this terminology, the interior fan of a Fano polyhedron
gives rise to a smooth Fano toric variety. 

Given a reflexive polyhedron, the toric variety $\PP_\Delta$ may be too
singular for the purposes of mirror symmetry. Because of this we 
will consider fans that refine the normal fan $\Na(\Delta) \subset \Nrm_\RR.$
Namely, we define a {\bf projective subdivision} of the 
fan $\Na(\Delta)$ to be a projective simplicial fan that refines the fan
$\Na(\Delta)$ and has its rays included in $\Delta^* \cap N.$
Such a projective subdivision determines in fact a regular
triangulation of $\Delta^* \cap \Nrm$ (the notion of a regular
triangulation will be defined in the next section). The toric
variety $X_\Sigma$ determined by a projective subdivision $\Sigma$ of
$\Na(\Delta)$ can have only orbifold singularities, and $\Delta$ is 
the polytope associated
to $-K_X.$ The important fact about the projective subdivisions
is that the induced map $f:X \to \PP_\Delta$ is crepant,
that is $f^*(K_{\PP_\Delta})= K_X.$ Also, the general member 
of the of the linear system $\vert -K_{\PP_\Delta} \vert$ is
a Calabi--Yau variety, while the general member of the linear
system $\vert -K_X \vert$ is a Calabi--Yau orbifold. For a complete
review of the proofs of these facts (due to Batyrev in this context),
see section 4.1. of \cite{CK}.
Our main interest will be the case when the Calabi--Yau
variety is smooth, so it does not intersect the orbifold locus
of the toric variety $X.$
Batyrev's mirror construction for hypersurfaces states that
the mirror of the Calabi--Yau hypersurface obtained in this way
is the hypersurface obtained by applying the same procedure
to the polytope $\Delta^*.$ 

The first mirror construction was proposed by Greene and Plesser 
\cite{GrPl} using the symmetries in the conformal field theory 
associated to the two geometrically distinct mirror manifolds. Based on 
Greene--Plesser orbifolding mirror construction, Batyrev 
proposed the general toric construction of mirror symmetric
Calabi--Yau manifolds. We review here the mirror symmetry construction 
in the case of Calabi--Yau complete intersections.
Suppose that $X=X_{\Sigma(\Delta^*)}=\PP_{\Delta}$
is an $n$-dimensional toric variety corresponding to a reflexive
polytope $\Delta$ and determined by the normal fan $\Sigma=\Na(\Delta)=
\Sigma(\Delta^*).$ Assume that the
divisors $E_1, \ldots, E_k$ determine a partition of the set 
$\Sigma(1)$  into a disjoint union of $k$ subsets $I_j, j=1, 
\ldots, k,$ such that 
\begin{equation}
E_j= \sum_{i \in I_j} D_i, \ j=1,\ldots, k.
\end{equation}
such that 
\begin{equation}
-K_X= \sum_{j=1}^k E_j.
\end{equation}
The decomposition $\Sigma(1)= I_1 \cup \ldots \cup I_k$
is called a {\bf nef--partition} if, for each j, $E_j$ is a 
Cartier divisor spanned by its global sections. If $\Delta_j$
is the polytope in $\Mrm$ corresponding to $E_j,$ then being a nef--partition
implies that $\Delta= \Delta_1+\ldots+\Delta_k$ (Minkowski sum).
The complete
intersection $W \subset X$ of generic sections of $\Oa(E_1)$, $\ldots,$ 
$\Oa(E_k)$ is a Calabi--Yau variety. If necessary, one can consider
a projective subdivision $\Sigma^\prime$ of $\Sigma(\Delta^*).$ The
proper transform of $W$ under the induced map 
$X_{\Sigma^\prime} \to X_{\Sigma}$ is again a Calabi--Yau
complete intersection.

The Batyrev--Borisov construction \cite{BB1} states that in order
to obtain the mirror family we should consider the polytopes
\begin{equation}
\nabla_j:= \hbox{Conv}(\{0\} \cup I_j) \subset \Nrm_{\RR}.
\end{equation}
Borisov \cite{Bor} proved that the Minkowski sum
\begin{equation}
\nabla=\nabla_1 + \ldots + \nabla_k
\end{equation}
determines a nef--partition of the reflexive polytope $\nabla,$
and that 
\begin{equation}
\begin{split}
\Delta^*= \hbox{Conv}(\nabla_1 \cup \ldots \cup \nabla_k), \
\nabla^*= \hbox{Conv}(\Delta_1 \cup \ldots \cup \Delta_k).
\end{split}
\end{equation}
This dual nef-partition defines the mirror family $M$ of the complete
intersection $W.$ This means that $M$ is the complete intersection
of generic sections of $\Oa (E_j^*),$ supported on the polytopes 
$\nabla_j$ $(j=1, \ldots, k).$ More explicitly, 
the mirror manifolds $M$ are birationally isomorphic to affine complete
intersections in $(\CC^*)^n=\hbox{Spec}\;\CC[t_1^{\pm 1},\ldots, 
t_n^{\pm 1}]$ defined by $k$ equations 
\begin{equation}\label{batmir}
1= \sum_{v_i \in I_j} a_i {\bf t}^{v_i}, \ j \in \{ 1, \ldots, k \},
\end{equation}
where $a_i$ are generic complex numbers and
${\bf t}^{v_i}, v_i \in \Sigma(1),$ are Laurent monomials in the 
variables $t_1, \ldots, t_n,$ $n=\hbox{rk} \Nrm.$

The following notion was introduced by Batyrev and Borisov \cite{BB1}
in their study of Calabi--Yau complete intersections 
in toric Fano varieties.

\begin{definition}\label{goren}
A cone $\sigma$ in $\RR^n$ is called a {\bf Gorenstein cone}
if it is generated by a finite set $\{v_1,\ldots,v_N\}\subset\ZZ^n$
\begin{equation}\label{eq:lambda}
\sigma=\RR_{\geq 0}v_1\,+\ldots+\,\RR_{\geq 0}v_N\,,
\end{equation}
such that there exists $v^{\vee}_0 \in \RR^{n \vee}$ with the property 
that \begin{equation}\label{gorencond}
\langle v^{\vee}_0,v \rangle =1, \ \hbox{for every generator} 
\ v \ \hbox{of} \  \sigma.
\end{equation}
It is called a {\bf reflexive Gorenstein cone} if both $\sigma$
and its dual $\check{\sigma}$ are Gorenstein cones,
i.e. there  exists a vector $v_0 \in \ZZ^n$ and a set
$\{v_1^\vee,\ldots, v_{N^\prime}^\vee \}\subset\ZZ^{n\vee}$
of generators of $\check{\sigma}$ such that
\begin{equation}
\langle v^{\vee}_0,v_0 \rangle =1, \ \hbox{for every generator} \ 
v^\vee \ \hbox{of}  \ \check{\sigma}.
\end{equation}
The vectors $v^{\vee}_0$ and $v_0$ are uniquely determined by
$\sigma,$ and the integer $\langle v^{\vee}_0, v_0 \rangle $ is called 
{\bf the index} of
$\sigma$ (or $\check{\sigma}$).
\end{definition}
We are now able to present the toric situation that will 
be the main focus of this work.

\begin{remark}\label{rem:toric}
The duality of nef--partitions can be described using
reflexive Gorenstein cones. Consider the nef--partition
$\Sigma(1)= I_1 \cup \ldots \cup I_k,$ and the Minkowski sum
described above $\nabla= \nabla_1 + \ldots + \nabla_k.$
By extending the lattice $\Nrm$ to  
$\overline{\Nrm}:= \Nrm \times \ZZ^k,$ we can define the cone 
$\sigma \subset \overline{\Nrm}$ by 
\begin{equation}\label{eq:goren}
\sigma:= \{(\sum_{j=1}^k \lambda_j v_j,\lambda_1, \ldots, \lambda_k) : 
v_j \in \nabla_j, \lambda_j \geq 0\}.
\end{equation}
It turns out \cite{BB1} that $\sigma$ is a reflexive Gorenstein cone
of index $k$ and its dual $\check{\sigma}$ (a reflexive Gorenstein cone of
index $k$ itself) is constructed by applying (\ref{eq:goren}) to the 
dual nef--partition and the Minkowski sum
$\Delta= \Delta_1 + \ldots + \Delta_k.$ 

In fact, the reflexive Gorenstein cone $\sigma$ is the support of 
a fan which determines the toric structure of the underlying 
space of a vector bundle $\EE \to X_{\Sigma}=\PP_{\Delta}.$
For any vector $v \in \Sigma(1)$ such that $v \in I_{i_0},$
$1 \leq i_0 \leq k,$ define the vector $\overline{v}:= (v,\delta_1^{i_0}, 
\ldots, \delta_k^{i_0}),$ where $\delta_i^j$ is the usual Kronecker 
delta function. The maximal cones of the fan $\overline{\Sigma}$ 
defining  $\EE$ are obtained by adjoining to the maximal cones of 
$\Sigma$ the vectors
$(0,e_1), \ldots, (0,e_k),$ where $e_i$ is the usual $i^{\rm th}$
standard basis element in $\ZZ^k.$  In the terminology of 
Batyrev and Borisov \cite{BB1}, this is an instance where the cones 
$\sigma$ and $\check{\sigma}$ are completely split, i.e. 
\begin{equation}\label{eq:splitvb}
\begin{split}
\EE &= \Oa (-E_1) \oplus \dots \oplus \Oa(-E_k), \\
\EE^\vee& = \Oa (E_1) \oplus \dots \oplus \Oa(E_k).
\end{split}
\end{equation}

The toric variety $X_\Sigma=\PP_\Delta$ (the zero section 
of the bundle $\EE=X_{\overline{\Sigma}}$) corresponds to
the closed toric orbit in $\EE$ determined by the cone in $\overline{\Sigma}$
generated by the vectors $(0,e_1), \ldots, (0,e_k).$
The Calabi--Yau complete intersection $W$ is obtained by intersecting 
the zero loci of 
$k$ generic sections of $\EE^\vee$ corresponding to $E_1,\dots, E_k,$ 
respectively (see sections 8 and 9 in \cite{Stienstra} for a 
detailed explanation of this construction). For example, in the 
hypersurface case $k=1,$ the toric
variety $\EE$ is the underlying variety of the canonical bundle of 
$\PP_\Delta$ and the corresponding Calabi--Yau variety is the zero set 
of a section of the anticanonical bundle $\EE^\vee.$ In fact, 
the different physical phases in the string theoretic investigations of 
Witten \cite{Witten1} and Aspinwall, Greene and Morrison \cite{AGM1}
do not correspond directly to different fans with the same support
as the fan $\Sigma=\Sigma(\Delta^*)$ defining $\PP_{\Delta},$ but
rather to different fans supported on $\sigma$ given by different regular 
triangulations of the reflexive Gorenstein cone. 
The next section is devoted to a brief review of regular triangulations
and polyhedral subdivisions of $\cA$--sets. 
\end{remark}

\begin{conv} In the rest of this work,
we will make a notational distinction between the generators of 
Gorenstein cones and $\cA$--sets (denoted by $\ov, \overline{w},$ etc.) and 
the generators of complete fans defining 
compact toric varieties (denoted by $v, w,$ etc.). 
\end{conv}

\subsection{The Secondary Fan and the Secondary Polytope}
\label{chap:triang}

A finite subset ${\cal A}= \{ \ov_1, \ldots, \ov_N \}$ of 
the integral lattice $\ZZ^n$ is called an {\bf $\cA$--set} if
it generates the lattice $\ZZ^n$
and there exists an element $\ov^{\vee}_0 \in \ZZ^{n\vee}$
such that 
\begin{equation}\label{hyperplane1}
\langle  \ov^{\vee}_0, \ov_i \rangle =1, \ (i=1,\ldots,N).
\end{equation}
Note that any Gorenstein cone determines an $\cA$--set. 
Denote by $\LL \subset \ZZ^N$ the lattice of relations among
the elements of $\cA,$
\begin{equation}
\LL:= \{ l=(l_1,\ldots,l_N) : l_1 \ov_1+ \ldots +l_N \ov_N=0 \}.
\end{equation}

A {\bf triangulation} of the set ${\cal A}= \{\ov_1, \ldots, \ov_N \}$ in 
$\RR^n,$ is a triangulation of the 
convex hull $\ \hbox{Conv}({\cal A})\ $ in $\RR^n$ such that all
vertices are among the the points $\ov_1, \ldots, \ov_N.$ Note, that
condition (\ref{hyperplane1}) shows that Conv$(\cA)$ has 
dimension $n-1.$ When there is no danger of confusion, given a subset 
$I \subset \{1,\ldots,N\},$ such that the vectors $\{v_i : i \in I \}$
are linearly independent, the simplex 
Conv$\{\ov_i : i \in I \}$ will simply be denoted by $I,$ and we will
use the notation $I^*$ to denote the set $\{ 1, \ldots, N \} \setminus I.$ 
Since we have to deal not only with simplices, but also with
polytopes, it will be convenient for us to use the 
notation $(Q,A),$ to denote a marked $(n-1)$-dimensional 
polytope with $Q=\hbox{Conv}(A),$ and $A \subset \cA.$ 

A {\bf regular triangulation} of $\Ta$ of $\cal A$ is a triangulation 
such that there exists a function on Conv($\cal A$) which is strictly 
convex with respect to the the triangulation $\Ta.$ There are other 
more geometrical ways of defining regular triangulations, see \cite{GKZ2}.
We describe a method for parameterizing the set of regular triangulations.
Let $w_1, \ldots, w_N \in \LL^{\vee}_{\RR}$ be the images of the standard 
basis vectors of $\RR^{N\vee}$ under the surjection 
$\RR^{N\vee} \to \LL^{\vee}_{\RR}$ dual to the inclusion 
$\LL \hookrightarrow \ZZ^N.$ This collection of vectors is called the 
{\bf Gale transform} of $\cA.$ Let $\cal B$ (resp. $\cal D$) be the collection
of those subsets $I$ and $J$ of $\{1, \ldots, N\}$ of cardinality $N-n$
(resp. $N-n-1$) for which the vectors $w_i$ ($i \in I$) are linearly
independent. For $I \in \cal B$ and $J \in \cal D$ we introduce
the cones 
\begin{eqnarray}\label{cone99}
{\cal C}_I &:=& \{ \sum_{i \in I} t_i w_i: t_i \in \RR_{\geq 0} \}, \\
\label{hyperplane2}
H_J &:=& \{ \sum_{i \in J} t_i w_i: t_i \in \RR \}.
\end{eqnarray}
The following properties hold: 

\begin{proposition} (cf. \cite{Stienstra})\label{prop:cone}
For every $w \in \LL^{\vee}_{\RR} \setminus 
\cup_{J \in {\cal D}} H_J$ the set 
\begin{equation}
\Ta^n:= \{ I : I^* \in {\cal B} \ \hbox{and} \  w \in {\cal C}_{I^*} \}
\end{equation}
is the set of maximal simplices of a regular triangulation $\Ta.$
\end{proposition}

The Gale transform can be given a more explicit description in
coordinates. Consider the $N \times n$ matrix $A=(v_{ij})$ whose
rows are the vectors $v_1, \ldots, v_N.$ We have that $\LL=\ker A 
\subset \ZZ^N,$ and let $B$ the $(N-n) \times N$ integer matrix
describing this kernel. The rows of $B$ are the vectors generating
the lattice $\LL$ in $\ZZ^N.$ If we identify the lattice 
$\LL^{\vee}_{\RR}$ with $\RR^{N-n \vee},$ the columns of $B$ 
give the vectors $w_1, \ldots, w_N.$ The following 
property holds:

\begin{proposition}\label{simplexvol}
If the cone ${\cal C}_I$ is generated by the 
vectors $w_i \in \LL^{\vee}_\RR,$ $i \in I, |I|=N-n,$
then the volume of the simplex 
generated by $w_i, i \in I,$ in $\RR^{N-n \vee} \cong \LL^{\vee}_\RR$
is equal to the volume of
the simplex generated by the vectors $v_j\; , j \in I^*,$
in $\RR^n.$
\end{proposition}
In particular, for a unimodular triangulation $\Ta$  the 
cones ${\cal C}_{I^*}$ have volume $1$ for all $I \in \Ta^n.$

For a regular triangulation $\Ta$ of $\Delta,$ we define
\begin{equation}\label{secon}
{\cal C}_{\Ta}:= \bigcap_{I \in \Ta^n} {\cal C}_{I^*}.
\end{equation}
Every $w \in {\cal C}_{\Ta} \setminus \cup_{J \in {\cal D}} H_J$ 
leads to the same triangulation $\Ta$. The collection of cones
${\cal C}_{\Ta} \subset \LL^{\vee}_{\RR},$ for all regular triangulations 
$\Ta$  of $\Delta,$ is a complete fan in $\LL^{\vee}_{\RR}$
called the {\bf secondary fan}. One of the important results of
Gelfand, Kapranov and Zelevinsky (cf. \cite{GKZ2}) is that
this fan is the normal fan of a polytope in $\LL_\RR$ called the
{\bf secondary polytope}. An interesting case is the case when 
all the cones in the secondary fan are generated by the vectors
in Gale transform (this always happens for $N-n \leq 2$). In that
case, we can always choose coordinates so that the cone
in the secondary fan corresponding to a unimodular triangulation 
is generated by the standard basis elements in $\RR^{N-n\vee}.$

Consider two regular triangulations $\Ta_1$ and $\Ta_2$ which correspond
to adjacent cones in the secondary fan. Equivalently, the vertices
of the secondary polytope corresponding to $\Ta_1$ and $\Ta_2$ are
joined by an edge. According to the general theory of regular 
triangulations \cite{GKZ2} the two triangulations differ by a 
modification (``perestroika'') along a circuit. The rest of
this section is a review of concepts discussed in chapter 7 of the
book \cite{GKZ2}.

A {\bf circuit} in $\cal A$ is a minimal independent subset 
$I \subset \{1,\ldots, N \}.$ In particular any circuit determines a relation 
$l=(l_1, \ldots, l_N) \in \LL,$  which has the form
\begin{equation}\label{circuit}
\sum_{i \in I_+ } l_i \ov_i + \sum_{i \in I_-} l_i \ov_i =0,
\end{equation}
where the two subsets $I_+:= \{ i: l_i > 0 \}$ and  $I_-:= \{ i: l_i < 0 \}$
are uniquely defined by the circuit up to replacing $I_+$ by $I_-$. The 
minimizing condition in the definition of a circuit implies that the 
sets Conv($I_+$) and Conv($I_-$) intersect in their common interior point. 
The polyhedron Conv($I$) admits exactly two triangulations $\Ta_+(I)$
and $\Ta_-(I)$ defined by the simplices Conv($I \setminus i$), for 
$i \in I_+$, respectively $i \in I_-$. Given a triangulation $\Ta$ of 
$\cA$ such that the induced triangulation on $I$ is $\Ta_+(I),$ one
can obtain a new triangulation $s_I(\Ta)$ of $\cA$ just by replacing 
all the simplices in $\Ta_+(I)$ by those in $\Ta_-(I)$. We review
the following definition.

\begin{definition}\label{def:support}
Let $\Ta$ be a triangulation of $\cA$ and let $I \subset\{1,\ldots, N \}$ 
be a circuit. We say that {\bf $\Ta$ is supported on $I$}, if the 
following conditions hold:
\begin{enumerate}
\item There are no vertices of $\Ta$ inside Conv$(I)$ except for the
elements of $I$ itself. 
\item The polytope Conv$(I)$ is a union of faces of the simplices 
of $\Ta.$
\item Let Conv$(J)$ and Conv$(J')$ be two simplices of maximal 
dimension of one of the two possible triangulations of 
Conv$(I).$ Then, for every subset $F \subset \{1,\ldots, N \} \setminus
I,$ the simplex Conv$(J \cup F)$ appears in $\Ta$ if and only if
the simplex Conv$(J' \cup F)$ appears. 
\end{enumerate}
\end{definition}
It turns out that the following characterization is true:
\begin{proposition}\label{prop:perestr}
Let $\Ta_1$ and $\Ta_2$ be two regular triangulations of $\cA.$ Then
the corresponding cones in the pointed secondary fan are adjacent 
(i.e. the corresponding vertices in the secondary polytope are
joined by an edge) if and only if there exists a circuit 
$I \subset \cA$ such that both $\Ta_1$ and $\Ta_2$ are supported on 
$I$ and $s_I(\Ta_1)=\Ta_2.$
\end{proposition}
The proof of this fact 
can be found in \cite{GKZ2}. Let ${\cal C}_1$ and ${\cal C}_2$ denote 
two adjacent cones in the 
secondary fan corresponding to the regular triangulations $\Ta_1$ and 
$\Ta_2$, respectively, with $s_I(\Ta_1)=\Ta_2$ and 
let $l=(l_1, \ldots, l_N) \in \LL_\RR$ be the induced relation
(\ref{circuit}). The common face of the cones ${\cal C}_1$ and 
${\cal C}_2$ is contained in a hyperplane in $\LL^{\vee}_{\RR}.$ Recall 
that we denoted by  
$w_1, \ldots, w_N \in \LL^{\vee}_{\RR}$ the elements of the Gale
transform of $\cA.$ We have the following straightforward 
lemma.

\begin{lemma} \label{lem} Assume that the hyperplane $H \subset 
\LL^{\vee}_{\RR}$ containing the common face of two adjacent cones
$\Ca_1$ and $\Ca_2$ in the secondary fan is given by some
$h \in \LL_\RR,$ $H:= \{ w \in \LL^{\vee}_{\RR} : \langle 
w,h \rangle =0 \}.$ Consider the subsets $I_+$ and $I_-$ in 
$\{ 1,\ldots,N \}$ defined by $I_+ := \{ i : \langle w_i,h \rangle
> 0 \}$ and $I_- := \{ i : \langle w_i,h \rangle < 0 \}.$ Then:
\begin{enumerate}
\item The vectors $w_i$ for $i \in \{1,\ldots, N\} 
\setminus (I_+ \cup I_-)$
generate the common face of the cones $\Ca_1$ and $\Ca_2,$ and
the set $I =I_+ \cup I_-$ is the circuit defining the transformation
$s_I$ between the corresponding regular triangulations 
$\Ta_1$ and $\Ta_2.$ 
\item \begin{equation}
\sum_{i \in I_+ } \langle w_i,h \rangle \ov_i + \sum_{i \in I_-} 
\langle w_i,h \rangle \ov_i=0.
\end{equation}
\end{enumerate}
\end{lemma}

Note that the hyperplane $H$ is one of the
hyperplanes $H_J,$ defined by (\ref{hyperplane2}) for some 
subset $J \subset \{1, \ldots, N \}$ with $|J|= N-n-1.$ We can
always assume then that
$h$ is chosen to be in integral in $\LL_{\RR},$ which means that
$h\in \LL \subset \LL_{\RR}.$ Moreover, if  we choose $h$ to be
saturated in $\LL$ (minimal vector in $\LL$ in that direction),
then $h$ is always part of an integral basis of $\LL$
(for any saturated vector in $\ZZ^{N-n},$ we 
can find a unimodular integer matrix having the coordinates of the vector 
on one row).

In fact, Gelfand, Kapranov and Zelevinsky gave an interpretation 
for all the faces of the secondary polytope in terms of polyhedral
subdivisions. Given an $(n-1)$-dimensional marked polytope $(Q,A),$ with 
$Q=\hbox{Conv}(A),$ $A \subset \cA,$ we have the following.

\begin{definition}\label{def:polysub}
A {\bf polyhedral subdivision} of $(Q,A)$ is a family 
of $(n-1)$-dimensional marked polytopes  
$(Q_j,A_j), j=1, \ldots, l,$ $A_j \subset A,$ such that
\begin{enumerate} 
\item Any intersection $Q_j \cap Q_k$ is a face (possible empty) of both 
$Q_j$ and $Q_k,$ and 
\begin{equation}
A_j \cap (Q_j \cap Q_k) = 
A_k \cap (Q_j \cap Q_k).
\end{equation}
\item The union of all $Q_j$ coincides with $Q.$
\end{enumerate}
\end{definition}
For two polyhedral subdivisions $\Sa=(Q_j,A_j)$ and $\Sa'=(Q'_k,A'_k),$ we 
shall say that $\Sa$ {\bf refines} $\Sa'$ if, for each $k,$ the collection
of $(Q_j,A_j)$ such that $ Q_j \subset Q_k,$ forms
a subdivision of $Q_k.$ The set of polyhedral subdivisions is then 
naturally a partially ordered set (poset) with respect to refinement. 
Triangulations are minimal elements of this poset. The maximal element 
is the subdivision $(\hbox{Conv}(\cA),\tilde{\cA}),$  
where $\tilde{\cA} \subset \cA$ is the minimal subset of vertices 
such that $\hbox{Conv}(\tilde{\cA})=\hbox{Conv}(\cA).$ 

A polyhedral subdivision $\Sa$ is called {\bf regular} 
if there exists a function on Conv($\cal A$) which is strictly 
convex with respect to the subdivision $\Sa.$ In the afore--mentioned 
book \cite{GKZ2}, it is shown that the poset of regular
polyhedral subdivisions is naturally isomorphic to the 
poset of faces (cones) of the secondary polytope (fan) determined by 
$\cA.$ In particular, for any face $F$ of the secondary polytope
there exists a well defined polyhedral subdivision $\Sa(F),$
and this operation is compatible with the two partial order relations. 
We will not explain here the explicit construction of the
polyhedral subdivision corresponding to a face (cone) in the 
secondary polytope (fan) in full generality. For later use, we now present
this construction for the edges of the secondary polytope.

\begin{definition}\label{def:sep}
Suppose that the regular triangulations $\Ta_1$ and $\Ta_2$ 
of $\cA$ are obtained from each other by a modification 
along a circuit $I, s_I (\Ta_1)= \Ta_2.$ We say that a 
subset $J \subset \cA \setminus I,$ is {\bf separating}
for $\Ta_1$ and $\Ta_2$ if, for some $i \in I,$ the set
$(I \setminus i) \cup J$ is the set of vertices of a 
simplex (of maximal dimension) of $\Ta_1.$
\end{definition}

Note that, because of definition \ref{def:support},
a separating set $J$ has the property that 
$(I \setminus i) \cup J,$ for $i \in I_+,$ are simplices 
in $\Ta_1$ and $(I \setminus i) \cup J,$ for $i \in I_-,$ are
simplices in $\Ta_2.$

\begin{proposition} \label{prop:edge} 
\sloppy{
The polyhedral subdivision 
$\Sa(\Ta_1, \Ta_2)$
corresponding to the edge of the secondary polytope
joining the vertices determined by $\Ta_1$ and $\Ta_2$ consists
of the simplices $(Conv(K),K)$ which $\Ta_1$ and $\Ta_2$
have in common and the polyhedra $(Conv(I \cup J), I\cup J)$
for all separating subsets $J \subset \cA \setminus I.$}
\end{proposition}
It is obvious that all the presented facts about triangulations and 
polyhedral subdivisions can be rephrased in the more often used language 
of fan refinements in toric geometry. The set $\cA$ determines 
a reflexive Gorenstein cone (definition \ref{goren}) which is the support 
of different (simplicial) fans determined by (triangulations) polyhedral
subdivisions. The explicit constructions described in the last
part of this section will play a distinguished r\^{o}le in 
our interpretation of Kontsevich's conjecture.

\subsection{Calabi--Yau Complete Intersections in Toric \\ Varieties}
\label{chap:cicy}

We now return to the specific toric situation described 
in remark \ref{rem:toric}. We start by presenting
the notions explained there in the language of $\cA$--sets and 
triangulations. Consider a reflexive polytope $\Delta \subset \Mrm_\RR,$ 
and assume that the fan $\Sigma \subset \Nrm_\RR$ defines a projective 
subdivision of the normal fan $\Na(\Delta).$  If $X=X_\Sigma$ is the 
corresponding toric Fano variety, let $W=E_1 \cap \ldots E_k$ be a 
Calabi--Yau complete intersection in $X,$ where the Cartier 
divisors $E_1, \ldots, E_k$ determine a nef--partition of the set 
$\Sigma(1)$  given by a disjoint union of $k$ subsets $I_j, j=1, 
\ldots, k.$ We have that  
\begin{equation}\label{sume}
E_j= \sum_{i \in I_j} D_i, \ j=1,\ldots, k.
\end{equation}
The Calabi--Yau condition implies that 
\begin{equation}
-K_X= \sum_{j=1}^k E_j.
\end{equation}
Consider the lattice $\overline{\Nrm}:=\Nrm \oplus \ZZ^k \cong \ZZ^n,$ 
$n=\hbox{rk} \Nrm \; +\; k.$ For each element $v \in \Sigma(1),$ such that
$v \in I_{i_0},$ $1 \leq i_0 \leq k,$ define the vector 
$\overline{v} \in \overline{\Nrm}$ given by $\overline{v}:= (v,\delta_1^{i_0}, 
\ldots, \delta_k^{i_0}),$ where $\delta_i^j$ is the usual Kronecker 
delta function. The $\cal A$--set corresponding to the Calabi--Yau
complete intersection $W \subset X$ is defined to be 
\begin{equation}\label{newA}
{\cal A}:= \{(0,e_1),\ldots,(0,e_k)\}\cup \{\overline{v}: v \in \Sigma(1)\},
\end{equation}
where $e_i$ is the usual $i^{\rm th}$ standard basis vector in $\ZZ^k.$
For the ease of notation, we denote the elements of the set ${\cal A}$ 
by $\{ \ov_1, \ldots , \ov_N\},$ such that
$\ov_j=(0,e_j),$ for $j=1,\ldots,k.$ It follows that
\begin{equation}
\ov^{\vee}_0= (\overset{1}1,\ldots,\overset{k}1,0,\ldots,0),
\end{equation}
where the first $k$ components are units. 

We now recall the combinatorial definition of the Stanley--Reisner ring.

\begin{definition}\label{reisner}
For any integral polytope $\Xi=\hbox{Conv}\{u_1,\ldots,u_d \}$
$ \subset \RR^n,$ with
$\{u_1,\ldots,$ $u_d \}$ $ \subset \ZZ^n,$
and any triangulation $\Ta$ of the polytope 
$\Xi,$ the 
{\bf Stanley--Reisner ring} of the 
set $\{u_1, \ldots, u_d \}$ is defined by
\begin{equation}
\Ra_{\;\Ta}:= \ZZ[\lambda_1, \ldots, \lambda_d]/ {\cal I},
\end{equation}
where the ideal $\cal I$ is the sum of the ideals 
${\cal I}_{lin}$ generated by the linear forms
\begin{equation}\label{linear}
\sum_{i=1}^d \langle u^{\vee}, u_i \rangle \lambda_i \ \ \hbox{for all} \ 
u^{\vee} \in \ZZ^{n\vee},
\end{equation}
and the ideal ${\cal I}_{mon}$ generated by the monomials
\begin{equation}\label{polyn}
\lambda_{i_1} \cdot \ldots \cdot \lambda_{i_s} \ \ \hbox{for} \
\hbox{\rm Conv}\{u_{i_1}, \ldots, u_{i_s}\} \ 
\hbox{not a simplex in} \ \Ta.
\end{equation}
\end{definition}

In the case when $X$ is smooth (i.e. the fan $\Sigma$ is smooth), 
it is a basic fact in the theory of 
toric varieties that
\begin{equation}
\Ra_{\;\Sigma}\cong A^* (X) \cong H^*(X,\ZZ).
\end{equation}

In the case when $X$ is singular with orbifold singularities 
(i.e. the fan $\Sigma$ is simplicial) it is still true that
\begin{equation}
\Ra_{\;\Sigma} \otimes \QQ \cong A^*(X)_\QQ \cong H^*(X,\QQ).
\end{equation}
Under these isomorphisms, the class of the toric divisor $D_j$
in $A^1(X)_\QQ \cong H^2(X,\QQ)$ corresponds to 
$\lambda_j \in \Ra_{\;\Sigma} \otimes \QQ,$ 
$k+1 \leq j \leq N.$

The triangulation determined by the simplicial fan $\Sigma$
induces a triangulation $\Ta$ of the set $\cal A.$ 
The maximal simplices of the triangulation $\Ta$ are obtained by 
adjoining the elements $\ov_1,\ldots,\ov_k$ to each maximal simplex 
(cone) of $\Sigma.$ Note that a smooth fan $\Sigma$ determines
a unimodular triangulation $\Ta.$
The following commutative diagram summarizes 
the relationship between $\cal A$ and $\Sigma.$
\begin{equation}
\begin{array}{ccccccccc}
0 & \longrightarrow & \LL \otimes \QQ & \longrightarrow & 
\QQ^N & \stackrel{\overline{B}}{\longrightarrow} & \overline{\Nrm}_\QQ& 
\longrightarrow & 0\\
\ & \ &  \downarrow & \ & \downarrow & \ & \downarrow & \ & \ \\
0 & \longrightarrow & H_2(X,\QQ)  & \longrightarrow & \QQ^{N-k}&
\stackrel{B}{\longrightarrow} & \Nrm_\QQ & \longrightarrow & 0  
\end{array}
\end{equation}
If $X$ is smooth, we have an analogous commutative diagram 
with coefficients in $\ZZ.$

The next proposition contains the essential facts about 
this diagram.
\begin{proposition}(page 97 in \cite{CK})\label{katzcox}
\
\begin{enumerate}
\item Each non--torsion class $\theta \in H_2(X,\ZZ)$ determines a 
relationship
\begin{equation}
\sum_{j=1}^k (E_j \cdot \theta)\ov_j=
\sum_{j=k+1}^N (D_j \cdot \theta)\ov_j.
\end{equation}
\item The leftmost vertical line in the above commutative diagram is 
an isomorphism, that is  
\begin{equation}
\LL \otimes \QQ \cong H_2(X,\QQ) \cong \Lambda \otimes \QQ,
\end{equation}
where $\Lambda$ is the lattice of relations corresponding to the fan
$\Sigma.$
\end{enumerate}
\end{proposition}
In conjunction with  this proposition, the construction of the 
triangulation $\Ta$ shows that the
elements $\lambda_1,\ldots, \lambda_k $ of $\Ra_{\;\Ta}$ are obtained
as linear combinations of the elements $\lambda_{k+1},\ldots, \lambda_N.$
The monomials of the type (\ref{polyn}) remain unchanged when
one replaces the ring $\Ra_\Sigma$ with the ring $\Ra_{\;\Ta}.$ The
following properties hold:

\begin{proposition}\label{toricoh}
\
\begin{enumerate}
\item  For a smooth toric variety $X_{\Sigma}$ we have that 
\begin{equation}
\Ra_{\;\Ta} \cong \Ra_\Sigma \cong A^*(X) \cong H^*(X,\ZZ).
\end{equation}
In general, for a simplicial fan $\Sigma,$ we have that 
\begin{equation}\label{iso:chow}
\Ra_{\;\Ta}\otimes \QQ \cong \Ra_\Sigma \otimes \QQ \cong 
A^*(X)_\QQ \cong H^*(X,\QQ).
\end{equation}
\item  The classes of the divisors $E_j$ correspond to the elements 
$-\lambda_j \in \Ra_{\;\Ta},$ $1 \leq j \leq k$ under the isomorphism 
(\ref{iso:chow}).
\item The ring $\Ha \otimes \QQ,$ with
$\Ha:= \Ra_{\;\Ta}/\hbox{Ann}(\lambda_1 \cdot \ldots \cdot
\lambda_k),$ is isomorphic to a quotient ring of the toric part of 
the cohomology ring $H^*(W,\QQ)$ 
of the Calabi--Yau complete intersection W.\footnote{I want to thank 
A.~Mavlyutov for an observation regarding a subtlety in this statement.}
\item When $W$ is smooth (so, it does not intersect the 
singular locus of $X_{\Sigma}$), the Todd class of $W,$
$\Td_W \in H^*(W,\QQ) \cong \Ha \otimes \QQ,$ is given by
\begin{equation}
\Td_W =\Bigl( \;\prod_{j=1}^k \frac{-\lambda_j} {1- e^{\lambda_j}}\;
\Bigr)^{-1}
\prod_{j=k+1}^N \frac{\lambda_j}{1- e^{-\lambda_j}}.
\end{equation}
\end{enumerate}
\end{proposition}
In general, it may happen that there exist 
non--toric even cohomology classes in 
$H^*(W,\ZZ).$ In that situation, $W$ may admit K\"ahler deformations
which do not come from K\"ahler deformations of the ambient
toric variety $X.$ Our analysis refers only to toric deformations.

\begin{proof} The first part contains standard facts in the theory 
of toric varieties. The only point to be made is that, since the vectors
$v_1, \ldots, v_k$ belong to all maximal simplices in the triangulation
$\Ta,$ no monomial in $\Ia_{mon}$ corresponding to $\Ta$ will
contain any of the elements $\lambda_1, \ldots, \lambda_k.$ The second 
part is a direct consequence of the way $\cA$ was obtained from 
$\Sigma.$ It implies that the following relation holds in 
$\Ra_\Ta$
\begin{equation}
\vert I_j \vert \; \lambda_j + \sum_{i\in I_j} \lambda_i =0, \ \hbox{for
$1 \leq j \leq k$}.
\end{equation}
Formula (\ref{sume}) implies then the desired statement. Part 3.
is a simple consequence of the fact that $W$ is a complete
intersection. For part 4., consider a smooth refinement fan $\Sigma'$
of the simplicial fan $\Sigma.$ Since $W$ does not intersect
the singular locus of $X_{\Sigma},$ the proper 
transform $f^*(W)$ of $W$ under the
induced map $f: X_{\Sigma'} \to X_{\Sigma}$ does not intersect the 
exceptional locus in $X_{\Sigma'}.$ This means that for any toric 
divisor $D$ in  $X_{\Sigma'}$ corresponding to a new ray that is
in $\Sigma^\prime$ but not in $\Sigma,$ we have that $f_* (D)=0.$ 
Moreover, since we only deal with smooth objects
in $X_{\Sigma'},$ we have that
\begin{equation}
\Td_{f^*(W)} =\Bigl( \;\prod_{j=1}^k \frac{f^*(E_j)} {1- e^{-f^*(E_j)}}
\;\Bigr)^{-1}
\prod_{i \in \Sigma^\prime(1)} \frac{D_i}{1-e^{-D_i}}.
\end{equation}
But the restriction $f_{\; \vert \; _{f^*(W)}} : f^*(W) \to W$
is an isomorphism of smooth varieties, 
so $f_*(\Td_{f^*(W)}) = \Td_W,$ and the statement follows.
\end{proof}

The relations (\ref{linear}) involved
in the definition of the Stanley--Reisner ideal can be also
expressed by the equality
\begin{equation}
\lambda_1 \ov_1 + \ldots + \lambda_N \ov_N=0,
\end{equation}
which holds in the $\ZZ$-module
$\ZZ[\lambda_1, \ldots, \lambda_N]/ \Ia_{lin}.$
The significance of this simple fact will become clear
in our discussion of GKZ theory.

For a general triangulation $\Ta,$ it is known that the rank of the 
$\QQ$--algebra $\Ra_{\Ta} \otimes \QQ$ is equal to the number of 
maximal simplices in the triangulation. In particular, for a 
unimodular triangulation the  rank is given by the volume
of the polyhedron $\Delta^*.$

\subsection{GKZ systems for Calabi--Yau Complete Intersections}
\label{chap:gkzcicy}

We review some general facts about GKZ systems. Let
${\cal A}= \{ v_1, \ldots, v_N \}$ be an $\cA$--set 
in $\ZZ^n$ and $\LL \subset \ZZ^N$ the lattice of relations
between the elements of $\cal A.$

The GKZ system with parameter $\beta=(\beta_1, \ldots, \beta_N) \in \CC^N$ 
is a system of partial differential equations for functions 
$\Psi(z),$ $z \in \CC^N$ given by the following equations:
\begin{eqnarray}
\label{gkz1}
\left(-\beta+ \sum_{i=1}^N\:\ov_i\,z_i\frac{\partial}{\partial z_i}\right)\;
\Psi&=&0,\\
\label{gkz2}
\left(\prod_{l_i>0} \left[\frac{\partial}{\partial z_i}\right]^{l_i}
\:-\:\prod_{l_i<0} \left[\frac{\partial}{\partial z_i}\right]^{-l_i}
\right)\Psi&=&0,\hspace{3mm}\hbox{ for } \ l\in\LL,
\end{eqnarray}
where (\ref{gkz1}) is a system of $n$ equations.
The definition of the $\cA$--set 
implies that $l_1 + \ldots +l_N =0,$ for any $l=(l_1, \ldots, l_N) \in \LL,$
so the  equations (\ref{gkz2}) are homogeneous. 

In a series of papers, Gelfand, Kapranov and Zelevinsky proved
that the solutions of this system of partial differential equations form
a locally constant sheaf of finite rank outside a hypersurface in 
$\CC^N.$ They also wrote solutions of the system in
the form of the so-called $\Gamma$-series 
\begin{equation} \label{gaser}
\sum_{l\in\LL}\;\prod_{i=1}^N \frac{\:z_i^{\gamma_i+l_i}}
{\Gamma (\gamma_i+l_i+1)},
\end{equation}
where $\Gamma$ is the usual $\Gamma$-function,
$l=(l_1,\ldots,l_N) \in\LL\subset\ZZ^N$.
The series depends on  an additional parameter
$\gamma=(\gamma_1,\ldots,\gamma_N)\in\CC^N$ which must satisfy
\begin{equation} 
\gamma_1 \ov_1+\ldots+\gamma_N \ov_N\,=\,\beta.
\end{equation}

The case that we will consider  is the case in which the 
$\Gamma$-series incorporate all the periods of a top holomorphic form
on a Calabi--Yau complete intersection in a toric Fano variety.
From the point of view of the GKZ theory, our interest will focus
on the case in which the value $\beta$ is given by 
\begin{equation}
\beta=(-1, \ldots, -1, 0, \dots, 0)=-\ov_1-\ldots-\ov_k.
\end{equation}
In this case, the $\Gamma$-series (\ref{gaser}) become
\begin{equation}\label{phi}
\Phi_{\lambda}(z):= 
\sum_{l\in\LL}\;
\prod_{i=1}^k \frac{z_i^{\lambda_i+l_i-1}} {\Gamma (\lambda_i+l_i)} 
\prod_{i=k+1}^N \frac{z_i^{\lambda_i+l_i}} {\Gamma (\lambda_i+l_i+1)},
\end{equation}
with the parameter $\lambda=(\lambda_1,\ldots,\lambda_N)\in\ \CC^N$
satisfying
\begin{equation}\label{phicond} 
\lambda_1 \ov_1+\ldots+\lambda_N \ov_N\,=0.
\end{equation}
Let $\Ca$ be the cone in the secondary fan corresponding to a 
regular triangulation $\Ta$ of the set $\cal A$
and define the function
\begin{equation}\label{phicone}
\Phi^{\Ca}_{\lambda}(z):=
\sum_{l\in\Ca^{\vee}\cap \LL }\;
\prod_{i=1}^k \frac{z_i^{\lambda_i+l_i-1}} {\Gamma (\lambda_i+l_i)} 
\prod_{i=k+1}^N \frac{z_i^{\lambda_i+l_i}} {\Gamma (\lambda_i+l_i+1)},
\end{equation}
where the parameter $\lambda=(\lambda_1,\ldots,\lambda_N)\in\CC^N$
satisfies (\ref{phicond}), 
and $\Ca^\vee \subset \LL_\RR$ is the dual of the cone 
$\Ca \subset \LL^{\vee}_\RR,$ given by 
\begin{equation}
\Ca^\vee :=\{ l :
\langle w,l \rangle \geq 0, \ \hbox{for all} \ w \in \Ca \}.
\end{equation}

One of the results obtained by Gelfand, Kapranov and Zelevinzky
was the actual convergence of the series $\Phi^{\Ca}_{\lambda}(z)$
in a neighborhood of some distinguished point determined 
by the regular triangulation $\Ta.$ We review their approach using our
notation. 

The cone $\Ca \subset \LL^{\vee}_\RR $ is obtained as the intersection 
of all cones $\Ca_{I^*} \subset \LL^{\vee}_\RR$ with $I \in \Ta^n.$
For some  $I \in \Ta^n,$ we choose
the generators $\tilde{l}_1, \ldots, \tilde{l}_{N-n}$ of the 
cone $\Ca_{I^*}^{\vee}\cap \LL$
and consider the function
\begin{equation}
\Phi^{I}_{\lambda}(z):= 
\sum_{l\in\Ca_{I^*}^{\vee}\cap \LL  }\;
\prod_{i=1}^k \frac{z_i^{\lambda_i+l_i-1}} {\Gamma (\lambda_i+l_i)} 
\prod_{i=k+1}^N \frac{z_i^{\lambda_i+l_i}} {\Gamma (\lambda_i+l_i+1)}.
\end{equation}
This series can be written as 
\begin{equation}\label{ser1}
\Phi_{\lambda}^{I}(z)= (\prod_{i=1}^{k} z_i)^{-1} \cdot 
\sum_{m_i \geq 0} c_m \ x_1^{m_1}\ldots x_{N-n}^{m_{N-n}},
\end{equation}
where the coordinates $ (x_1, \ldots, x_{N-n})$ are given by 
\begin{equation}
x_i:=\prod_{j=1}^{N} z_j^{\tilde{l}_{ij}+\lambda_{ij}}, \ \hbox{for}\ 
i=1,\ldots, N-n.
\end{equation}
According to \cite{GKZ1}, for fixed values $(\lambda_1, \ldots,\lambda_N)
\in \CC^N$
and some  $I \in \Ta^n$ the series of type
(\ref{ser1}) are convergent in some neighborhood of the 
point $x_1=\ldots= x_{N-n}=0,$ and the regularity of the
considered triangulation comes into play and implies that all of them 
have a common domain of convergence (proposition 1.2 in \cite{GKZ1}).
If $p: \RR^{N\vee} \to \LL^{\vee}_\RR $ is the canonical
projection dual to the inclusion $\LL \hookrightarrow \ZZ^N,$ 
then the domain of convergence contains the subset 
of $\CC^N$ defined by
\begin{equation}\label{convz}
\big\{ (z_1, \ldots, z_N) : p(-\log |z_1|, \ldots, -\log |z_N|) \in \Ca +c
\big\} \ \hbox{for some} \ c \in \Ca.
\end{equation}

Recall (\ref{secon}) that the cone $\Ca$ in the secondary fan
corresponding to the regular triangulation $\Ta$ is defined
as the intersection of cones of type $\Ca_{I^*},$ for all $I$ maximal 
simplices in the triangulation. Then 
\begin{equation}
\Ca^{\vee}= \bigcup_{I \in \Ta^n} \Ca_{I^*}^{\vee}.
\end{equation}
This equality of sets implies the following inequality
\begin{equation}
\vert \Phi_{\lambda}^{\Ca}(z) \vert \leq \sum_{I \in \Ta^n}
\vert \Phi_{\lambda}^{I}(z) \vert,
\end{equation}
which means that the series $\Phi_{\lambda}^{\Ca}(z)$ is convergent
in a neighborhood of the point corresponding to the cone $\Ca$ in the 
toric variety determined by the secondary fan. One of the features 
of mirror duality is that
this toric variety provides a compactification of the moduli space
of complex structures on the mirror manifold of $W.$ Our discussion
implies that, by interpreting $z$ as a parameter on this complex 
moduli space, the series 
$\Phi^{\Ca}_{\lambda}(z)$ is convergent around a (generally singular)
boundary point given explicitly by the toric structure of this 
moduli space. 


In fact, following an idea used extensively by Yau and
his collaborators, Givental and 
many others, we want to view the function $\Phi_{\lambda}(z)$ as 
a formal power series. In order to be able to do that, we 
operate a change of variables given by 
\begin{equation}\label{def:tcoor}
z_j=e^{2 \pi i t_j}, \ \hbox{for all $j$, $ 1\leq j \leq N$}.
\end{equation}

For any formal parameter $\lambda,$ and any integer $m \in \ZZ,$
we can define the following formal function 
\begin{equation}
\Gamma(m + \lambda):=
\begin{cases} 
\Gamma(\lambda + 1)\cdot (\lambda +1)\cdot \ldots \cdot 
(\lambda+m-1)& \hbox{if $m \geq 1$},\\
\Gamma(\lambda + 1) \cdot (\lambda \cdot \ldots \cdot (\lambda +m))^{-1}& 
\hbox{if $m \leq 0$}.
\end{cases}
\end{equation}

With this convention, the function $(\Gamma(\lambda +m))^{-1}$ 
can be written as formal power series in 
$\lambda,$ where it is understood that the usual expansion
around $\lambda=0$ is used for $\Gamma(\lambda +1).$
Note that the formal definition will continue to hold for analytic
functions in the appropriate domain of values for $\lambda.$

The $\Gamma$-series 
\begin{equation}
\Phi_{\lambda}(t)=  \sum_{l\in\LL}\;
\prod_{i=1}^k \frac{e^{2 \pi i t_i(\lambda_i+l_i-1)}} 
{\Gamma (\lambda_i+l_i)} 
\prod_{i=k+1}^N 
\frac{e^{2 \pi i t_i (\lambda_i+l_i)}} {\Gamma (\lambda_i+l_i+1)}
\end{equation}
can be viewed as an element of the algebra of formal power series 
in the variables $\lambda_i,t_i,$ $1 \leq i \leq N,$ that is 
\begin{equation}
\Phi_{\lambda}(t) \in (\ZZ[[\lambda_1,\ldots, \lambda_N, t_1, \ldots, t_N]]
/ I_{lin}) \otimes \CC,
\end{equation}
where the ideal $I_{lin}$ is given by (\ref{linear}). 

For any $l \in \LL,$ and $\lambda=(\lambda_1, \ldots, \lambda_N)$
such that $ \lambda_1 v_1 + \ldots +\lambda_N v_N =0,$
define the formal function $P_{\lambda}(l)$ to be the coefficient 
of the series $\Phi_{\lambda}$
\begin{equation}\label{gammacoeff}
P_{\lambda}(l) :=  
\prod_{i=1}^k \frac{1} {\Gamma (\lambda_i+l_i)} 
\prod_{i=k+1}^N \frac{1} {\Gamma (\lambda_i+l_i+1)}.
\end{equation}

We prove a simple (but important) property.

\begin{proposition}\label{propdual}
Given a regular triangulation $\Ta$ and the corresponding cone
in the secondary fan $\Ca$, assume that $\lambda_1, \ldots,\lambda_N$ 
are the generators of the Stanley--Reisner ring of the triangulation
$\Ra_{\Ta}.$ For any $l \in \LL$ define
\begin{eqnarray}\label{ij}
I&:=& \{i: l_i \leq 0, 1\leq i \leq k \} \ \hbox{and} \nonumber\\
J&:=& \{j: l_j < 0, k+1 \leq j \leq N \}.
\end{eqnarray}

\begin{enumerate}
\item If $P_{\lambda}(l) \not= 0$ in the algebra $\Ra_{\Ta},$ then 
\begin{equation}\label{nonzero}
{\rm Conv}(I \cup J) \ \hbox{is a simplex in} \ \Ta.
\end{equation}

\item Any $l \in \LL$ such that (\ref{nonzero}) holds,
has the property that
\begin{equation}\label{dualcone}
l \in \Ca^{\vee}.
\end{equation}
\end{enumerate}
\end{proposition}

\begin{proof} Indeed, for any  $l \in \LL,$ we have that 
\begin{equation}
P_{\lambda}(l) = (\prod_{i \in I} \sin(\pi \lambda_i)) 
(\prod_{j\in I} \sin(\pi \lambda_j))
\;Q_l(\lambda),
\end{equation}
for some holomorphic function $Q_l(\lambda)$ (we have used the 
classical identity $\Gamma(x) \; \Gamma(1-x) =\frac{\pi}{\sin(\pi x)}$). 
Property (\ref{nonzero})
is then an immediate consequence of the relations (\ref{polyn})
in the definition of the Stanley--Reisner ring $\Ra_{\Ta}.$ To check
part 2., assume that for some  $l \in \LL$ we have that 
${\rm Conv}(I \cup J)$ is a simplex in  $\Ta.$ 
Then there exists a simplex $K \in \Ta^n$
such that $I \cup J \subset K.$ 
It follows  that $l_i \geq 0,$ for all $i \in K^*,$ so 
(\ref{dualcone}) follows directly from
\begin{equation}
\Ca^{\vee}= \bigcup_{K \in \Ta^n} \Ca_{K^*}^{\vee}.
\end{equation}
\end{proof}
As a direct consequence of this proposition we obtain the 
following equality for the $\Ra_\Ta \otimes \CC$-valued formal
functions
\begin{equation}
\Phi_{\lambda}(t) = \Phi^{\Ca}_{\lambda}(t).
\end{equation}
Denote by $\Ma$ the toric variety determined by the secondary fan. 
For any regular triangulation $\Ta,$ denote by $p_\Ta$ the
point corresponding to the corresponding big cone $\Ca$ in $\Ma.$
The toric variety $\Ma$ has a distinguished projective
embedding induced by the secondary polytope, whose
existence was proved by Gelfand, Kapranov and Zelevinsky \cite{GKZ2}.
As mentioned in the introduction, the toric variety $\Ma$
compactifies the so--called simplified polynomial moduli
space of the mirror Calabi--Yau space $M$ \cite{AGM2}.
There are various other ways to compactify the complex
moduli space of the mirror, but here we will consider only
the compactification induced by the secondary fan.

Following work of Yau and collaborators \cite{Yau} and
Givental \cite{Giv2}, Stienstra proved the following results.
An explicit solution basis to the GKZ system for the case
of a unimodular triangulation is also the subject of theorem
3.6.23 in the recent book \cite{SST}.
\begin{proposition}\cite{Stienstra}\label{stienstra}
\begin{enumerate} 
\item For any regular triangulation $\Ta$ of $\cA,$ the 
isomorphism 
\begin{equation}\label{isom}
\Ra_{\Ta} / Ann(\lambda_1 \cdot \ldots \cdot \lambda_k) \cong
(\lambda_1 \cdot \ldots \cdot \lambda_k) \cdot \Ra_{\Ta}
\end{equation}
implies that $\Phi^{\Ca}_{\lambda} (t)$ takes values in the 
algebra $\Ha \otimes \CC$ with
\begin{equation}\label{hdef}
{\cal H}:=
\Ra_{\Ta} / Ann(\lambda_1 \cdot \ldots \cdot \lambda_k). 
\end{equation}
 \item The map:
\begin{eqnarray}\label{inj}
\hbox{Hom} ({\cal H}, \CC)
 &\to& \ \hbox{Solution space of the GKZ system about} \ p_\Ta \nonumber \\
F  &\mapsto& F \cdot \Phi_{\lambda}(t)=F \cdot \Phi_{\lambda}^{\Ca}(t)
\end{eqnarray}
is injective.
\end{enumerate}
\end{proposition}
Part 1 of the proposition is a consequence of 
lemma 7, pg.18, in Stienstra's paper \cite{Stienstra} which shows 
that any regular 
triangulation $\Ta$ such that $v_1, \ldots, v_k$ belong to all
maximal simplices in $\Ta,$ has the property that
any $l = (l_1, \ldots,l_N)  \in \Ca^{\vee},$ 
is such that $l_j \leq 0,$ for $j=1,\ldots,k.$ The equality of 
the $\Ha$--valued series 
$\Phi_{\lambda}(t)=\Phi_{\lambda}^{\Ca}(t)$ is a consequence 
of proposition \ref{propdual}. combined with that lemma.
Part 2 looks a bit unsatisfactory, but in practice it is all that we need.
Gelfand, Kapranov and Zelevinsky proved that the dimension of the 
solution space is given by the volume of the polyhedron 
$\hbox{Conv}(\cA).$ In the case of
a unimodular triangulation, this volume is given by 
the number of maximal simplices in the triangulation. In general, 
the number of maximal simplices in the triangulation $\Ta$ is equal
to the rank of the vector space $\Ra_{\Ta} \otimes \CC
\cong H^*(X, \CC).$ The injective map (\ref{inj}) shows how to match 
the ``toric''
cohomology classes in $H^*(W,\ZZ)$ with an appropriate sub-module of 
the GKZ $\cal D$-module. This sub-module is the Picard--Fuchs 
$\cal D$-module of the Calabi--Yau complete intersection. The 
additional differential equations defining this sub-module
can be found using the Griffiths--Dwork method and this has been
pursued in the toric context for example in \cite{Yau}. We won't
need the complete set of differential equations here. 
In most explicit examples, one can usually find by direct inspection 
the complete system of differential equations. 

For a unimodular triangulation $\Ta,$ the $\Ra_\Ta \otimes \CC$--valued
series $\Phi_\lambda$ provides a complete system of solutions to the GKZ
system. This is not true anymore for a non--unimodular triangulation
$\Ta,$ since the rank of the vector space $\Ra_{\Ta} \otimes \CC$
(the number of maximal simplices in $\Ta$) obviously 
drops. However, it turns out that in certain cases when the Calabi--Yau
complete intersection is smooth, the last statement 
of proposition \ref{stienstra} provides
a complete set of solutions to the Picard--Fuchs equations
for the periods of the mirror manifold around the boundary point 
$p_\Ta$ in $\Ma$ (the periods form a subset of the solutions to the 
GKZ system). Of course, when the triangulation becomes even 
coarser, the above result, although still true, fails to provide
enough solutions to the Picard--Fuchs equations. In an extreme 
(but interesting) situation this rank can be  1 (this the case 
for the so-called Landau--Ginzburg point). One of goals in the next 
section will be to overcome this difficulty and construct a complete 
set of of solutions corresponding to a general regular triangulation 
using analytic continuation.


\section{Analytic Continuation and Monodromy \\ Formulae}
\label{chap:mon}
We now explain our approach towards constructing solutions
to the GKZ system about special points in the (singular) toric 
variety $\Ma$ determined by the secondary fan. 
These points
are determined by the maximal dimensional cones and are
fixed under the maximal torus action.
The method will involve an analytic continuation process based on 
Mellin--Barnes integral representations of hypergeometric functions.

According to proposition \ref{toricoh}, the ring $\Ha$ defined
by (\ref{hdef}) describes the toric part of the even cohomology 
of the Calabi--Yau complete intersection $W.$ 
Part 2 of proposition \ref{stienstra} shows how the ring $\Ha$ plays 
the r\^ole of a bookkeeping device for the solutions of the 
GKZ system around the point $p_\Ta$ in $\Ma$ corresponding to the regular
triangulation $\Ta.$ This approach is nothing else but 
the generalization of the classical Frobenius method for 
ordinary differential equations and, as mentioned already, has been used 
extensively by Yau and collaborators, Givental and many others. It is 
tempting to try to use the analogue of the ring $\Ha$ for any regular
triangulation to obtain a complete set of solutions to the Picard--Fuchs
equations of the mirror Calabi--Yau $M.$ As mentioned at the end
of the previous section, this cannot work for coarser triangulations
$\Ta$ for dimensional reasons. Assume that a regular triangulation $\Ta_0$ 
of $\cA$ is obtained as in section \ref{chap:cicy} from a simplicial fan
$\Sigma \subset \Nrm_\RR,$ and assume that the corresponding nef--partition
defines a {\it smooth} Calabi--Yau complete intersection 
$W \subset X_\Sigma.$ Our method of analytic continuation will show 
that in fact it is natural to retain the ring $\Ha \subset H^*(W, \ZZ)$ 
corresponding to $\Ta_0$ 
as a bookkeeping device for describing the solutions around {\it any} 
toric fixed point $p_{\Ta}$ corresponding to a general regular 
triangulation $\Ta.$ In all the examples that will be considered later in
this work, it is true that  $\Ha \cong H^*(W, \ZZ),$ but in 
general, the two rings may differ by some non--toric classes.
The use of the $t$--coordinates (defined by (\ref{def:tcoor}))
will
allow us to work with analytic 
functions instead of multivalued functions.

Consider a sequence of big cones $\Ca_0, \ldots, \Ca_D$ in the secondary
fan such that each $\Ca_d$ corresponds to some regular triangulation
$\Ta_d$ of  $\cA,$ and, for all $d,$ $0 \leq d \leq D-1,$ the cones 
$\Ca_{d}$ and $\Ca_{d+1}$ are adjacent. Such a sequence corresponds
to a path along the edges of the secondary polytope. 
We assume that $\Ca_0$ corresponds to the triangulation 
$\Ta_0$ chosen above. Denote by
$p_d$ the special point determined by the cone $\Ca_d$ in the toric 
variety $\Ma$ defined by the secondary fan. In the $t$-coordinates, 
$z_j=e^{2 \pi i t_j}, j=1,\ldots,N,$ the convergence domain
(\ref{convz}) of the series $\Phi_{\lambda}^{\Ca_d}(t)$ is 
\begin{equation}\label{def:convt}
U_d:=\big\{ (t_1, \ldots, t_N) : p (\Im t_1, \ldots, \Im t_N) 
\in \Ca_d + c
\big\} \ \hbox{for some} \ c \in \Ca_d,
\end{equation}
where $p: \RR^{N\vee} \to \LL_{\RR}^\vee$ is the projection
dual to the inclusion $\LL \hookrightarrow \ZZ^N.$

For a complex manifold $U,$ we denote by $\Oa(U)$ the ring of 
germs of analytic functions on $U.$ Consider the parameters
$\lambda_1, \ldots, \lambda_N,$ as being the coordinates
on some space $\CC^N.$ 
The intersection of all the hyperplanes
in  $I_{lin}$ is an  $(N-n)$-dimensional linear subspace
$H \cong \CC^{N-n} \subset \CC^N.$ 
Given an open neighborhood $V$ of the origin in $H,$ there exists a 
natural projection $\Pi: \Oa(V) \to \Ha$ 
such that for any 
holomorphic function $\phi \in \Oa(V),$ the function
$\Pi \cdot \phi$ is obtained by considering the 
power series expansion at the origin of $H$ of 
the function $\phi.$ We will also use a ``family'' version of this 
construction. Namely, for a function $\Phi_\lambda(t) \in \Oa(W \times U)$
the projection $\Pi \cdot \Phi_\lambda(t)$ is obtained by considering
the $\Oa(U)$--valued expansion of $\Phi_\lambda(t)$ at the origin
of $H.$
We also need to consider a special type
of meromorphic functions on $H$. Consider a finite subset $S$ 
of the ideal $I_{lin}.$ Each relation in $S$ determines 
a linear form $\alpha$ on $\CC^N$ that vanishes on
the subspace $H.$ We denote by $\Ma_S (H)$ the set  
of meromorphic functions on $H$ with {\it finite order}
poles contained in the set 
\begin{equation}
H \cap \bigcup_{q \in \ZZ \setminus \{0\}, \alpha \in S} \ker (\alpha +q) 
\; .
\end{equation}
Implicitly, any function in $\Ma_S (H)$ is analytic
in $V,$ for some neighborhood $V$ of the origin 
in $H,$ so we can also view the map $\Pi$ as acting on
$\Ma_S (H).$
 
If $\lambda_1, \ldots, \lambda_N$ represent the coordinate
hyperplanes in $\CC^N,$ and the relations in $I_{lin}$ define 
the linear $(N-n)$-dimensional subspace $H,$ we define the
finite set of linear forms $S \subset I_{lin}$ by
\begin{equation}
S:=\{ h_1\lambda_1 + \ldots  + h_N \lambda_N: h_j \in \ZZ, 
\vert h_j \vert \leq m,\; j=1,\ldots,N \}.
\end{equation}
The proof will show the natural number $m$ only depends
on the combinatorial structure of the set of vectors $\cA$
given by the ideal $I_{lin}$ (or, equivalently, by the lattice 
$\LL$). 

We need another technical definition. Recall from section \ref{chap:triang}
that the elements of the Gale transform of the set 
$\cA = \{\ov_1, \ldots, \ov_N \}$
were denoted by $w_i \in \LL^\vee,$ $1 \leq i \leq N.$ 
A set $I \subset \{1, \ldots, N \}$
is said to be {\bf saturated} if any $i \in \{ 1, \ldots, N \}$ such 
that there exists $j \in I$ with $w_j = t w_i, t > 0,$ has the 
property that $i \in I.$

Recall (\ref{cone99}) that 
\begin{equation}
{\cal C}_J := \{ \sum_{j \in J} t_j w_i: t_j \in \RR_{\geq 0} \}.
\end{equation}
The proof of the theorem will use the following assumption. 

\begin{condition}\label{cond1}
There exist no two subsets $J^\prime$ and $J''$ of $\{1, \ldots, N \}$ 
of cardinality $N-n-1,$ with both $w_j, j \in J^\prime$ and
$w_j, j \in J''$ collections of linearly independent vectors,
such that the cone $\Ca_{J^\prime}$ is a proper subcone of
$\Ca_{J''}.$ 
\end{condition}

Note that the condition is automatically satisfied for $N-n=1,2.$ 
At the time of this writing, the author does not know how to remove this 
assumption.

We prove the following:

\begin{thm}\label{gensol}
For any $d,$ $1 \leq d \leq D,$ there exists an   
``analytic continuation'' operator $\Pa_d : \Ma_S (H) \to 
\Ma_S(H),$ such that 
the following properties hold:
\begin{enumerate}
\item For all saturated $I \subset \{1, \ldots, N \}$ such that 
$I$ is not a simplex in the triangulation $\Ta_d$ (this is equivalent to  
$\prod_{j \in I}\lambda_j \in I_{mon}^{\Ta_d}$) and any 
function $\phi \in \Ma_S(H),$  
\begin{equation}\label{p1}
\Pi \cdot \Pa_d \cdot (\prod_{j \in I}\sin(\pi\lambda_j))\; 
\phi(\lambda) = 0.
\end{equation}
\item The map:
\begin{eqnarray}\label{p4}
\hbox{Hom} ({\cal H}, \CC)
 &\to& \ \hbox{Solution space of the GKZ system about} \ p_d \nonumber \\
F  &\mapsto& F \cdot \Pi \cdot \Pa_d \cdot
\Phi_{\lambda}^{\Ca_d}(t) 
\end{eqnarray}
is injective.
\end{enumerate}
\end{thm}

The theorem provides the general form of the solution space of the GKZ system 
around a general boundary point in the toric 
variety given by the secondary fan. It shows that the analytic
continuation of the $\Ha$--valued series  $\Pi \cdot \Phi_{\lambda}^{\Ca_0}$
which gives the solutions to the GKZ system around the point
$p_0$ (cf. proposition \ref{stienstra}) is the $\Ha$--valued series 
$\Pi \cdot \Pa_d \cdot \Phi_{\lambda}^{\Ca_d}$ which gives the 
solutions to the 
GKZ system around the point $p_d.$ One should note that 
condition (\ref{p1}) is equivalent to the combinatorial fact that
the product $\prod_{j \in I}\lambda_j$ is zero in the Stanley--Reisner 
ring $\Ra_{\Ta_d}$ of the triangulation $\Ta_d.$ 
Hence part 1 of the theorem shows exactly which terms in the formal
power series $\Phi_{\lambda}$ contribute to the solutions of the
GKZ system around the point $p_d.$ 
The proof of the theorem will also provide the strategy and the explicit 
formulae that will help us perform the desired monodromy calculations 
in the next section.

\begin{proof} The proof is by induction. For $d=0$ the operator 
$\Pa_0$ is just the identity operator. 
The properties 1 and 2 follow directly  from propositions 
\ref{propdual} and \ref{stienstra}. 

To prove the induction step, we consider
two adjacent cones $\Ca_d$ and $\Ca_{d+1}$ and let $\Ta_d$ and 
$\Ta_{d+1},$ be the corresponding regular triangulations. 
Assume that the analytic continuation operator $\Pa_d$ 
corresponding to the cone $\Ca_d$ exists and has the 
expected properties. 
 
Proposition \ref{prop:perestr}
shows that there exists a circuit $I$ in $\cal A$ such that
the two triangulations are supported on $I$ and $s_I (\Ta_d)=\Ta_{d+1}.$
According to lemma 1. the hyperplane in $\LL_{\RR}^\vee$ 
separating the cones $\Ca_d$ and
$\Ca_{d+1}$ is defined by an element 
$h=(h_1,\ldots,h_N) \in \LL$ with $h_j= \langle w_j, h \rangle,$
where $w_j$ are the elements of the Gale transform of $\cal A$
for $j=1,\ldots,N.$ We choose the element $h$ to be minimal
in the lattice $\LL,$ that is so that no rational non--zero multiple 
of $h$ with smaller length is in $\LL.$ 

According to proposition \ref{katzcox} there
exists a class $\theta \in H_2 (X, \ZZ)$ such that 
\begin{eqnarray}\label{theta2}
h_j &=& - E_j \cdot \theta, \ \hbox{for} \ j=1,\ldots, k, \nonumber\\
h_j &=& D_j \cdot \theta, \ \hbox{for} \ j=k+1, \ldots, N,
\end{eqnarray}
where $D_j$ are the toric divisors on the ambient toric variety $X_\Sigma$
defined by the simplicial fan $\Sigma.$ As usual,
define the subsets $I_+= \{ j: h_j > 0 \}$ and $I_-= \{j: h_j < 0 \},$
so that the triangulation $\Ta_d$ contains all the simplices of type
Conv$\{I \setminus j: j \in I_+ \}$ and the triangulation $\Ta_{d+1}$ 
contains all the simplices of type Conv$\{I \setminus j: j \in I_- \}.$

As mentioned in the remarks after lemma \ref{lem}, 
any saturated $h$ can be included as an element of an
integral basis for the lattice $\LL.$ The induced decomposition 
\begin{equation}\label{decomp}
\LL= \langle h \rangle \oplus \LL^\prime
\end{equation}
implies that any $l=(l_1, \ldots, l_N) \in \LL$ can be written as 
\begin{equation}\label{ldecomp}
l_j= m \cdot h_j + l_j^\prime, \ \hbox{with} \ m \in \ZZ, 
l^\prime \in \LL^\prime, j=1,\ldots,N.
\end{equation}
Since the parameters $\lambda_1, \ldots,
\lambda_N$ are coordinates on $\CC^N$
the decomposition (\ref{decomp}) leads to some choice of
coordinates on the linear subspace $H \subset \CC^N$ 
such that any $\lambda=(\lambda_1, \ldots, \lambda_N) \in H$  
can be written as 
\begin{equation}
\lambda_j= \mu \cdot h_j + \mu_j^\prime, \ \hbox{with} \ 
\mu^\prime \in H^\prime, j=1,\ldots,N.
\end{equation}

Each $\lambda_j$ is naturally identified with an element
in $H^2(X, \ZZ)$ and as a consequence of (\ref{theta2}) the parameter 
$\mu$ represents the class in $H^2(X, \ZZ)$ that is dual to 
$\theta \in H_2 (X, \ZZ).$

Since the big cones $\Ca_d$ and $\Ca_{d+1}$ have a common face, their
dual cones $\Ca^{\vee}_d$ and $\Ca^{\vee}_{d+1}$ in $\LL_\RR$ 
project according to the decomposition (\ref{decomp}) 
onto the same cone  $\Ca_{d,d+1}^\vee$ in $\LL^{\prime}_\RR.$
The existence of such a cone is a consequence of the important
result of Gelfand, Kapranov and Zelevinsky which states that
the secondary fan is the normal fan of some convex 
polytope called the secondary polytope. The
hyperplane that contains the common face of the
cones $\Ca_d$ and  $\Ca_{d+1}$ corresponds to an edge 
of the secondary polytope joining the vertices corresponding
to the regular triangulations $\Ta_d$ and $\Ta_{d+1}.$ The 
cone $\Ca_{d,d+1}^\vee$ can be taken to be the cone obtained
by intersecting the secondary polytope with a section
normal to the edge. 

Let $I_0^\prime$ the set of elements in $I_0$ that are
generators of the cone $\Ca_{d,d+1}^\vee$ which is the  
common face of the of the cones
$\Ca_d$ and $\Ca_{d+1}$. For any element $l^\prime \in 
 \Ca_{d,d+1}^\vee$ we have that 
\begin{equation}
l_j=l_j^\prime \geq 0 \quad \hbox{for all $j \in I_0$},
\end{equation}
and the following facts hold
\begin{eqnarray}\label{lconecond11}
l \in \Ca^{\vee}_d &\quad& \hbox{if and only if $l_j \geq 0$
for some $j \in I_+$}, \\ \label{lconecond21}
l \in \Ca^{\vee}_{d+1} &\quad& \hbox{if and only if $l_j \geq 0$
for some $j \in I_-$}.
\end{eqnarray}
Hence there exist two integers $m_1$ and $m_2$ such that 
for any $l \in \LL$, decomposed as in (\ref{ldecomp}) with 
fixed $l^\prime \in \Ca_{d,d+1}^\vee,$
\begin{eqnarray}\label{lconecond12}
l \in \Ca^{\vee}_d &\quad& \hbox{if and only if} \quad m \geq m_1,\\
\label{lconecond22}
l \in \Ca^{\vee}_{d+1} &\quad& \hbox{if and only if} \quad m \leq m_2.
\end{eqnarray} 
We emphasize here that, in general, the integers $m_1$ and $m_2$ 
{\it depend} (piecewise) linearly on the value $l^\prime.$ We have seen in 
the previous section that the series $\Phi_{\lambda}^{\Ca_d}(t)$
(see (\ref{phicone})) is convergent in the neighborhood 
of the point $p_d$ for any complex parameters $\lambda_1,
\ldots, \lambda_N$ such that $\lambda_1 \ov_1 + \ldots 
+ \lambda_N \ov_N= 0.$ The induction hypothesis implies that 
for any evaluation $F \in \hbox{Hom}(\Ha, \CC)$ the function
$F \cdot \Pa_d \cdot \Phi_{\lambda}^{\Ca_d}$ gives the solutions 
to the GKZ system in a neighborhood of the point $p_d.$
We have that 
\begin{equation} \label{form0}
\begin{split}
&\Pa_d \cdot \Phi_{\lambda}^{\Ca_d}(t)=\\
&=\Pa_d \cdot \sum_{l^\prime \in \Cpr} 
(\sum_{m \geq m_1} 
\prod_{j=1}^N \frac{1}{\Gamma(h_j(m+\mu)+l_j^\prime +\mu_j^\prime+1)} 
\; x^{m+\mu})
\; \prod_{j=1}^N e^{2 \pi i t_j (l_j^\prime+\mu_j^\prime)},
\end{split}
\end{equation}
where 
\begin{equation}\label{def:x}
x:= \prod_{j=1}^N e^{2 \pi i t_j h_j}.
\end{equation}
The coordinate $x$ is  well defined in the complex
plane considered with a branch cut and can be viewed as the
coordinate on the rational curve joining the points $p_d$ and
$p_{d+1}$ in the toric variety $\Ma$ determined by the 
secondary fan (see also section 5.3 in \cite{AGM3}). According to our
discussion in the previous section, the series $\Phi_{\lambda}^{\Ca_d}(t)$
is absolutely convergent in a neighborhood of the point $p_d,$
and there we can operate without restraint any change in the order of
the summation. For notational reasons, in the formula (\ref{form0}) 
we applied the change 
\begin{equation}\label{changel}
l_j^\prime \mapsto l_j^\prime -1, \ \hbox{for} \ 1 \leq j \leq k,
\end{equation}
which amounts to having $l^\prime \in \LL^\prime$ vary in a 
a cone denoted by $\Ca^{\prime \,\vee }$ 
obtained by a translation by a constant integer vector
of the cone $\Ca_{d,d+1}^\vee.$ 

In general, the series
\begin{equation}
\sum_{m \geq m_1} 
\prod_{j=1}^N \frac{1}{\Gamma(h_j(m+\mu)+l_j^\prime +\mu_j^\prime+1)} 
\; x^{m+\mu}
\end{equation}
is convergent in an annular region with a branch cut.
As mentioned before, the branch cut can be 
replaced by a restriction on the real part of the 
$t$--parameters.
To pursue the analytic continuation we use 
the classical method of Mellin--Barnes integral 
representations. The following property holds (see \cite{bateman}, page 49, 
and \cite{WW}, \S 14.5).

\begin{lemma}\label{lemma:melbar}
Consider the integral
\begin{equation}
\int_{\gamma- i \infty}^{\gamma +i \infty}
 \frac{\prod_{j=1}^p \Gamma(A_j s + a_j)}
{\prod_{j=1}^q \Gamma(B_j s + b_j)} \; y^s \; ds,
\end{equation}
with $\gamma, A_j, B_j$ all real. The path of integration 
is parallel to the imaginary axis for $|s|$--large, but it 
can be curved elsewhere so that it avoids the poles of the 
integrand. Introduce the following notations
\begin{equation}
\begin{split}
\alpha&:=\sum_{j=1}^p \vert A_j \vert - \sum_{j=1}^q \vert B_j \vert, \quad
\beta:= \sum_{j=1}^p  A_j  - \sum_{j=1}^q  B_j,\\
\eta&:= \Re (\sum_{j=1}^p  a_j  - \frac{1}{2} p - \sum_{j=1}^q  b_j
+\frac{1}{2} q), \quad
\rho:= (\prod_{j=1}^p \vert A_j \vert^{-A_j}) 
(\prod_{j=1}^q \vert B_j \vert^{B_j}),
\end{split}
\end{equation}
and
\begin{equation}
s=v + iu, \ y=R e^{i\theta}.
\end{equation}
Then the absolute value of the integrand has the asymptotic form 
\begin{equation}
e^{-\frac{1}{2} \alpha \pi \vert u \vert } \; \vert u \vert\
^{\beta v +\eta}\; R^\gamma\; e^{\theta u}\; \rho^{-\gamma}
\end{equation}
when $|u|$ is large (this is a consequence of Stirling's formula).
Therefore, if 
\begin{equation}\label{converg}
\alpha > 0,
\end{equation}
the integral is absolutely convergent (and defines
an analytic function of $y$) in any domain contained in  
\begin{equation}
\vert \arg y \vert < \min(\pi, \frac{\alpha\pi}{2}).
\end{equation}
Moreover, if $\beta=0,$ the integral is equal to the sum
of the residues on the right of the contour for $|y| < \rho,$
and to the sum of the residues  on the left of the contour for
$|y| > \rho$ (these facts are obtained by closing the contour to
the right, respectively to the left, with a semicircle of
radius $r \to \infty$). 
\end{lemma}
To make use of the lemma, we write
\begin{equation}\label{form1}
\begin{split}
&\Phi_\lambda^{\Ca_d}(t)=\sum_{l^\prime \in \Cpr} 
(\sum_{m \geq m_1} 
\prod_{j=1}^N \frac{1}{\Gamma(h_j(m+\mu)+l_j^\prime +\mu_j^\prime+1)} 
\; x^{m+\mu})
\; \prod_{j=1}^N e^{2 \pi i t_j (l_j^\prime+\mu_j^\prime)}=\\
&=\sum_{l^\prime \in \Cpr} \sum_{m \geq m_1} \prod_{j \in I_-}
\frac{\sin(\pi (-h_j (m+\mu) - l_j^\prime - \mu_j^\prime))}{\pi}
\Gamma (-h_j (m+\mu) - l_j^\prime - \mu_j^\prime)\cdot\\
&\cdot \prod_{j \notin I_-} \frac{1}
{\Gamma (h_j (m+\mu) + l_j^\prime + \mu_j^\prime+1)}
\; x^{m+\mu} \prod_{j=1}^N e^{2 \pi i t_j (l_j^\prime+\mu_j^\prime)}=\\
&=\sum_{l^\prime \in \Cpr}  
\prod_{j \in I_-}
\frac{\sin(\pi (-h_j \mu - l_j^\prime - \mu_j^\prime))}{\pi}
\prod_{j=1}^N e^{2 \pi i t_j (l_j^\prime+\mu_j^\prime)}
\cdot \\
&\cdot  \sum_{m \geq m_1}    (-1)^{m (\sum_{j \in I_-}h_j)}
\frac{\prod_{j \in I_-}
\Gamma (-h_j (m+\mu) - l_j^\prime - \mu_j^\prime)}
{\prod_{j \notin I_-}  \Gamma (h_j (m+\mu)+ l_j^\prime+ \mu_j^\prime+1)}
 \; x^{m+\mu}.
\end{split}
\end{equation}

In order to analytically continue the function $\Phi_\lambda^{\Ca_d}(t)$
from a neighborhood of $p_d$
to a neighborhood of $p_{d+1},$ we consider the following 
Mellin--Barnes integral
(compare to the last line in the previous formula)
\begin{equation}\label{int:barnes}
e^{i \pi \epsilon \mu} \frac{1}{2\pi i}
\int_{m_2 - i\infty}^{m_2 +i \infty}
\frac{\prod_{j \in I_-} \Gamma (-h_j (s+\mu) - l_j^\prime - 
\mu_j^\prime) \Gamma(-s) \Gamma(1+s)}
{\prod_{j \notin I_-}  \Gamma (h_j (s+\mu) + l_j^\prime + \mu_j^\prime+1)}
(e^{-i\pi\epsilon} x)^{s+\mu}\, ds
\end{equation}
with 
\begin{equation}\epsilon:=
\begin{cases}
1& \hbox{if}\ \sum_{j \in I_-}h_j \ \hbox{is even},\\
0& \hbox{if}\ \sum_{j \in I_-}h_j \ \hbox{is odd}.
\end{cases}\label{def:epsilon}
\end{equation}
The path of integration avoids all the poles of the integrand,
is parallel to the imaginary axis at large $|s|,$
and indented such that the poles of the function 
$\prod_{j \in I_-} \Gamma (-h_j (s+\mu) - l_j^\prime - 
\mu_j^\prime)$ lie on the left of the path. Such a choice of the contour
with a ``finite'' indentation is possible, given the properties 
(\ref{lconecond21}) and (\ref{lconecond22}) that define $m_2.$ Although 
the integer $m_2$ depends (piecewise) linearly on
$l^\prime \in \Cpr,$ we can choose indentations that are ``similar''
for all $l^\prime.$ With such a choice, 
the change of variable $s \mapsto s+m_2,$ allows us to  
rewrite the previous integral as 
\begin{equation}\label{integral1}
\begin{split}
& e^{i \pi \epsilon \mu} \frac{1}{2\pi i} 
\int_{- i\infty}^{+i \infty}
\frac{\prod_{j \in I_-} \Gamma (-h_j (s+m_2+\mu) - l_j^\prime - 
\mu_j^\prime)}
{\prod_{j \notin I_-}  \Gamma (h_j (s+m_2+\mu) + l_j^\prime + \mu_j^\prime+1)}
\cdot \\
&\cdot \Gamma(-(s+m_2))\, \Gamma(1+s+m_2)\, (e^{-i\pi\epsilon} x)^{s+m_2+\mu} 
ds,
\end{split}
\end{equation}
where the contour of integration is now {\it independent} of $l^\prime$
and coincides with the imaginary axis for large $|s|.$

The relation $h_1 \ov_1 +\ldots + h_N \ov_N=0, $ 
combined with condition (\ref{hyperplane1}), give that 
\begin{equation}
-\sum_{j \in I_-} h_j= \sum_{j \notin I_-} h_j.
\end{equation}
This shows that, in the notation of the lemma \ref{lemma:melbar},
$\alpha=2$ and $\beta=0$. So the integral (\ref{integral1}) is
absolutely convergent in any domain included in
\begin{equation}
\vert \arg(x) - \epsilon \pi \vert < \pi
\end{equation}
and this defines the branch cut in the $x$--plane. We also 
conclude that the path of integration can  
can be closed both to the right for $|x|$ small, and to the left for
$|x|$ big. 

An important remark to make is that the residues at poles 
$s=m \in \ZZ, m+m_2 < m_1$ of the function $\Gamma(-(s+m_2)) \Gamma(1+s+m_2)$ 
can be discarded for the purposes of the analytic continuation, since
the corresponding terms are canceled out by the operator $\Pa_d$. 
Indeed, the poles $s=m \in \ZZ,$ $m+m_2 < m_1,$  have residues 
which give factors of the form
\begin{equation}
\begin{split}
P_{\lambda}(l)&=\frac{\prod_{j \in I_-}\sin(\pi (-h_j \mu - l_j^\prime - 
\mu_j^\prime))\Gamma (-h_j (m+m_2+\mu) - l_j^\prime - \mu_j^\prime)}
{\prod_{j \notin I_-}  \Gamma (h_j (m+m_2+\mu) + l_j^\prime + \mu_j^\prime +1)}
\\
&= \frac{1}{\prod_{j=1}^N \Gamma (h_j (m+m_2+\mu) + l_j^\prime + 
\mu_j^\prime +1)}.
\end{split}
\end{equation}
The inequality $m+m_2 < m_1$ and conditions (\ref{lconecond11}) and 
(\ref{lconecond12}) show that the function $P_{\lambda}(l)$ can be written as
\begin{equation}
P_{\lambda}(l)= \prod_{j \in I_{+}} \sin(\pi(h_j \mu+ \mu_j^\prime))\; 
\phi(\lambda),
\end{equation}
with $\phi \in \Oa(H).$ We know that the triangulation $\Ta_d$ 
contains simplices of type 
$\hbox{Conv}\{I \setminus j : j \in I_+ \},$
so $I_+$ is not a simplex in $\Ta_d.$ Clearly $I_+$ is 
saturated, so the induction hypothesis 
implies then that
\begin{equation}
\Pi \cdot \Pa_d \cdot P_{\lambda}(l) = 0.
\end{equation}
This shows the residues at integer points  $s=m, m+m_2 < m_1$ 
do not contribute to the result of the analytic continuation of the
solutions to the Picard--Fuchs equations, 
and this is true independently of the location of the poles on the left or
right side of the path of integration.

We now study the evaluation of the integral
(\ref{integral1}) using residues, when we close the contour to the right
(for $|x| < \rho$), or to the left (for $|x| > \rho$). 
Assume first that {\bf \large (\dag)}
all the poles $ s=m \in \ZZ$ of the function 
$\Gamma(-(s+m_2)) \Gamma(1+s+m_2)$  with $m+m_2 \geq m_1,$
are located on the right of the contour of integration. 

This means that
all the poles $s \geq m_1$ of $\Gamma(-s) \Gamma(1+s)$ are located
on the right of the contour of integration in the integral
(\ref{int:barnes}), and the sum of all their residues is
\begin{equation}
\sum_{m \geq m_1}    (-1)^{m (\sum_{j \in I_-}h_j)}
\frac{\prod_{j \in I_-}
\Gamma (-h_j (m+\mu) - l_j^\prime - \mu_j^\prime)}
{\prod_{j \notin I_-}  \Gamma (h_j (m+\mu)+ l_j^\prime+ \mu_j^\prime+1)}
 \; x^{m+\mu}, 
\end{equation}
which is exactly the last line in the formula (\ref{form1}). 
The series $\Phi_{\lambda}^{\Ca_d}(t)$ written as in (\ref{form1})
is absolutely convergent in the neighborhood $U_d$ of the point $p_d$
corresponding to $|x|$--small. Lebesgue's dominated
convergence theorem shows that, for $t \in U_d$ (i.e. $|x|$--small), 
we can express $\Phi_{\lambda}^{\Ca_d}(t)$ by the following 
absolutely convergent integral
\begin{equation}\label{int:theta}
\Phi_{\lambda}^{\Ca_d}(t)= \frac{1}{2\pi i} \int_{- i\infty}^{+i \infty}
\Theta_{\lambda} (t,s) \;ds,
\end{equation}
with 
\begin{equation}\label{def:theta}
\begin{split}
\Theta_{\lambda} (t,s):=&
\sum_{l^\prime \in \Cpr}  
\prod_{j \in I_-}
\frac{\sin(\pi (-h_j \mu - l_j^\prime - \mu_j^\prime))}{\pi}
\prod_{j=1}^N e^{2 \pi i t_j (l_j^\prime+\mu_j^\prime)}e^{i \pi \epsilon \mu}
\cdot\\
&\cdot \frac{\prod_{j \in I_-} \Gamma (-h_j (s+m_2+\mu) - l_j^\prime - 
\mu_j^\prime)} {\prod_{j \notin I_-}  
\Gamma (h_j (s+m_2+\mu) + l_j^\prime + \mu_j^\prime+1)} \cdot\\
&\cdot \Gamma(-(s+m_2))\; \Gamma(1+s+m_2)\; (e^{-i\pi\epsilon} x)^{s+m_2+\mu}.
\end{split}
\end{equation}

Consider the coordinates $x_1, ,\ldots, x_{N-n-1}, 
x_{N-n}=x$ in $\Ma$ such that the rational  curve joining the points 
$p_d$ and $p_{d+1}$ is given by $x_1=\ldots= x_{n-n-1}=0,$ with the 
point $p_d$ given by $x_1=\dots=x_{N-n-1}=x=0$ (see \cite{AGM3}
for a detailed discussion of coordinates). We argue now that
the integral in the formula (\ref{int:theta}) is in fact 
absolutely convergent in a ``tubular'' neighborhood of the 
rational curve given by $x_1=\dots=x_{N-n-1}=0.$ Indeed, choose some
$\xi >0$ such that $U_d$ contains a product of $(N-n)$ complex disks 
with radii equal to $\xi.$ For an arbitrary coordinate $x,$ choose a real 
number $R >0$ such that $\vert x/R \vert < \xi.$ We can write
\begin{equation}
(\prod_{j=1}^N e^{2 \pi i t_j l_j^\prime}) \, x^{m_2}=
(\prod_{j=1}^{N-n-1} x_j^{p_j}) \, x^{m_2}= 
(R^{m_2} \prod_{j=1}^{N-n-1} x_j^{p_j}) \, \big(\frac{x}{R}\big)^{m_2},
\end{equation} 
where the integers $p_j$ are determined by the integers 
$l_j^\prime.$ The important fact to remember is that the
integer $m_2$ depends (piecewise) linearly on the integers
$l_j^\prime,$ i.e. on the integers $p_j.$ This means that we
can choose the coordinates $x_1, \ldots, x_{N-n-1}$ to be
small in absolute value, such that there exist
complex numbers $\tilde{x}_1,\ldots, \tilde{x}_{N-n-1},$ with
$|\tilde{x}_j| < \xi$ for $1 \leq j \leq N-n-1,$ and
\begin{equation}
(R^{m_2} \prod_{j=1}^{N-n-1} x_j^{p_j})= 
\prod_{j=1}^{N-n-1} \tilde{x}_j^{p_j}.
\end{equation}
If $\tilde{t}_j, 1\leq j \leq N,$ are the $t$--coordinates of the point
$(\tilde{x}_1, \ldots, \tilde{x}_{N-n-1}, x/R) \in U_d,$ it follows that
for $s=i u, u \in \RR,$ we have that
\begin{equation}
e^{2 \pi i t_j (l_j^\prime+\mu_j^\prime)} x^{s+m_2+\mu} 
= R^{i u +\mu} e^{2 \pi i \tilde{t}_j (l_j^\prime+\mu_j^\prime)}
\big(\frac{x}{R}\big)^{s+m_2+\mu}.
\end{equation}
Since the function $\Theta_\lambda(\tilde{t},s)$
is a convergent series of absolutely integrable functions,
Lebesgue's dominated convergence theorem implies indeed 
that the integral
\begin{equation}
\frac{1}{2\pi i} \int_{- i\infty}^{+i \infty}
\Theta_{\lambda} (t,s) \;ds
\end{equation}
is absolutely convergent in a tubular neighborhood of the  
curve joining $p_d$ and $p_{d+1}$ in $\Ma.$

This means that we can indeed perform the analytic continuation
using the above integral, so that we can cross over to the 
convergence region $U_{d+1}$ around the point $p_{d+1}.$ There 
we can close the contour of integration of the integral
(\ref{int:barnes}) (or (\ref{integral1})) to the left.
As shown above, any residues of the function 
$\Gamma(-s) \Gamma(1+s)$ can be discarded, so 
we only have to account for the poles of the functions 
$\Gamma (-h_j (s+\mu) - l_j^\prime - \mu_j^\prime),$ 
$j \in I_- \; .$ Hence, for $l^\prime$ fixed and $|x|$ large,
the last line in the formula (\ref{form1}) becomes
\begin{equation}\label{form3}
\sum_{m \leq m_2} \hbox{Res}_{\nu}  
\frac{ \pi e^{i \pi \epsilon \mu}}
{\sin(\pi(m+ \nu - \mu))}
\frac{\prod_{j \in I_-} \Gamma (-h_j (m+\nu) - l_j^\prime - \mu_j^\prime)}
{\prod_{j \notin I_-}  \Gamma (h_j (m+\nu) + l_j^\prime + \mu_j^\prime+1)}
\ (e^{-i\pi \epsilon} x)^{m+\nu}.
\end{equation}

We now analyze the case when the assumption {\bf \large (\dag)} does not hold, 
so there exist poles of the function $\Gamma(s) \Gamma(1-s)$ on the left
of the contour in the integral (\ref{int:barnes}). As shown above, we 
only have to account for the poles $s \geq m_1.$ For fixed
$l^\prime,$ only a finite number of them can be on the
left of the contour. The function obtained after
evaluating the residues at
these poles is merely a polynomial in $x$ (of degree depending on
$l^\prime$ though), and this means that the corresponding piece of 
$\Phi_\lambda^{\Ca_d}(t)$ can be manipulated freely in the $x$--direction, 
that is in a ``tubular'' neighborhood of the rational curve 
joining $p_d$ and $p_{d+1}$ in $\Ma.$ This shows that the
residues at all the poles of the function $\Gamma(s) \Gamma(1-s),$ whether 
located on the right or on the left of the contour, can be added 
together to obtain the series $\Phi_\lambda^{\Ca_d}(t).$

It follows that the analytic continuation of the series
$\Pa_d \cdot \Phi_{\lambda}^{\Ca_d}$ is given by
\begin{equation}
\begin{split}
& \Pa_d \cdot  \hbox{Res}_{\nu} \sum_{l^\prime \in \LL^\prime}  
\prod_{j \in I_-}
\frac{\sin(\pi (-h_j \mu - l_j^\prime - \mu_j^\prime))}{\pi}
\prod_{j=1}^N e^{2 \pi i t_j (l_j^\prime+\mu_j^\prime)} \cdot \\
&\cdot \sum_{m \leq m_2} 
\frac{ \pi e^{i \pi \epsilon \mu}}
{\sin(\pi(m+ \nu - \mu))}
\frac{\prod_{j \in I_-} \Gamma (-h_j (m+\nu) - l_j^\prime - \mu_j^\prime)}
{\prod_{j \notin I_-}  \Gamma (h_j (m+\nu) + l_j^\prime + \mu_j^\prime+1)}
\ (e^{-i\pi \epsilon} x)^{m+\nu} =\\
&= \hbox{Res}_{\nu}\; 
 e^{i \pi \epsilon (\mu-\nu)} \ \frac{\pi}{\sin(\pi(\nu -\mu))} 
\prod_{j \in I_-}
\frac{\sin(\pi (-h_j \mu -  \mu_j^\prime))}
{\sin(\pi (-h_j \nu  -  \mu_j^\prime))}
\cdot\\
&\cdot 
\sum_{l^\prime \in \LL^\prime}
 \sum_{m \leq m_2} \ \prod_{j=1}^N 
\frac{e^{2 \pi i t_j (l_j^\prime+\mu_j^\prime)}}
{\Gamma (h_j (m+\nu) + l_j^\prime + \mu_j^\prime+1 )}
\; x^{m+\nu} .
\end{split}
\end{equation}
The symbol $\Rz_\nu$ will be explained in detail below.
The characterization (\ref{lconecond22}) of the 
elements $l \in \Ca_{d+1}^\vee$ leads us to conclude that
the analytic continuation of $\Pi \cdot \Phi_{\lambda}^{\Ca_0}$ 
is $\Pi \cdot \Pa_{d+1} \cdot \Phi_{\lambda}^{\Ca_{d+1}}$
with $\Pa_{d+1}$ defined by 
\begin{eqnarray}
\Pa_{d+1} \cdot \Phi_{\lambda}^{\Ca_{d+1}}:=\Pa_d \cdot
 \hbox{Res}_{\nu} 
\frac{\pi  e^{i \pi \epsilon (\mu-\nu)}} {\sin(\pi(\nu -\mu))} 
\prod_{j \in I_-}
\frac{\sin(\pi (-h_j \mu -  \mu_j^\prime))}
{\sin(\pi (-h_j \nu  -  \mu_j^\prime))} \Phi_{\lambda}^{\Ca_{d+1}},
\nonumber \\ \label{pj1}
\end{eqnarray}
where the choice of the generators $\nu, \mu_j^\prime$ is 
determined by the combinatorial structure of
the transition from the cone $\Ca_d$ to the cone $\Ca_{d+1}$
(these are just different choices of coordinates on the
linear subspace $H \subset \CC^N$). In what follows we
show that this is indeed the analytic continuation, which
means that $\Pa_{d+1} \cdot \Phi_{\lambda}^{\Ca_{d+1}}$ gives
the solutions of the GKZ system in a neighborhood
of the point $p_{d+1}$. 


The operators $\Pa_d$ will not depend on the various
choices of coordinates involved. They
describe the analytic continuation procedure for the solutions
to the holonomic GKZ system, and the process is unique given the
path defined by the sequence of adjacent 
cones in the secondary fan. Such a path always exists, because
it is a well known fact that the singularities of the GKZ
system have real codimension $2$.
One needs to analytically continue along paths 
which are consistent with the the branch cuts mentioned above. In terms
of the $t$--coordinates, these branch cuts impose conditions
on the real part of the parameters $t.$

Care is required in defining the symbol Res$_\nu$ used in the
expressions (\ref{form3}). A direct calculation shows that
\begin{equation}
\begin{split}
& e^{i \pi \epsilon (\mu-\nu)} 
\frac{\pi}{\sin(\pi(\nu -\mu))} \prod_{j \in I_-}
\frac{\sin(\pi (-h_j \mu -  \mu_j^\prime))}
{\sin(\pi (-h_j \nu  -  \mu_j^\prime))}  \\
&= e^{2 \pi i \epsilon^\prime (\mu-\nu)} 
\frac{2\pi i}{e^{2 \pi i (\nu-\mu)}  - 1} 
\prod_{j \in I_-}
\frac{1- e^{2 \pi i (h_j \mu +  \mu_j^\prime)}}
{1-e^{2 \pi i (h_j \nu  +  \mu_j^\prime)}},
\label{form4}
\end{split}
\end{equation}
where $\epsilon^\prime$ is the integer defined by
\begin{equation}
\epsilon^\prime:= \frac{\epsilon +1 - \sum_{j \in I_-} h_j}{2}.
\end{equation}
The poles that we have to consider are the poles
of the function given by (\ref{form4}) other than $\nu=\mu$, 
such that we avoid over-counting in (\ref{form3}),
that is the complex values $\nu \in A_-(\mu^\prime),$ $\nu \not= \mu$ with
\begin{equation}\label{poles}
A_-(\mu^\prime):= \{\nu : (-h_j \nu  - \mu_j^\prime) \in  
\{ 0,1,\ldots,-h_j-1 \} \
\hbox{for some $j \in I_-$} \}.
\end{equation}
The fact that, according to conditions (\ref{lconecond21}) and 
(\ref{lconecond22}), $m=m_2$ is the largest integer such that
$l=-h_j m - l_j^\prime \leq 0,$ for some $j \in I_-$, implies that by taking 
residues at $\nu$ with $\nu \in A_-(\mu^\prime),$ and summing 
over all integers 
$m \leq m_2$ we have taken care of all the required poles of the functions
involved in (\ref{form3}) without over-counting. 

In order to deal with the issue of residues, we need the following
lemma.
\begin{lemma}\label{hai} 
For any function
$\phi \in \Ma_S(H),$ we have that
\begin{eqnarray}\label{hai2}
& &\Rz_{\nu} 
\frac{\pi}{\sin(\pi(\nu -\mu))} \prod_{j \in I_-}
\frac{\sin(\pi (-h_j \mu -  \mu_j^\prime))}
{\sin(\pi (-h_j \nu  -  \mu_j^\prime))} \; \phi (\nu, \mu^\prime)
\nonumber \\ \label{form10}
&=&\prod_{j \in I_-}
\sin(\pi(-h_j \mu - \mu_j^\prime)) \;
\widetilde{\phi}(\mu, \mu^\prime) - \phi(\mu, \mu^\prime),
\end{eqnarray}
for some function $\widetilde{\phi} \in \Ma_S(H).$ 
\end{lemma}
The underlying message of this lemma is that the sum over
the residues ``around'' the origin of $H$ (which is an integration
over some appropriate contour, as we will see below) has the 
feature of ``curing'' the poles of the integrand around the origin.
Fundamentally, this is a consequence of the classical Cauchy's theorem.

\begin{proof} (of the lemma) The main point of this lemma 
is to show that the sum of residues $\Rz_\nu$ defined as above
is an analytic function in $(\mu,\mu_j^\prime)$ in a 
neighborhood of the origin in $H.$ In order to see this,
assume first that the complex numbers $\mu_j^\prime$ are fixed and close
to zero. The fact that $\phi \in \Ma_H(S)$ implies that
the values of $\nu$ where the function $\phi (\nu, \mu^\prime)$
has residues (for fixed $\mu^\prime$) are non--zero, so
we can choose $\mu$ small in an open neighborhood of the origin such that
$\phi (\nu, \mu^\prime)$ is analytic in that open set. 
Assume first that $\mu \not\in A_-(\mu^\prime).$ That means that
the sum of residues $\Rz_\nu$ can be replaced by a contour of
integration $C_1$ which leaves outside the value $\mu,$ but encloses all 
the required poles in $A_-(\mu^\prime).$ If the function $\phi$ has poles 
given by the affine form $-h_j \nu  -  \mu_j^\prime +k , k\not= 0,$ 
for some $j \in I_-$, the contour will enclose those that are in 
$A_-(\mu^\prime),$ but will avoid all the others. We can write:

\begin{equation}\label{residue1}
\begin{split}
&  2 \pi i \; \Rz_{\nu} 
\frac{\pi}{\sin(\pi(\nu -\mu))} \prod_{j \in I_-}
\frac{\sin(\pi (-h_j \mu -  \mu_j^\prime))}
{\sin(\pi (-h_j \nu  -  \mu_j^\prime))} \; \phi (\nu, \mu^\prime) \\
&=\int_{C_1}
\frac{\pi}{\sin(\pi(\nu -\mu))} \prod_{j \in I_-}
\frac{\sin(\pi (-h_j \mu -  \mu_j^\prime))}
{\sin(\pi (-h_j \nu  -  \mu_j^\prime))} \; \phi (\nu, \mu^\prime)\; d\nu\\ 
&=
\prod_{j \in I_-}\sin(\pi (-h_j \mu -  \mu_j^\prime))
\int_{C_2}
\frac{\pi}{\sin(\pi(\nu -\mu))} \prod_{j \in I_-}
\frac{\phi (\nu, \mu^\prime)}
{\sin(\pi (-h_j \nu  -  \mu_j^\prime))} d\nu \\
&- 2 \pi i \; \Rz_{\nu=\mu} \frac{\pi}{\sin(\pi(\nu -\mu))} \prod_{j \in I_-}
\frac{\sin(\pi (-h_j \mu -  \mu_j^\prime))}
{\sin(\pi (-h_j \nu  -  \mu_j^\prime))} \; 
\phi (\nu, \mu^\prime),
\end{split}
\end{equation}
where the contour $C_2$ is chosen such that it encloses
the contour $C_1$ {\it and} the value $\nu=\mu,$ but no
other poles of the functions involved. 
The fact that some values in the set
$A_-(\mu^\prime)$ might coincide among themselves, or with
$\mu,$ does not change the analyticity of the function
defined by the integral over the contour $C_2.$ In other words,
the function of $\mu$ defined for fixed $\mu^\prime$
by the integral over $C_2$ is analytic also for 
$\mu \in A_-(\mu^\prime).$ Since 
for small values of $\mu^\prime_j$ the point $\nu=\mu$ 
is not a residue of $\phi(\nu,\mu^\prime_j),$ the last 
term in the formula (\ref{residue1}) is simply equal to 
$\phi(\mu,\mu^\prime).$ Recall that we fixed the complex
numbers $\mu^\prime_j.$ However small variations of them
around the origin have no effect on the contours of 
integration involved.
This shows indeed that the equality
(\ref{hai2}) holds in a neighborhood of the origin in $H,$
and $\widetilde{\phi}(\mu, \mu^\prime)$ defined by
\begin{equation}
\widetilde{\phi}(\mu, \mu^\prime):=\frac{1}{2\pi i}
\int_{C_2}
\frac{\pi}{\sin(\pi(\nu -\mu))} \prod_{j \in I_-}
\frac{1}
{\sin(\pi (-h_j \nu  -  \mu_j^\prime))} \; \phi (\nu, \mu^\prime)\; d\nu
\end{equation}
is analytic around the origin.
The integral over $C_1$ can be analytically 
continued by an appropriate continuous ``movement'' of the contour.
However, this procedure will fail for certain values $\mu, 
\mu_j^\prime.$ Cauchy's theorem shows that the function
of $\mu$ and $\mu^\prime$ defined by integration over the contour $C_1$
will have again poles of finite order, and, combinatorially
we see that these poles are included in the hyperplanes of
the form $\ker(\alpha +k), \alpha \in S$ for non--zero integers $k.$
\end{proof}

\

We have constructed the operator 
$\Pa_{d+1}$  which has all the required properties with  
the exception of the first property in the 
statement of the theorem. We turn our attention
to this question next. Assume that $J \subset \{ 1,\ldots,N \}$  
is saturated and does not define a simplex in the triangulation 
$\Ta_{d+1}.$

A possible case is that $J$ defines a simplex in the triangulation
$\Ta_d.$ We know that the change in simplices between $\Ta_d$ and
$\Ta_{d+1}$ is given by the circuit $I_+ \cup I_- , $ and the 
only way in which the simplex $J$ can be in $\Ta_d$ but not 
in $\Ta_{d+1}$ is if it contains a simplex of type
$(I_+ \cup I_-) \setminus \{ i \},$ for some $i \in I_+$. But
$I_- \subset (I_+ \cup I_-) \setminus \{ i \},$ 
so the simplex $J$ contains the simplex $I_-$. The definition (\ref{pj1}) of
the operator $\Pa_{d+1}$ shows that
\begin{equation}
\begin{split}
& \Pa_{d+1} \cdot \prod_{j \in J} \sin(\pi \lambda_j) \phi(\lambda)
= \Pa_d \cdot  \hbox{Res}_{\nu} 
\frac{\pi  e^{i \pi \epsilon (\mu-\nu)}} {\sin(\pi(\nu -\mu))} \\
& \prod_{j \in I_-}
\frac{\sin(\pi (-h_j \mu -  \mu_j^\prime))}
{\sin(\pi (-h_j \nu  -  \mu_j^\prime))}
\prod_{j \in J} \sin(\pi(-h_j \nu  -  \mu_j^\prime)) \phi(\lambda),
\label{fort}
\end{split}
\end{equation}
with $\lambda$'s written in the appropriate coordinates.
But $I _- \subset J,$ so all the poles in the formula (\ref{fort})
are canceled, and the result is zero. Note that in the analytic 
continuation procedure terms that are canceled immediately in this way 
show up when we perform a transition from $\Ca_{d+1}$ back to 
$\Ca_d.$

Assume now that $J$ does not define a simplex neither in
$\Ta_{d+1},$ nor in $\Ta_d$. By using lemma \ref{hai} we can write 
that 
\begin{equation}
\begin{split}
\Pa_{d+1} \cdot \prod_{j \in J} \sin(\pi \lambda_j) \phi(\lambda) 
& = \Pa_d \cdot (\prod_{j \in I_-}
\sin(\pi(-h_j \mu - \mu_j^\prime)) \;
\widetilde{\phi}(\mu, \mu^\prime) - \\
& -\prod_{i\in J} \sin(\pi(-h_j \mu -  \mu_j^\prime)) \phi(\lambda) ),
\label{form12}
\end{split}
\end{equation}
for some function $\widetilde{\phi} \in \Ma_S(H).$ The second 
term in the formula (\ref{form12}) is canceled by $\Pa_d$ according
to the induction hypothesis, since $J$ is saturated and it is 
not a simplex in $\Ta_d.$
We show now that in fact $\Pa_d$ cancels the first term in the 
formula (\ref{form12}), too.

The two corresponding cones $\Ca_d$ and $\Ca_{d+1}$ in 
$\LL_{\RR}^\vee$ have a common face which is a cone $\Ca_{d,d+1}$ of 
dimension $N-n-1,$ and let's denote by $I_0$ the set of vectors of the Gale 
transform included in this cone and by $H_{d,d+1}$ the
hyperplane in $\LL_{\RR}^\vee$ containing the cone $\Ca_{d,d+1}.$

Assume first that for any element $j \in J$ the corresponding
element in the Gale transform $w_j$ is not contained in the
hyperplane $H_{d,d+1}.$ This means that
$J \subset I_+ \cup I_-$, and, since  $J$ is not a simplex
in either one of $\Ta_d$ and $\Ta_{d+1},$ the only possibility
is that $J$ contains $I_+ \cup I_-$. But now we are back to the
previous case, and formula (\ref{fort}) proves again the 
required statement.

Assume now that $J$ contains an element $j$ 
such that the corresponding vector of the Gale 
transform has the property $w_j \in H_{d,d+1}.$ 
Moreover, assume that $w_j$ is not contained in any cone of type
$\Ca_{K^*}$ (see proposition \ref{prop:cone}), for 
$K \in \Ta_d^n.$ Equivalently
this means that $j$ belongs to any simplex of the triangulation
$\Ta_d$ (otherwise, there would be a simplex $K \in \Ta_d^n$
such that $j \in K^*,$ that is $\Ca_d \subset \Ca_{K^*}$
which contradicts the assumption). So the element $j$
will not be missing from any of the simplices of $J,$ 
and because $J$ is not a simplex in $\Ta_d$ we obtain
again that $J$ contains $I_-$ so again we can apply formula 
(\ref{fort}) to prove the desired statement. 

The condition \ref{cond1} that was imposed in the beginning 
of this section shows that
the last case we have to consider is the case when
$J$ contains an element $k \in I_0$ such that $w_k$ is 
a generator of the cone $\Ca_{d,d+1}$ but there 
is no other larger cone in the hyperplane $H_{d,d+1}$
that includes the cone $\Ca_{d,d+1}$ as a proper subset.
Denote by $I_k$ the set of elements $j \in I_0$
such that $w_j= t w_k,$ for some positive real number $t >0.$
The point here is then that for any simplex $K \in \Ta_d^n,$ 
(which means $\Ca_d \subset \Ca_{K^*} \subset \LL^\vee_\RR$)
the cone $ \Ca_{K^*}$ has to contain at least one of the vectors 
$w_j$ for some $j \in I_- \cup I_k$ (equivalently 
$(I_- \cup I_k) \cap K^* \not= \emptyset$).
Hence $I_- \cup I_k $ is not included in any 
maximal dimensional simplex of $\Ta_d,$ so it cannot be a simplex 
in $\Ta_d.$ Now we can apply lemma \ref{hai}. In particular,
for any $j \in I_k \subset I_0$ the element $\sin(\pi \lambda_j)= 
\sin(\mu_j^\prime)$ is independent of $\nu$ (the 
integration parameter), so it remains unchanged after integration and 
we can write that

\begin{equation}
\begin{split}
\Pa_{d+1} \cdot \prod_{j \in J} \sin(\pi \lambda_j) \phi(\lambda)
& = \Pa_d \cdot (\prod_{j \in I_- \cup I_k }
\sin(\pi(-h_j \mu - \mu_j^\prime)) \;
\widetilde{\phi}(\mu, \mu^\prime) - \\
& - \prod_{i\in J} \sin(\pi(-h_j \mu -  \mu_j^\prime))\phi(\lambda) ),
\label{form13}
\end{split}
\end{equation}
for some function $\widetilde{\phi} \in \Ma_S(H).$ Again, because
neither $I_- \cup I_k$ nor $J$ are simplices in $\Ta_d,$ and they
are saturated, the induction hypothesis implies that $\Pa_d$
cancels both terms in the formula (\ref{form13}).
\end{proof}

The inductive analytic continuation formula obtained in the previous
proof is contained in the following corollary
(see (\ref{pj1}) and (\ref{form4})).

\begin{corollary}\label{cor:contform}
The analytic continuation  $\Pa_{D+1} \cdot \Phi_\lambda^{\Ca_{D+1}}$
(the solution to the GKZ system in a neighborhood of $p_{D+1}$)
of the series $\Pa_D \cdot \Phi_\lambda^{\Ca_D}$ (the solution
to the GKZ system in a neighborhood of $p_D$) is given by 
\begin{equation}\label{form133}
\Pa_D \cdot
\Rz_\nu e^{2 \pi i \epsilon^\prime (\mu-\nu)} 
\frac{2\pi i}{e^{2 \pi i (\nu-\mu)}  - 1} 
\prod_{j \in I_-}
\frac{1- e^{2 \pi i (h_j \mu +  \mu_j^\prime)}}
{1-e^{2 \pi i (h_j \nu  +  \mu_j^\prime)}} \; 
\Phi_{\nu,\mu'}^{\Ca_{D+1}},
\end{equation}
where the integer $\epsilon^\prime$ is given by 
\begin{equation}\epsilon':=
\begin{cases}
1+ (L/2) & 
\hbox{if}\ L  \ \hbox{is even},\\
(1+L)/2 &
\hbox{if}\ L \ \hbox{is odd},
\end{cases} \label{def:eps1}
\end{equation}
with $L:= \sum_{j \in I_+}h_j = - \sum_{j \in I_-}h_j.$ The symbol
$\Rz_\nu$ represents the sum over the poles $\nu$
with $\nu \in A_-(\mu^\prime),$ and
$\nu \not= \mu.$ The finite set $A_-(\mu^\prime)$
is defined by
\begin{equation}\label{def:A-0}
A_-(\mu^\prime):= \{\nu: (-h_j \nu  - \mu_j^\prime) \in  
\{ 0,1,\ldots,-h_j-1 \} \
\hbox{for some $j\in I_-$} \}.
\end{equation}

\end{corollary}

We are now ready to apply the methods developed in the previous 
proof to compute the monodromy transformations. The main inspiration
for analyzing the transitions between different cones in the 
secondary fan comes from string theory. The works of Witten \cite{Witten1}
and Aspinwall, Greene and Morrison 
\cite{Paul},\cite{AGM1},\cite{AGM2},\cite{AGM3} showed 
that these transitions are related to the so--called ``phase
transitions'' in string theory. From the physics point of view,
the smooth Calabi--Yau manifolds (like those corresponding to unimodular 
triangulations) are viewed on equal footing with the 
singular phases (corresponding to non--unimodular triangulations).
In fact, the study of the phase in which the whole Calabi--Yau space
collapses to a point (the so--called Landau--Ginzburg orbifold point) 
inspired Greene and Plesser \cite{GrPl} to propose the first mirror symmetry
construction, also known as the orbifolding construction. Our 
discussion is closest in style to the methods pursued in 
\cite{Paul} and \cite{AGM3}. The paths to be considered are paths
in the toric variety determined by the secondary fan. In the terminology
introduced by Aspinwall, Greene and Morrison,
this toric variety compactifies the space of complexified 
K\"ahler forms of the Calabi--Yau complete intersection $W$
and identifies it with the moduli space of complex structures
of the mirror Calabi--Yau space. This allows a 
detailed study of the discriminant locus (the singular complex
structures locus, and the same time the singular locus of the
GKZ system) using the general methods developed by Gelfand, 
Kapranov and Zelevinsky \cite{GKZ2} who showed that this locus 
has generically complex codimension one and admits a stratification whose 
combinatorial structure is incorporated 
in the data provided by the set $\cA$ (see \cite{GKZ2}). The
so--called ``principal component'' of the discriminant locus and 
its parameterization (the Horn uniformization) will later play a 
distinguished r\^ole in our analysis. The loops that we will 
consider are similar to the loops studied in \cite{AGM3}, namely
they can be continuously deformed to loops included in the rational curves 
in the moduli space connecting limit points determined by
adjacent cones in the secondary fan. The web
of such rational curves gives the secondary polytope (see Figure 3
in \cite{AGM3} for a nice example). The main difference compared to
the loops considered in \cite{AGM3} is in the fact that the loops
that we work with need 
{\it not} be included in the rational curves. This is a natural
feature since we are interested in the monodromy of the 
full GKZ system, not only on its restriction to the 
rational curves. 

More precisely, if we denote by $x$ the coordinate 
on a rational curve joining two limit points (as we did in the
previous proof), then
we can assume that the other components of the 
discriminant intersect this curve at $x=0, \infty$ and at some
other non--zero finite value $x=a.$ The exact value of $a$
will not be important to us, although Gelfand, Kapranov and
Zelevinsky (\cite{GKZ2}) developed a general theory of discriminants which 
provides a combinatorial
interpretation of this value. 
One can rescale the coordinate $x$ such that this value 
becomes $a=1$ (see \cite{AGM3} for a discussion on this, and the
relationship with the monomial divisor mirror map). For 
each transition between two adjacent cones, the loop 
that we will consider will go around the components of
the discriminant locus passing through $x=a.$ Since none of
the power series that we use converges around $x=a,$
we will in fact obtain an analytic continuation formula
along such a loop by considering a composition of loops
around the points $x=0$ and $x=\infty$ (see figure \ref{fig:cont2}).

The sequence of transitions that we will study is
of the type $\Ca_D \to \Ca_{D+1} \to \Ca_D.$ According
the the previous result, the analytic
continuation of the series $\Phi_{\lambda}^{\Ca_0}$ to the 
neighborhood of the point $p_D$ is given by $\Pa_D \cdot
\Phi_{\lambda}^{\Ca_D},$ and the
corollary \ref{cor:contform}. gives that further analytic 
continuation to the neighborhood of 
the point $p_{D+1}$ is given by 

\begin{equation}\label{form20}
\Pa_D \cdot
 \hbox{Res}_{\nu} e^{2 \pi i \epsilon^\prime (\mu-\nu)} 
\frac{2\pi i}{e^{2 \pi i (\nu-\mu)}  - 1} 
\prod_{j \in I_-}
\frac{1- e^{2 \pi i (h_j \mu +  \mu_j^\prime)}}
{1-e^{2 \pi i (h_j \nu  +  \mu_j^\prime)}} \; 
\Phi_{\nu,\mu'}^{\Ca_{D+1}},
\end{equation}
where we write $\Phi_{\nu,\mu'}^{\Ca_{D+1}}$ for 
$\Phi_{\lambda}^{\Ca_{D+1}}$ to show explicitly that 
the parameter $\lambda$ is written in the coordinates $\nu, \mu^\prime.$
The integer $\epsilon^\prime$ is defined by 
\begin{equation}
\epsilon^\prime:= \frac{\epsilon +1 - \sum_{j \in I_-} h_j}{2}.
\end{equation}
The analytic continuation corresponding to the transition 
from $\Ca_{D+1}$ back to $\Ca_D$ will be then given by 
\begin{eqnarray}
\Pa_D &\cdot&
\hbox{Res}_{\nu} e^{2 \pi i \epsilon^\prime (\mu-\nu)} 
\frac{2 \pi i}{e^{2 \pi i (\nu-\mu)}  - 1} 
\prod_{j \in I_-}
\frac{1- e^{2 \pi i (h_j \mu +  \mu_j^\prime)}}
{1-e^{2 \pi i (h_j \nu  +  \mu_j^\prime)}} \nonumber \\
& &\hbox{Res}_{\xi} e^{2 \pi i \epsilon'' (\nu-\xi)} 
\frac{2 \pi i }{e^{2 \pi i (\xi-\nu)}  - 1} 
\prod_{j \in I_+}
\frac{1- e^{2 \pi i (h_j \nu +  \mu_j^\prime)}}
{1-e^{2 \pi i (h_j \xi  +  \mu_j^\prime)}} \; 
\Phi_{\xi, \mu^\prime}^{\Ca_D}\nonumber  \; , \\
\label{form21}   
\end{eqnarray}
where $\epsilon''$ is the integer defined by 
\begin{equation}
\epsilon'':= \frac{\epsilon +1 - \sum_{j \in I_+} h_j}{2}
= \frac{\epsilon +1 + \sum_{j \in I_-} h_j}{2}.
\end{equation}
Of course we know that the expression (\ref{form21})
is nothing else than $\Pa_D \cdot \Phi_{\lambda}^{\Ca_D},$ 
since the loop followed in the rational curve joining the 
points $p_D$ and $p_{D+1}$ is trivial. Let's analyze this formula 
in more detail. As we have seen, the two sums of residues 
can be expressed as contour integrals, with specially 
chosen contours, say $C_\nu$ and $C_\xi,$ respectively.  
The previous formula can be written as
\begin{equation}\label{form22}
\begin{split}
\Pa_D \cdot \int_{C_\nu} \int_{C_\xi} 
e^{2 \pi i \epsilon^\prime (\mu-\nu)} 
\frac{1}{e^{2 \pi i (\nu-\mu)}  - 1} 
\prod_{j \in I_-}
\frac{1- e^{2 \pi i (h_j \mu +  \mu_j^\prime)}}
{1-e^{2 \pi i (h_j \nu  +  \mu_j^\prime)}}\\
 e^{2 \pi i \epsilon'' (\nu-\xi)} 
\frac{1}{e^{2 \pi i (\xi-\nu)}  - 1} 
\prod_{j \in I_+}
\frac{1- e^{2 \pi i (h_j \nu +  \mu_j^\prime)}}
{1-e^{2 \pi i (h_j \xi  +  \mu_j^\prime)}} \; 
\Phi_{\xi,\mu^\prime }^{\Ca_D} \; d\xi d\nu.
\end{split}
\end{equation}
Again we need to exercise care in choosing the contours
$C_\nu$ and $C_\xi.$ The contour $C_\nu$ is initially chosen 
such that it encloses only the poles included in the set 
$A_-(\mu^\prime)$ as defined in (\ref{poles}) and 
corresponding to the integers $h_j,$ $j \in I_-$, and 
does not include the value $\nu=\mu.$ It is not hard
to see that nothing is lost if the contour $C_\nu$ does
not inclose those values $\nu$ that belong to 
$A_-(\mu^\prime) \cap A_+(\mu^\prime),$ where $A_+(\mu^\prime)$ is
the analogous set corresponding to the integers $h_j, j \in I_+$. Indeed,
formula (\ref{form22}) shows that those poles are canceled anyway.
Then the contour $C_\xi$ can be chosen so that it encloses all the 
poles $\xi$ included
in the set $A_+(\mu^\prime),$ and lies outside the 
contour $C_\nu.$ The principal merit of this choice of contours
is that $C_\xi$ is not dependent on the parameter $\nu,$ at least
for values $\nu$ in an appropriate domain. This allows us
to change the order of integration in formula (\ref{form22})
which can rewritten as
\begin{equation}\label{form23}
\begin{split}
 \Pa_D  \cdot  \int_{C_\xi} & \int_{C_\nu} 
 e^{2 \pi i \epsilon^\prime (\mu-\nu)} 
\frac{1}{(e^{2 \pi i (\nu-\mu)}  - 1)(e^{2 \pi i (\xi-\nu)}  -1) } 
\frac{\prod_{j \in I_+} (1- e^{2 \pi i (h_j \nu +  \mu_j^\prime)})}
{\prod_{j \in I_-} (1- e^{2 \pi i (h_j \nu +  \mu_j^\prime)})} \\
& e^{2 \pi i \epsilon'' (\nu-\xi)} 
\frac{\prod_{j \in I_-} (1- e^{2 \pi i (h_j \mu +  \mu_j^\prime)})}
{\prod_{j \in I_+} (1-e^{2 \pi i (h_j \xi  +  \mu_j^\prime)})} \; 
\Phi_{\xi, \mu^\prime}^{\Ca_D} \; d\nu d\xi.
\end{split}
\end{equation}
At this point, one can see that the contour $C_\nu$ can be 
``moved'' again to enclose all the values $\nu \in A_-(\mu^\prime)$
without affecting the final result. We investigate the contour integral
\begin{equation}\label{form24}
\begin{split}
\int_{C_\nu} e^{2\pi i (\epsilon''-\epsilon')\nu} 
\frac{1}{(e^{2 \pi i (\nu-\mu)}  - 1) 
(e^{2 \pi i (\xi-\nu)}  - 1)} 
\frac{\prod_{j \in I_+} (1- e^{2 \pi i (h_j \nu +  \mu_j^\prime)})}
{\prod_{j \in I_-} (1- e^{2 \pi i (h_j \nu +  \mu_j^\prime)})} \; d\nu.\\
\end{split}
\end{equation}
Note that 
\begin{equation}\label{def:p}
\epsilon''-\epsilon'= \sum_{j \in I_-} h_j= - \sum_{j \in I_+} h_j < 0
\end{equation}
After applying the change of variable $t=e^{2 \pi i \nu},$ the integral
(\ref{form24}) becomes 
\begin{equation}\label{form25}
\begin{split}
& \int_{C_t} t^{\epsilon''-\epsilon} \ 
\frac{1}{e^{-2 \pi i\mu} \; t -1} \; 
\frac{t}{e^{2 \pi i\xi} -t}\;
\frac{\prod_{j \in I_+} (1 -  e^{2 \pi i \mu_j^\prime} \; t^{h_j}  )}
{\prod_{j \in I_-} (1 -  e^{2 \pi i \mu_j^\prime} \; t^{h_j}  )}
\; \frac{dt}{2\pi i t}\\
= & \int_{C_t} \frac{1}{e^{-2 \pi i\mu} \; t -1} \;
\frac{t}{e^{2 \pi i\xi} -t}\;
\frac{\prod_{j \in I_+} (1 -  e^{2 \pi i \mu_j^\prime} \; t^{h_j}  )}
{\prod_{j \in I_-} ( t^{-h_j}-  e^{2 \pi i \mu_j^\prime}  )}
\; \frac{dt}{2\pi i t},
\end{split}
\end{equation}
where the contour $C_t$ is chosen so that it encloses all the poles
in the $t$--plane of the functions in the integrand, {\it with the
exception} of the possible poles 
$t=0, \infty, e^{2\pi i \mu}, e^{2\pi i \xi}.$ Note the that $t$--plane is 
very well adapted for such an
integration, and, ultimately this is a consequence of the fact
that we had to avoid over-counting when integrating along the 
contour $C_\nu.$ By direct power counting, we see 
that $t=0$ and $t=\infty$ are not poles after all. The residue 
at the point $t= e^{2\pi i \xi}$ cancels out all the poles
of the integrand in (\ref{form23}) (with respect to $\xi$). It follows
that the only pole that needs be taken into account is $t=e^{2\pi i \mu}.$
Recall though that formula (\ref{form23}) is just another way of writing
$\Pa_D \cdot \Phi_{\xi, \mu^\prime}^{\Ca_D},$ since
we have performed analytic continuation along a homotopically 
trivial loop, and as a consequence of the present remarks 
this result is obtained solely as a contribution
of the pole $t=e^{2\pi i \mu}$ in the integral (\ref{form25}).

\iffigs
\begin{figure}
  \centerline{\epsfxsize=7cm\epsfbox{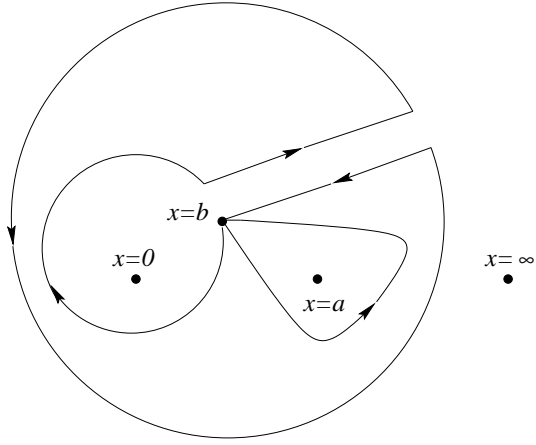}}
  \caption{The loop based at $x=b$ in the $x$--plane along which 
the analytic continuation is performed. The special values are
$x=0,$ $x=a$ and $x=\infty.$}
  \label{fig:cont2}
\end{figure}
\fi

We are now ready to consider a non--trivial loop in the $x$--plane 
(the rational curve joining $p_D$ and $p_{D+1}$) as shown in 
figure \ref{fig:cont2}. The effect of the loops around $x=0$ and $x=\infty$
on the the formula (\ref{form23}) is the appearance of 
two new factors in the integrand. The general form of the $\Phi$ series
shows that the clockwise
loop around $x=0$ introduces the factor $e^{-2 \pi i \mu},$
while the counterclockwise loop around $x=\infty$ brings in the 
factor $e^{2 \pi i \nu}.$ The direct analog of the 
integral (\ref{form25}) that we need to study in this case is 
\begin{equation}
\int_{C_t} \frac{t}{e^{2 \pi i \mu}} \;
\frac{1}{e^{-2 \pi i\mu} \; t -1} \;
\frac{t}{e^{2 \pi i\xi} -t}
\frac{\prod_{j \in I_+} (1 -  e^{2 \pi i \mu_j^\prime} \; t^{h_j}  )}
{\prod_{j \in I_-} ( t^{-h_j}-  e^{2 \pi i \mu_j^\prime}  )}
\; \frac{dt}{2\pi i t}.
\end{equation}
As before the poles $t=0$ and $t=e^{2\pi i \xi}$ are just
apparent singularities of the integrand, so they do not  contribute 
at all to the result. The residue
at $t=e^{ 2 \pi i \mu}$ gives exactly the same answer as in 
(\ref{form25}), which means that, as before, the contribution
of this residue to the analytic continuation of $\Pa_D \cdot
\Phi_{\mu, \mu^\prime}^{\Ca_D}$ is just  $\Pa_D \cdot
\Phi_{\mu, \mu^\prime}^{\Ca_D}.$
 
Finally, we are left with the pole $t=\infty.$ We obtain 
\begin{equation}
\begin{split}
& 2\pi i \; \Rz_{t=\infty} 
\frac{t}{e^{2 \pi i \mu}} \;
\frac{1}{e^{-2 \pi i\mu} \; t -1} \;
\frac{t}{e^{2 \pi i\xi} -t}
\frac{\prod_{j \in I_+} (1 -  e^{2 \pi i \mu_j^\prime} \; t^{h_j}  )}
{\prod_{j \in I_-} ( t^{-h_j}-  e^{2 \pi i \mu_j^\prime}  )}
\; \frac{1}{2\pi i t} \\ 
&=  2\pi i \; \Rz_{t=\infty} 
\frac{t}{t - e^{2 \pi i \mu} } \;
\frac{t}{e^{2 \pi i\xi} -t}
\frac{\prod_{j \in I_+} (1 -  e^{2 \pi i \mu_j^\prime} \; t^{h_j}  )}
{\prod_{j \in I_-} ( t^{-h_j}-  e^{2 \pi i \mu_j^\prime}  )}
\; \frac{1}{2\pi i t}  
\end{split}
\end{equation}\begin{equation}
\begin{split}
&= 2\pi i \;  \Rz_{s=0} 
\frac{1}{1 -  s e^{2 \pi i \mu} } \;
\frac{1}{s e^{2 \pi i\xi} -1}
\frac{\prod_{j \in I_+} (s^{h_j} -  e^{2 \pi i \mu_j^\prime}  )}
{\prod_{j \in I_-} ( 1 -  e^{2 \pi i \mu_j^\prime} \;s^{-h_j}  )}
\; \frac{1}{2\pi i s} \\
&= -  \prod_{j \in I_+} (-e^{2 \pi i \mu_j^\prime} ) .
\end{split}
\end{equation}
Using formula (\ref{form23}), we conclude that
the analytic continuation of  
the series $\Pa_D \cdot \Phi_{\lambda}^{\Ca_D}$ along the
counterclockwise loop around $x=a$ (see figure \ref{fig:cont2})
is given by 
\begin{equation}\label{form:mon0}
\begin{split}
\Pa_D \cdot \Phi_{\mu,\mu^\prime }^{\Ca_D} \ \mapsto \
& \Pa_D \cdot \Big( \Phi_{\mu,\mu^\prime}^{\Ca_D}
- e^{2 \pi i \epsilon' \mu} \prod_{j \in I_-}
(1- e^{2 \pi i (h_j \mu +\mu_j^\prime)})  \cdot \\
& \int_{C_\xi}e^{-2 \pi i \epsilon'' \xi}
 \prod_{j \in I_+} \frac{ - e^{2 \pi i \mu_j^\prime}  }
{1- e^{2 \pi i (h_j \xi + \mu_j^\prime)}}
\; \Phi_{\xi,\mu^\prime}^{\Ca_D} \; d\xi  \; \Big).
\end{split}
\end{equation}
Since $\epsilon''=\epsilon' - \sum_{j \in I_+} h_j,$ we see that
\begin{equation}
e^{2 \pi i \epsilon' \mu} e^{-2 \pi i \epsilon'' \xi}
\prod_{j \in I_+} \frac{- e^{2 \pi i \mu_j^\prime} }
{1- e^{2 \pi i (h_j \xi + \mu_j^\prime)}} 
= e^{2 \pi i \epsilon'(\mu- \xi)}
\prod_{j \in I_+} \frac {1}{1- e^{-2 \pi i (h_j \xi + \mu_j^\prime)}}.
\end{equation}
We now summarize the result of our calculation.

\begin{thm}\label{thm:mon}
The analytic continuation of the series $\Pa_D \cdot 
\Phi_{\lambda}^{\Ca_D}$ along a loop around 
the point $x=a$ (as shown in figure \ref{fig:cont2})
homotopic to a loop included in the rational curve joining the points 
$p_D$ and $p_{D+1},$ is given by the following formula
\begin{equation}\label{form:mon}
\begin{split}
\Pa_D \cdot \Phi_{\mu,\mu^\prime }^{\Ca_D} \ \mapsto \
& \Pa_D \cdot \Big( \Phi_{\mu,\mu^\prime}^{\Ca_D}
- \prod_{j \in I_-} (1- e^{2 \pi i (h_j \mu +\mu_j^\prime)})\cdot \\
& \int_{C_\xi} 
e^{2 \pi i \epsilon'(\mu- \xi)}
\prod_{j \in I_+} \frac{1}{1- e^{-2 \pi i (h_j \xi + \mu_j^\prime)}}
\; \Phi_{\xi,\mu^\prime}^{\Ca_D} \; d\xi \; \Big),
\end{split}
\end{equation}
where the integer $\epsilon'$ is given by 
\begin{equation}\epsilon':=
\begin{cases}
1+ (L/2) & 
\hbox{if}\ L  \ \hbox{is even},\\
(1+L)/2 &
\hbox{if}\ L \ \hbox{is odd},
\end{cases} \label{def:eps}
\end{equation}
with $L:= \sum_{j \in I_+}h_j = - \sum_{j \in I_-}h_j.$ The contour 
$C_\xi$ is chosen such that it encloses only the poles $\xi$
with $\xi \in A_+(\mu^\prime),$ where the finite set $A_+(\mu^\prime)$
is defined by
\begin{equation}\label{def:A+}
A_+(\mu^\prime):= \{\xi: (h_j \xi  + \mu_j^\prime) \in  
\{ 0,1,\ldots,h_j-1 \} \
\hbox{for some $j\in I_+$} \}. 
\end{equation}
\end{thm}


\section{Automorphisms of Derived Categories \\ and Monodromy
Transformations} \label{chap:4}

In this chapter we use the results of the previous sections 
to prove some cases of Kontsevich's conjecture. In fact,
in section \ref{chap:functors}, we state a refinement of Kontsevich's 
proposal (conjecture \ref{newconj}) specific to the toric context 
considered in this work.
We study a smooth Calabi--Yau complete intersection 
$W=E_1 \cap \ldots E_k \subset X_{\Sigma}=\PP_{\Delta},$ with 
$X_{\Sigma}=\PP_{\Delta}$ a toric Fano variety described by 
the fan $\Sigma \subset \Nrm_{\RR}$ (or, equivalently, by the reflexive 
polytope $\Delta \subset \Mrm_{\RR}$) and its mirror Calabi--Yau variety 
$M.$ We refer to remark \ref{rem:toric} for the detailed description of the 
combinatorics of the toric structure involved. We mention that the 
complete fan $\Sigma \subset \Nrm_{\RR}$ with $\Sigma(1)= 
\{ v_{k+1}, \ldots, v_N \},$ determines a fan 
$\overline{\Sigma} \subset \Nrm_{\RR} \times \RR^k$ with 
$\overline{\Sigma}(1)= \cA= \{ \ov_1, \ldots, \ov_k, \ov_{k+1}, 
\ldots, \ov_N \},$ and the vectors 
$\ov_1, \ldots, \ov_k$ corresponding to the nef--partition
$E_1, \ldots, E_k.$

As described in section \ref{chap:triang}, the phase transitions of the 
Calabi--Yau space $W$ are associated with modifications of the fan structure 
(triangulations) of the $\cA$-set given by circuits in $\cA.$ We need to 
impose a numerical condition for all the circuits that will be used 
in this chapter. 

For any positive integer $d \geq 1,$ define the set $A_d \subset \QQ,$ by
\begin{equation}\label{def:A}
A_d:= \{ \dfrac{\xi}{d} : 0 \leq \xi \leq d-1, \xi \in \ZZ \}.
\end{equation}

Given a circuit in $\cA$ of the form 
\begin{equation}
\sum_{j \in I_+} q_j \ov_j= \sum_{j \in I_-'} q_j \ov_j +
\sum_{j \in I_-''} d_j \ov_j,
\end{equation}
with $(I_+ \cup I_-') \subset \{k+1, \ldots, N\},$ 
$I_-'' \subset \{1, \ldots , k\},$ $q_j, d_j > 0,$ we assume that for 
any non--zero rational number $r \in \QQ,$ the following is true:

\begin{condition}\label{cond2}
\begin{equation}
{\rm Card}\{ j : r \in A_{q_j}, j \in I_+ \} \leq
{\rm Card}\{ j : r \in A_{d_j}, j \in I_-'' \}.
\end{equation}
\end{condition}
In other words, any non--zero rational number belongs to more sets of type
$A_{d_j},$ $j \in I_-'',$ than sets of type $A_{q_j},$ $j \in I_+.$
The reason for imposing  this condition will become clear
a bit later in our analysis. Essentially, it will allow us to reduce 
the study of the residues of the functions of type
\begin{equation}
\frac{\prod_{j\in I_-''} (1- e^{-d_j x})}
{\prod_{j\in I_+} (1- e^{-q_j x})}
\end{equation}
to the study of the residue at $x=0.$ We should note that this condition
is satisfied for all the examples considered in the literature
(for example, those in \cite{BatVS}, \cite{Yau}). In particular,
it is trivially satisfied for the case of a primitive relation
($q_j=1,$ for all $j \in I_+$). However, at the time of this 
writing, I am not able to see why such a condition would be true in the 
general case. If this is the case, then it probably has to do with
the smoothness of $W$ and with the fact that the divisors $E_j$ are 
Cartier-- for an illustration of how the latter fact
might be used, see the proof of the lemma 3.5.6 in \cite{CK}.

\subsection{Automorphisms of Derived Categories and 
\\Fourier--Mukai Functors}\label{chap:functors}

We review some notions about derived categories and their
automorphisms. For a detailed account of Verdier's formalism 
of derived categories, see \cite{verdier}, \cite{gelmanin}.
As mentioned in the introduction, in 
Kontsevich's point of view  \cite{Kont1},\cite{Kont2},
mirror symmetry should be viewed as an equivalence between 
``Fukaya's $A_\infty$ category'' of a Calabi--Yau manifold $M$
and the bounded derived category of coherent sheaves of the 
mirror Calabi--Yau manifold $W.$ Using the monodromy formulae discussed 
in the previous sections, we will be able to confirm some of the 
consequences of Kontsevich's conjecture. 

In this section, we present some constructions
of automorphisms of the bounded derived category of coherent sheaves 
of a smooth projective variety. For a smooth projective 
variety $W,$ $\Der (W)$ denotes the bounded derived category 
of coherent sheaves on $W.$ The category itself and its automorphisms
for a general smooth projective variety have been the subject of study 
in recent years, after Mukai \cite{Mukai1}, \cite{Mukai2} had studied 
the derived category of coherent sheaves for the case of abelian varieties
and $K3$ surfaces. Since our intention is to 
compare the action on the cohomology of a Calabi--Yau manifold $W$ of the 
automorphisms of its derived category of coherent sheaves, with the 
previously obtained action of the monodromy transformations on the 
cohomology of the mirror Calabi--Yau manifold $M,$ we will focus our 
attention on the group of automoprhisms of derived categories 
and Fourier--Mukai functors. 

Any object $\Ea^\bullet$ (also called a {\bf kernel})
in $\Der(W \times W)$ 
determines an exact functor $\Phi_{\Ea^\bullet}$ from $\Der(W)$ to $\Der(W)$ 
(called a {\bf Fourier--Mukai functor}) defined by 
\begin{equation}
\label{eq:auto1}
\Phi_{\Ea^\bullet}(\cdot):={\bf \rm R}p_{2_*}(\Ea^\bullet \Ltensor p_1^*(\cdot)).
\end{equation}
Bondal and Orlov \cite{BO} proved that the functor $\Phi_{\Ea^\bullet}$ is fully 
faithful if and only if the following orthogonality conditions are satisfied:
\begin{eqnarray}
& & \Hom ^i_X (\Phi_{\Ea^\bullet}(\Oa_{t_1}), \Phi_{\Ea^\bullet}(\Oa_{t_2})) =0,\
\hbox{for every $i$ and $t_1 \not= t_2,$}\\
& & \Hom ^0_X (\Phi_{\Ea^\bullet}(\Oa_{t}), \Phi_{\Ea^\bullet}(\Oa_{t})) =\CC,\\
& &  \Hom ^i_X (\Phi_{\Ea^\bullet}(\Oa_{t}), \Phi_{\Ea^\bullet}(\Oa_{t})) =0,\ 
\hbox{unless $0 \leq i \leq \dim W,$}
\end{eqnarray}
where $t, t_1, t_2$ are points in $W$, and $\Oa_{t_j}$ the 
corresponding skyscraper sheaves. Moreover, Bridgeland \cite{Bridge}
proved that $\Phi_{\Ea^\bullet}$ is in fact an equivalence of categories
precisely when  
\begin{equation}
\Phi_{\Ea^\bullet}(\Oa_{t}) \otimes \omega_W \cong \Phi_{\Ea^\bullet}(\Oa_{t}),
\end{equation}
for any point $t \in W,$ where $\omega_W$ is the canonical
sheaf on $W.$ Obviously, this condition is satisfied if $W$ is
Calabi--Yau. 

We now present some  examples of kernels and their corresponding 
actions in cohomology. According to expos\'{e} IV in
\cite{SGA}, there exists an additive map from $\Der(W)$ to the
complex algebraic Grothendieck group of coherent sheaves $K_\bullet(W)$ 
on $W.$ A map $f$  
is said to be {\bf additive} if $f(\Ea^\bullet)=f(\Ea^{\prime\bullet})+
f(\Ea^{\prime\prime\bullet}) $ for any distinguished triangle 
$\Ea^{\prime\bullet} \to \Ea^\bullet \to \Ea^{\prime\prime\bullet} \to
\Ea^{\prime\bullet}[1]$ 
in the triangulated 
category $\Der(W).$ Each complex in $\Der (W)$
has a canonical image in $K_\bullet(W),$ and each exact functor $\Der (W) 
\to \Der (W')$ induces a canonical homomorphism $K_\bullet(W) \to 
K_\bullet(W').$ In fact, any morphism $W \to  W'$ induces a 
commutative diagram 
\begin{equation}
\begin{split}
\xymatrix{
\Der(W) \ar[r] \ar[d] & \Der(W') \ar[d]\\
K_\bullet(W) \ar[r] & K_\bullet(W')
}
\end{split}
\end{equation}
Furthermore, for any exact functor $\Der (W) \to \Der (W'),$
the Chern character ring homomorphism $ch : K_\bullet(W) \to H^*(W,\CC)$
induces a canonical map in cohomology  
$H^*(W,\CC) \to H^*(W',\CC).$ In fact, Mukai \cite{Mukai1} proved 
that the set of Fourier--Mukai functors acting on $\Der (W)$ induced by
kernels in $\Der (W \times W)$ admits also a natural 
monoidal (with unity) multiplicative structure. 
In what follows, for a complex 
$\Ea^\bullet \in \Der(W),$ we will denote by $ch(\Ea^\bullet)$ 
the Chern character of the corresponding element in  $K_\bullet(W).$

For any smooth Calabi--Yau manifold $W$ consider 
the diagonal embedding 
$j:\Delta \to W \times W$. This embedding is the graph of the identity
map on $W$, so $j_* {\cal O}_{\Delta}$ acts as the
identity on $\Der (W)$. For any line bundle $\cal L$ on $W \cong \Delta$, 
the kernel 
\begin{equation}
p_2^*(\La) \otimes j_* {\cal O}_{\Delta} =
j_* (\La)
\label{eq:linebundle}
\end{equation}
determines an action on $H^{*}(W,\CC)$ given by 
\begin{equation}
\label{eq:line}
\gamma \mapsto ch(\cal L) \cdot \gamma.
\end{equation}

In string theoretic language, this transformation corresponds to an 
automorphism of the physical theory given by an integral shift 
in the B--field. Indeed, the effect of transformation 
(\ref{eq:line})
on $H^*(W,\CC)$ is mirrored by the monodromy action on $H^*(M,\CC)$
of a loop around a divisor passing through the corresponding 
maximal unipotent point in the moduli space of complex structures 
on $M.$ If the maximal unipotent point is given by $z_1=\ldots=z_{N-n}=0,$
then the action of the loop defined by $z_{j_0} \mapsto e^{2 \pi i t} z_{j_0}, 
0 \leq t \leq 1,$ and $z_j=$ constant, for $1 \leq j \leq N-n, j \not=j_0,$
on the series $\Phi_\lambda(z)$ (given by by (\ref{phi})) is expressed by  
\begin{equation}\label{eq:trran}
\Phi_\lambda(z) \mapsto e^{2\pi i \lambda_{j_0}} \Phi_\lambda(z).
\end{equation}
According to \cite{CDGP}, the flat ($\sigma$-model) coordinates 
on the moduli space of complexified K\"ahler forms on $W$
can be defined by 
\begin{equation}\label{def:bf}
B_j + i J_j := \frac{\Phi_j (z)}{\Phi_0 (z)}, \; 1 \leq j \leq N-n,
\end{equation}
where $\Phi_0 (z)= \Phi_\lambda (z)_{|\,\lambda=0}$ is the unique 
analytic period at $z_1=\ldots=z_{N-n}=0,$ 
and $\Phi_j (z)= \Big(\frac{\partial}{\partial \lambda_j} \Phi_\lambda (z)
\Big)_{|\,\lambda=0}$ are the  
periods with logarithmic monodromy at that point.
The fact that the transformation (\ref{eq:trran}) (the mirror of the 
transformation (\ref{eq:line}))
induces an integral shift 
in the B--field coordinates is then a direct consequence of 
the definition (\ref{def:bf}).
For a thorough discussion of this issue, see \cite{MTr}.

In what follows, when there is no possible confusion, we will denote the 
sheaf $j_* {\cal O}_{\Delta}$ by $\Oa_\Delta.$

A line bundle $\cal L$ on $W$ determines also the complex:
\begin{equation}
\label{eq:cplex}
\ldots \to 0 \to p_2^*({\cal L}) \otimes p_1^* ({\cal L^{\vee}})
\to {\cal O}_{\Delta} \to 0 \to \ldots.
\end{equation}
(with ${\cal O}_{\Delta}$ at the $0$-th place).

The Grothendieck--Riemann--Roch formula gives the 
action of this complex on $H^{*}(W,\CC)$
\begin{equation}
\label{eq:linecoh}
\gamma \mapsto \gamma - (\int_W \gamma \cdot ch({\cal L}^{\vee}) \cdot 
{\rm Todd}_W) \cdot ch({\cal L}),
\end{equation}
where ${\rm Todd}_W$ is the Todd class of $W$.

In particular, for the trivial line bundle, we obtain a complex 
representing the ideal sheaf of the diagonal ${\cal I}_{\Delta}:$
\begin{equation}\label{eq:toddcplex}
\ldots \to 0 \to  {\cal O}_{W \times W} \to {\cal O}_{\Delta} 
\to  0 \to \ldots.
\end{equation}
The corresponding action in cohomology is given by:
\begin{equation}
\label{eq:todd}
\gamma \mapsto \gamma - (\int_W \gamma \cdot {\rm Todd}_W) \cdot {\rm 1}_W.
\end{equation}
According to Kontsevich, this is the kernel that should correspond
in general to a loop based at the large complex structure limit point which
surrounds once the principal component of the discriminant locus. 

We now describe a general construction of kernels for the case
of Calabi--Yau complete intersections in toric varieties. Later,
we will discuss a conjectural construction of the corresponding 
loops in the moduli space of complex structures of the mirror 
Calabi--Yau manifold $M.$ This will provide a more specific
interpretation of Kontsevich's conjecture for the case
of complete intersections in toric varieties. 
As a particular case,
the identification of the kernel (\ref{eq:toddcplex}) 
with a loop around the principal component of the discriminant locus 
will be given a natural explanation. 

As described in sections \ref{chap:toricgen} and \ref{chap:triang}, the 
smooth Calabi--Yau toric complete intersection 
$W=E_1 \cap \ldots \cap E_k \subset \PP_{\Delta}$ undergoes phase 
transitions encoded in the structure of the corresponding secondary polytope
(fan). Recall that the poset of faces (cones) of the secondary polytope 
(fan) is isomorphic to the poset of
regular polyhedral subdivisions of the $\cA$ set generating the
Gorenstein cone $\sigma.$ The smooth phase is determined by
a regular triangulation $\Ta_0$ of the set $\cA$ defining a toric 
structure on the space of the bundle $\EE_{\Ta_0} \longrightarrow
X_{\Sigma_0}=\PP_{\Delta},$ where $\Sigma_0$ is the fan defining the 
complete toric variety $\PP_{\Delta}.$ 

Consider a face $F$ of the secondary polytope 
with one vertex corresponding to the regular triangulation $\Ta_0.$ Under 
the poset equivalence described above, there exists a polyhedral subdivision 
$\Sa(F)$ of the $\cA$--set that is refined by 
the triangulations corresponding
to the vertices of the face $F.$ The polyhedral subdivision
$\Sa(F)$ defines a fan structure supported on the Gorenstein cone
$\sigma,$ and we have a toric morphism $f : \EE_{\Ta_0} \to
 \EE_{\Sa(F)}$ that induces a corresponding toric map 
$f: X_{\Sigma_0} \to X_{\Sigma_F}$ as well as a regular map
$f: W \to W'.$ The  
fan $\Sigma_F$ defining the complete toric variety $X_{\Sigma_F}=
f(X_{\Sigma_0})$ can be described explicitly. Namely, $f(X_{\Sigma_0})$
is the toric orbit $V(\tau') \subset \EE_{\Sa(F)}$ defined 
by the minimal cone $\tau'$ in $\Sa(F)$ that contains the image 
of the cone $\tau$ in $\Ta_0,$ with
\begin{equation}
\ovl{\tau}:= \RR_{\geq 0} \ov_1 + \ldots + \RR_{\geq 0} \ov_k.
\end{equation}
Recall that $\ov_1, \ldots, \ov_k$ belong to any
simplex in the triangulation $\Ta_0,$ and consequently 
$V(\ovl{\tau})= X_{\Sigma_0}$. Note that, in general, 
$X_{\Sigma_F}$ and $W'$ could be singular and of lower dimension 
than $X_{\Sigma_0}$ and $W,$ respectively. 

As we will see in section
\ref{chap:generalcase}, in the case when $F$ is an edge in the
secondary polytope, the toric map $X_{\Sigma_0} \longrightarrow X_{\Sigma_F}$
is an elementary contraction in the sense of Mori theory. These
transformations have been investigated in the toric context by M. Reid 
\cite{reid}. The interpretation should come as no 
surprise, as the Mori cone of the toric variety $X_{\Sigma_0}$
(generated by homology classes of effective curves on $X_{\Sigma_0}$)
is dual to its K\"ahler cone. But the K\"ahler cone
of the toric variety $X_{\Sigma_0}$ is isomorphic to the big cone 
in the secondary fan corresponding to the the triangulation $\Ta_0,$ so 
the edges of the secondary polytope starting at the vertex corresponding to
$\Ta_0$ are in one-to-one correspondence with the extremal rays
in the Mori cone. Therefore, the Mori cone of $X_{\Sigma_0}$ 
can be naturally identified with the ``big'' cone in the 
secondary polytope at the vertex given by the triangulation 
$\Ta_0$ (see remark \ref{rem:mori} for a more detailed 
discussion). 

Another interesting situation is the case when the face $F$ of the 
secondary polytope is the polytope itself. Then $X_{\Sigma_F}=W'$ is a 
point (the so-called Landau--Ginzburg point) and this is exactly 
the face which induces the automorphism (\ref{eq:toddcplex})
in the derived category. As a prelude to the refinement
to Kontsevich's conjecture that
we are about to present, note here that, under mirror symmetry,
we expect that the automorphism given by the maximal face of the 
secondary polytope matches with the 
automorphism given by a loop based at $p_{\Ta_0}$ that goes
around the principal component of the 
discriminant locus. By Horn uniformization (see remark \ref{rem:discr} 
below), this component 
is birationally isomorphic to a projective space and it is included in the 
open toric stratum determined in the moduli space of complex
structures of the mirror $M$ {\it exactly} by the polytope itself
(or, equivalently, by the origin viewed as a subcone in the
secondary fan). 

Let $E_0 \subset X_{\Sigma_0}$ and 
$Z_0 \subset X_{\Sigma_F}$ be the loci where the toric map 
$X_{\Sigma_0} \longrightarrow X_{\Sigma_F}$ is not an isomorphism. 
As it will be explained in section \ref{chap:generalcase}, the fiber of
this elementary contraction is a weighted projective space. In fact
\begin{equation}\label{eq:complint}
E_0= \bigcap_{j \in I_-'} D_j,
\end{equation}
where $D_j$ are toric divisors in $X_{\Sigma_0}$ (the peculiar
choice of the notation $I_-'$ for the index set will fit in well with 
the context of section \ref{chap:generalcase}). The fibered product 
$E_0 \times_{Z_0} E_0$ can also be described torically.
If $E$ is the intersection 
\begin{equation}
E:=W\, \cap \, E_0,
\end{equation}
and $Z:=f(E),$ let $Y$ be the fibered product 
\begin{equation}
Y:= E \times_{Z} E.
\end{equation} 
In general, the exceptional locus of
the restriction of the morphism $f$ to $W$ is contained in $E,$
but may be smaller than $E.$

By naturality, the following diagram is commutative. 
\begin{equation}\label{diagr:big}
\begin{split}
\xymatrix{
  &   & W \times W \ar[lldd]_{p_2} \ar[rrdd]^{p_1}  & & \\
  &   & \overset{\;}Y \ar@{^{(}->}[u]^i \ar[dl]_{q_2} \ar[dr]^{q_1}
& &\\
W &\;E \ar@{_{(}->}[l]^{j_2} \ar[dr]^{r_2}& 
  & E\; \ar@{^{(}->}[r]_{j_1} \ar[dl]_{r_1} &W\\
  &                         & Z &                         &}
\end{split}
\end{equation}  

Denote by $\Delta' \cong E $ the diagonal in $Y= E \times_Z E.$ 
Consider the complex of sheaves $\Fa^{\prime\bullet}$ 
on $W \times W$ given by 
\begin{equation}\label{eq:trivialcplx}
\ldots \to 0 \to \Oa_Y \to 0 \to \ldots.
\end{equation}
The exact position of the sheaf $\Oa_Y$ will automatically be 
specified by our discussion. Define also
the complex of sheaves $\Fa^{\prime\prime\bullet}$ on $W\times W$ 
given by
\begin{equation}\label{eq:toddnewcplex}
\begin{split}
\ldots & \to 0 \to \bigoplus_{k \in I_-'} \Oa_{\Delta} 
(-\sum_{j \in I_-' \setminus \{k\} } D_j ) \to 
\ldots \to \bigoplus_{j' < j'' \in I_-'} \Oa_{\Delta} (-D_{j'} - 
D_{j''}) \to \\
& \to \bigoplus_{j \in  I_-'} \Oa_{\Delta} (-D_{j'}) \to \Oa_\Delta
\to \Oa_{\Delta'} \to  0 \to  \ldots
\end{split}
\end{equation}
(with the first non-zero sheaf at the $0$-th place, or, equivalently,
with the sheaf $\Oa_{\Delta'}$ at the $|I_-'|$-th place). 
This complex is a truncated Koszul resolution of the sheaf 
${\cal O}_{\Delta'},$ associated to the complete intersection
$\Delta' \cong E$ in $\Delta \cong W$ (see (\ref{eq:complint})).
The only non-zero cohomology sheaf is at
the $0$-th place and it is the sheaf of sections of the line bundle 
$\La^\vee,$ with 
\begin{equation}\label{def:lbun}
\La=\Oa(\sum_{j \in I_-'} D_j).
\end{equation}
If $I_-'=\emptyset$ (i.e. $E=W$), then $\La=\Oa_W$ and 
$\Fa^{\prime\prime\bullet}$ is the complex
\begin{equation}
\ldots \to 0 \to \Oa_\Delta\to 0 \to \ldots.
\end{equation}
 
We also note here that, since $W$ is Calabi--Yau and $E$ is a complete
intersection in $W,$ the adjunction formula gives that
\begin{equation}
\La_{|\, E}=K_E , 
\end{equation}
where $K_E$ is the canonical bundle of $E.$

Consider now the morphism in the derived category 
$\Der (W \times W)$ 
\begin{equation}
\Fa^{\prime\bullet} \to  \Fa^{\prime\prime\bullet},
\end{equation}
with the only non-zero map given by the restriction 
$\Oa_Y \to \Oa_{\Delta'}$ (note that this fact specifies 
the place of the sheaf $\Oa_Y$ in the complex $\Fa^{\prime\bullet}$). 

A fundamental property of the 
derived category is the fact that it is a triangulated category. 
In particular, this means that the above morphism can 
be completed to a distinguished triangle in $\Der(W \times W)$
\begin{equation}
\Fa^{\prime\bullet} \to  \Fa^{\prime\prime\bullet} \to
\Fa^{\bullet} \to \Fa^{\prime\bullet}[1].
\end{equation}
In our concrete example, the complex $\Fa^{\bullet}$ 
is obtained by the mapping cone construction and has the 
form
\begin{equation}\label{eq:fbullet}
\begin{split}
\ldots & \to 0 \to \bigoplus_{k \in I_-'} \Oa_{\Delta} 
(-\sum_{j \in I_-' \setminus \{k\} } D_j ) \to 
\ldots \to \bigoplus_{j' < j'' \in I_-'} \Oa_{\Delta} (-D_{j'} - 
D_{j''}) \to \\
& \to \bigoplus_{j \in  I_-'} \Oa_{\Delta} (-D_{j'}) \to 
\Oa_Y \oplus \Oa_\Delta \overset{\tau} \to \Oa_{\Delta'} \to  0 \to  \ldots,
\end{split}
\end{equation}
with the sheaf $\Oa_{\Delta'}$ placed at the $|I_-'|$-th place.
For a local section $(f,h)$ of the sheaf $\Oa_Y \oplus \Oa_\Delta,$ the 
map $\tau$ is given by 
\begin{equation}
\tau(f,h) = f_{|\, \Delta'} + h_{|\, \Delta'}.
\end{equation}

The refinement to Kontsevich's conjecture that we are about  
to describe will use the action on $H^*(W,\CC)$
induced by a twisted version of the complex $\Fa^\bullet.$ It will
be denoted by $\Ea^\bullet (F)$ to indicate clearly the fact 
that this complex is determined by the face $F$ of the secondary 
polytope associated to the toric context.

\begin{definition}\label{def:bigdef}
The complex $\Ea^\bullet (F) $ in $\Der(W \times W)$ 
is defined by 
\begin{equation}\label{eqn:facecomplex}
\Ea^\bullet (F):= p_2^* (\La) \otimes \Fa^\bullet,
\end{equation}
where $\Fa^\bullet \in \Der(W \times W)$ is the complex
(\ref{eq:fbullet}).
\end{definition}

\begin{remark}\footnote{I want to thank R.~P. Thomas for correcting an 
earlier version of this remark.}
In the case $E=W,$ the line bundle $\La(F)$ is trivial and $Y$
is a subvariety of $ W\times W,$ so the complex 
$\Ea^\bullet (F)$ reduces to the restriction map 
\begin{equation}
\ldots \to 0 \to \Oa_Y \to \Oa_\Delta \to 0 \to \ldots.
\end{equation}
Note that, as expected, when $Z$ is a point (i.e. when $Y= W \times W$), 
$\Ea^\bullet (F)$ is the complex (\ref{eq:toddcplex}). 
In the case $|I_-'|=1,$
the complex $\Ea^\bullet(F)$ is written as  
\begin{equation}
\ldots \to 0 \to p_2^* (\La) \otimes (\Oa_Y \oplus \Oa_\Delta) 
\to p_2^* (\La) \otimes \Oa_{\Delta'} \to  0 \to  \ldots,
\end{equation}
with the first non-zero sheaf at the $0$-th place.
\end{remark}

By Grothendieck--Riemann--Roch, the action induced on $H^*(W,\CC)$
by the complex $\Ea^\bullet (F)$ is given by
\begin{equation}\label{eq:grr}
\gamma \mapsto p_{2_*} \big( ch (\Ea^\bullet(F)) \cdot \Td_W \cdot 
p_1^* (\gamma)\big).
\end{equation}

\begin{lemma} \label{lemma:actionn}
Given a complex $\Ia^\bullet$ on 
$Y= E \times_Z E,$ the induced actions on $H^*(W,\CC)$
of the Fourier--Mukai functors $\Phi_\Ia$ and $\Phi_{i_*{\Ia}},$
with $i : Y \hookrightarrow W,$ are equal.
\end{lemma}

Note that $\Phi_{\Ia^\bullet}$ is defined as a 
Fourier--Mukai functor induced by the 
kernel $\Ia^\bullet$ on $Y$ by using the projections $j_1 \circ q_1$ and
$j_2 \circ q_2$ (in the notation of the diagram (\ref{diagr:big})).
The lemma states that the action in cohomology induced by the kernel 
$i_*{\Ia^\bullet} \in \Der(W \times W)$ ($i_*$ denotes here the 
corresponding derived functor) can also be computed as 
the action in cohomology
induced of the kernel $\Ia^\bullet \in \Der (Y).$

\begin{proof}
We make repeated use of the Grothendieck--Riemann--Roch
formula and of the properties of the Gysin map in cohomology. For 
$\Ha \in \Der (W),$ we have that $ch(\Phi_{\Ia^\bullet} (\Ha)) \cdot \Td_W$ 
is given by 
\begin{equation}\label{eq:GRR}
\begin{split}
& (j_2 \circ q_2)_* \Big( (j_1 \circ q_1)^* (ch(\Ha)) \cdot ch(\Ia^\bullet) 
\cdot \Td_Y \Big) =\\
& (p_2 \circ i)_* \Big( (p_1 \circ i)^* (ch(\Ha)) \cdot ch(\Ia^\bullet) \cdot
\Td_Y \Big)= \\
& p_{2_*} \circ i_* \Big(i^*(p_1^* (ch(\Ha))) \cdot ch(\Ia^\bullet) \cdot 
\Td_Y \Big)= \\
& p_{2_*} \Big( p_1^* (ch(\Ha)) \cdot ch(i_*(\Ia^\bullet)) 
\cdot \Td_{W\times W} \Big).
\end{split}
\end{equation} 
The last expression is exactly $ch(\Phi_{i_*{\Ia^\bullet}}(\Ha)) \cdot \Td_W,$ 
which proves the lemma.
\end{proof}

\begin{proposition} \label{prop:action}
The induced action on $H^*(W, \CC)$
of the automorphism $\Phi_{\Ea^\bullet(F)}$ is given by 
\begin{equation}
\gamma \mapsto \gamma - \prod_{j \in I_-'} (1-e^{[D_j]})
q_{2_*} \Bigl( (j_1 \circ q_1)^* (\gamma) \cdot  
\Td_Y \cdot q_2^* \bigl( (\Td_E)^{-1} \bigr)
\end{equation}
\end{proposition}


\begin{proof} According to the previous lemma, the induced action on 
$H^*(W,\CC)$ of the Fourier--Mukai functor determined by the 
complex $\Fa^{\prime\bullet} \in \Der(W\times W)$ is identical to the action 
Fourier--Mukai functor determined the same complex viewed 
as a complex of sheaves on $Y.$ The lemma shows then that, for $\gamma
\in  H^*(W,\CC),$ this action is given by 
\begin{equation}\label{eq:actt}
\gamma \mapsto \pm (j_{2_*} \circ q_{2_*}) \bigl( (j_1 \circ q_1)^* (\gamma)
\cdot \Td_Y \cdot q_2^* \bigl( (\Td_E)^{-1} \bigr) \bigr),
\end{equation}
where the sign is determined by the parity of $|I_-'|.$

But the map $j_2$ embeds $E$ in $W$ as a complete 
intersection of the divisors $D_j, j \in I_-',$ so again, by 
Grothendieck--Riemann--Roch, it follows that the action of the complex 
$\Fa^\bullet$ on $H^*(W,\CC)$ can be written as 
\begin{equation} 
\gamma \mapsto e^{-[\La]} \cdot \gamma -
\prod_{j \in I_-'} (e^{-[D_j]}-1) 
\Bigl(q_{2_*} \bigl( (j_1 \circ q_1)^* (\gamma) \cdot
\Td_Y \cdot q_2^* \bigl( (\Td_E)^{-1} \bigr) \bigr) \Bigr).
\end{equation}
This formula is implied by (\ref{eq:actt}) and by  
the fact that the complex $\Fa^{\prime\prime\bullet}$ is equivalent
to the sheaf $\Oa_\Delta(\La^\vee)$ placed at the $0$-th position.
The proposition follows as a direct consequence of definition \ref{def:bigdef}.

\end{proof}

\begin{remark} \label{rem:discr}
To make a precise statement regarding the more general automorphisms 
on the B--side, we need to review a few results about discriminants
and their relationship to $\cA$--sets following Gelfand, Kapranov and
Zelevinsky \cite{GKZ2}. 

In what follows, we use the notation introduced in 
section \ref{chap:triang}. Given an $\cA$--set in $\ZZ^n$  with
$\cA= \{ \ov_1, \ldots, \ov_N \},$ one 
can consider the set of polynomials $p$ of the form 
$p=\sum_{1 \leq j \leq N} 
c_j x^{\ov_j}$ $\in \CC[x_1, \ldots, x_n].$ Define 
$\nabla_{\cA}$ to be the Zariski closure of
the polynomials $p$ as above such that there exists some
$y \in (\CC^*)^n$ with the property that $f$ is singular at $y.$
By definition, the discriminant $\Delta_{\cA} \in \ZZ[c_1, \ldots, 
c_N]$ is the irreducible polynomial (defined up to a sign)
whose zero set is given by the union of the irreducible
components of codimension $1$ of $\nabla_{\cA}.$ In fact, we shall
be interested only in the case when $\nabla_{\cA}$ is a
hypersurface, and in that case $\nabla_{\cA}$ is irreducible
(for the case codim$\nabla_{\cA} >1,$ one sets  $\nabla_{\cA}=1$).
Then $\nabla_{\cA}$ is called the principal component of the
discriminant locus. The Horn uniformization theorem (proved
in this case by Kapranov, see theorem 3.1 in
\cite{GKZ2}) states that the discriminantal hypersurface
$\nabla_{\cA} \subset (\CC^*)^n  $ is birationally equivalent
to the projective sphere $\PP^{n-1}.$ Using the combinatorial
structure of $\cA,$ it is possible to write down a rational
parameterization for $\nabla_{\cA}$ (we will do this explicitly 
in section \ref{chap:2par}). 

The full discriminant $\Delta$ is obtained by taking the 
product
\begin{equation}
\Delta:= \prod_{\Gamma \subset Q} \Delta_{\cA \cap \Gamma},
\end{equation}
where $\Gamma$ are non-empty faces of $Q=Conv (\cA).$ It turns
out \cite{GKZ3} that the full discriminant $\Delta$ 
defines the singular locus of the GKZ system of differential
equations, so in particular it defines the locus in the moduli space 
of complex structures of singular mirror Calabi--Yau manifolds $M.$ 

Any face $F$ of the secondary polytope determines a polyhedral 
subdivision $\Sa(F)=(Q_j, \cA_j)$ of $(Q,\cA)$ (see section
\ref{chap:triang}). If $\LL$ is the lattice of relations 
among the elements of $\cA,$ for any $j,$ define
\begin{equation}\label{eq:newrel}
\begin{split}
\LL_{\cA_j}&:=\{ l=(l_i) \in \LL\; : \; l_i=0,\; \text{for} \; 
\ov_i \notin \cA_j \},\\
\LL_F& := \oplus_j \LL_{\cA_j}.
\end{split}
\end{equation}
Note that the both the face $F$ of the secondary polytope
as well as the sublattice $\LL_F$ are naturally embedded 
in $\LL.$ In fact, according to corollary 2.7 on page 230 in \cite{GKZ2}, 
the affine span of the face $F$ is a parallel translation of 
$\LL_F,$ so $\dim(F)=\dim(\LL_F).$ Since the normal fan of
secondary polytope is the secondary fan, the subfan of the secondary fan
that is dual to the face $F$ determines a $\dim(F)$--dimensional 
toric orbit in the moduli space of complex structures corresponding
to some subset of vectors $\cA' \subset \cA.$ But the principal 
component $\nabla_{\cA'}$ of the discriminant locus corresponding
to $\cA'$ is included in this open toric orbit and 
(again by Horn uniformization) is 
birationally isomorphic to a projective space $\PP^{\dim(F)-1}.$
In general, this projective space will be part of the compactification
of the full discriminant locus corresponding to $\cA.$ 
Consider a loop based near the large complex structure point 
corresponding to the triangulation $\Ta_0$ that is included in the 
toric stratum of dimension $\dim(F)$
determined by $F.$ Since, in general,
this toric stratum is part of the full discriminant locus,
we may have to ``separate'' the loop from the toric stratum, so let
$\gamma_F$ be a loop which is homotopic to the previous one,
but does not intersect the discriminant locus.
\end{remark}

The previous discussion applies now verbatim to the toric map 
$X_{\Sigma_0} \to X_{\Sigma_F},$ whose existence is a consequence
of the fact that the triangulation $\Ta_0$ (corresponding to the 
smooth Calabi--Yau phase $W$) is a refinement of the polyhedral 
subdivision $\Sa(F).$ The definition of the kernel $\Ea^\bullet(F)$
is completely analogous. We can now state a refinement of Kontsevich's 
conjecture.

\begin{newconj} \label{newconj} For any face $F$ of the secondary 
polytope with one vertex corresponding to the smooth phase $W,$
the action on $H^*(W,\CC)$ induced by the kernel $\Ea^\bullet(F)$ (given 
by definition \ref{def:bigdef}) is identified under mirror 
symmetry with the monodromy action on $H^*(M,\CC)$ determined
by the loop $\gamma_F.$
\end{newconj}

This conjecture is consistent with all the constructions 
described here. The next sections will be devoted to a detailed study 
of this identification. Besides the one and two parameter cases that 
will be studied in detail, by using the monodromy calculations of 
section \ref{chap:mon}, we will be able to prove the conjecture for 
the case when $F$ is an edge in the secondary polytope. 

\subsection{The One--Parameter Case}\label{chap:1par}

Let $W=E_1 \; \cap \ldots \; \cap E_k \; $ be a generic {\it smooth} 
Calabi--Yau complete intersection in some weighted projective space 
$X=\PP(q_1,\ldots,q_{n+1}),$ with $E_j$ hypersurfaces of homogeneous
degree $d_j$ for 
$j,$ $1 \leq j \leq k.$ We can always assume that the g.c.d. of 
$q_1, \ldots, q_{n+1}$
is $1,$ and define $q:= q_1 +\ldots +q_{n+1}.$ According
to lemma 3.5.6. in \cite{CK}, $X$ is Fano if and only
if $q_j \vert q$ for all $j.$ The Calabi--Yau condition implies that
\begin{equation}
d_1 + \ldots +d_k= q.
\end{equation}
In general, a weighted projective space is singular, but, in our
case, $W$ does not intersect the singular locus of $X.$ Possibly after
a change of coordinates (which may produce a branched cover of the
original weighted projective space), we can assume that $q_{n+1}=1.$


We place our discussion in the context of sections \ref{chap:cicy} 
and \ref{chap:gkzcicy}, so we consider the  
set $\cA= \{ \ov_1, \ldots , \ov_{k+n+1} \} \subset \ZZ^{k+n},$ 
describing torically this situation. As before, the vectors 
$\ov_1, \ldots, \ov_k$ have their first $k$ components equal to be the 
standard basis vectors in $\ZZ^k$ and the rest of 
$n$  components equal to zero. The vectors 
$\ov_{k+1}, \ldots,  \ov_{k+n+1}$ are obtained from a system of vectors 
$v_{k+1}, \ldots ,$ $v_{k+n+1}$ that generate a fan in $\ZZ^n$ defining 
torically the weighted projective space $\PP(q_1,\ldots,q_n,1).$
This fan is the interior point fan of the reflexive polyhedron determined 
by the vectors $v_j.$ The hypersurfaces $E_1, \ldots, E_k$ correspond to a
partition of the $n+1$ vectors generating the fan of 
$\PP(q_1,\ldots,q_n,1)$ in subsets 
$\{v_{k+1}, \ldots, v_{k+d'_1} \},\ldots,$ $\{ v_{k+d'_{k-1}+1}, \ldots, 
v_{k+n+1}\},$ $d'_1 + \ldots + d'_k= n+1.$
One possible choice for the vectors $v_j$ is 
\begin{equation}\label{def:wpr}
\begin{split}
v_j &= (0, \ldots, \overset{j-k}1, \ldots , 0), \ 
k+1 \leq j \leq k+n, \\
v_{k+n+1} &= (-q_1, \ldots, -q_n).
\end{split}
\end{equation}
The standard basis vector $e_i$ in $\ZZ^k$ is adjoined to a vector $v_j$ 
if the divisor in $\PP(q_1,\ldots,q_n,1)$ 
corresponding to $v_j$ contributes to the hypersurface $E_i,$
$1\leq i \leq k.$ The lattice of relations $\LL$ among the elements
of $\cA$ is generated by the relation
\begin{equation}\label{onlyrel}
q_1 \ov_{k+1}+ \ldots + q_n \ov_{k+n} + \ov_{k+n+1}
= d_1 \ov_1+ \ldots + d_k \ov_k .
\end{equation}
In this case $\cA$ admits two possible regular triangulations. 
The triangulation describing the weighted projective space itself
has maximal simplices of
the form $\{ 1, \ldots, k+n+1 \} \setminus \{j \},$ for all $j,$
$k+1 \leq j \leq k+n+1,$ while the triangulation corresponding
to the Landau--Ginzburg point has maximal simplices 
of the form $ \{ 1, \ldots, k+n+1 \} \setminus \{j \},$ for all $j,$
$1 \leq j \leq k.$

Using the description (\ref{batmir}) of Batyrev's mirror construction
in section \ref{chap:toricgen}, the Calabi--Yau varieties $M_x$
in the mirror family can be represented as compactifications of 
complete intersections of the 
affine hypersurfaces in $(\CC^*)^{n+1},$ given by 
\begin{equation}\label{eq:affine}
\begin{split}
t_1 + \ldots + t_{d'_1} &= 1 \\
\ldots & \\
t_{d'_1+\ldots d'_{k-1} +1} + \ldots + t_{n+1} &= 1\\
t_1^{q_1} \cdot \ldots \cdot t_n^{q_n} \cdot t_{n+1} &= x.
\end{split}
\end{equation}
Hence, we obtain a one--parameter family of mirror Calabi--Yau
manifolds $M_x.$ The value $x=0$ is the large
complex structure limit point, and, as first remarked by 
Batyrev and Van Straten \cite{BatVS}, the unique analytic
solution around $x=0$ to the Picard--Fuchs equations is given by 
\begin{equation}\label{eq:series}
\sum_{m \geq 0} \frac{\prod_{j=1}^k (d_j m)!}
{\prod_{j=1}^{n+1}(q_j m)!} \;  x^m.
\end{equation}
To make things more clear, we describe a  concrete
example. 

\begin{example} The quintic family in $\PP^4$ is defined by the equation
\begin{equation}
y_1^5+ \ldots + y_5^5 - x^{1/5} y_1 \cdot \ldots \cdot y_5 =0.
\end{equation}
If we change the $y$--coordinates to $\ZZ_5^4$--invariant $t$--coordinates 
(page 479 in \cite{Giv2}) using the change of variable
\begin{equation}
t_j:= x^{1/5} \frac{y_j^5}{y_1 \cdot \ldots \cdot y_5}, \ 1 \leq i \leq 5,
\end{equation}
then, according to the Greene--Plesser orbifolding construction,
an affine part of the mirror quintic $M_x$ is given by the equations
\begin{equation}
\begin{split}
t_1+\ldots+ t_5 & =1 \\
t_1\cdot \ldots \cdot t_5 &= x.
\end{split}
\end{equation}
The Calabi--Yau manifold $M_x$ is smooth for all 
$x \not= 0,\infty, 5^{-5}.$ The unique analytic solution at $x=0$ to the 
Picard--Fuchs equations corresponding to this family of 
Calabi--Yau manifolds is given by 
\begin{equation}
\sum_{m \geq 0} \frac{(5m)!}{(m!)^5} \; x^m.
\end{equation}

\end{example}

In general, the Calabi--Yau spaces $M_x$ defined by (\ref{eq:affine})
are smooth for all values of the 
complex parameter $x$ different from $0, \infty,$ and
$a,$ with 
\begin{equation}
a=\prod_{j=1}^k d_j^{\;-d_j} \prod_{j=1}^{n} q_j^{q_j}.
\end{equation}
The convergence radius of the series 
(\ref{eq:series}) is  $\prod_{j=1}^k d_j^{\;-d_j}\prod_{j=1}^{n} 
q_j^{q_j}.$

The discussion in the previous sections shows that we can write
a series that incorporates all the solutions to the 
Picard--Fuchs equations describing the variations of complex structure 
of $M_x.$ Let $\Psi_\lambda (x)$ be the series defined by
\begin{equation}
\Psi_\lambda (x):= \sum_{m \geq 0} 
\frac{\prod_{j=1}^k \Gamma(d_j(m+\lambda)+1)}
{\prod_{j=1}^{n+1} \Gamma(q_j(m+\lambda)+1)} \; x^{m+\lambda}.
\end{equation}
It is not hard to see that the series $\Psi_\lambda$ is just
an adaptation of the series $\Phi_\lambda$ considered before-- if we
introduce the integer $\epsilon$ by (compare to (\ref{def:epsilon}))
\begin{equation}\epsilon:=
\begin{cases}
1& \hbox{if $q$ is odd},\\
2& \hbox{if $q$ is even},
\end{cases}\label{def:epsilon1}
\end{equation}
we have that
\begin{equation}\label{eq:phipsi}
\begin{split}
\Phi_\lambda(e^{i \pi \epsilon} x)&= \sum_{m \geq 0} 
\frac{e^{i \pi \epsilon (m+\lambda)}}
{\prod_{j=1}^k \Gamma(-d_j(m+\lambda))} \; 
\frac{1}{\prod_{j=1}^{n+1} \Gamma(q_j(m+\lambda)+1)} \; x^{m+\lambda}\\
&=e^{i\pi\epsilon \lambda} \prod_{j=1}^k \frac{\sin(-\pi d_j\lambda)}
{\pi} \; \Psi_\lambda (x)\\
&= e^{i \pi (\epsilon + q) \lambda}
\prod_{j=1}^k \frac{e^{-2\pi i d_j \lambda}-1}
{2 \pi i } \; \Psi_\lambda (x).
\end{split}
\end{equation}
The cohomology of the smooth $(n-k)$-dimensional complete intersection
$W$ is $H^*(W,\CC)= \CC[\lambda] / (\lambda^{n+1-k}),$ 
so we obtain that the solutions around $x=0$ to 
the Picard--Fuchs equations are $\Psi_0(x), \ldots, \Psi_{n-k}(x)$ given by
\begin{equation}\label{eqcigamma}
\Psi_\lambda(x)= \sum_{j=0}^{n-k} \Psi_j(x) \lambda^j,
\end{equation}
where this formula is viewed as an equality that holds in the ring
$H^*(W,\CC).$ As we will see shortly, the 
parameter $2 \pi i \lambda$ should be interpreted as the class of the
restriction to $W$ of a degree one divisor in the weighted projective 
space (the divisor corresponding to the vector $v_{k+n+1}$).
The Stanley--Reisner ring of triangulation corresponding to 
$X=\PP(q_1,\ldots, q_n,1)$ is $\ZZ[\lambda] / (\lambda^{n+1})$ and 
we have that $H^*(X,\QQ) \cong \QQ[\lambda] / (\lambda^{n+1}),$ and, 
(using the previous notation) the ring $\Ha$ is equal in this case to
$H^*(W,\ZZ)$ $\cong \ZZ[\lambda] / (\lambda^{n+1-k}).$

We can now compute the monodromy of the $H^*(W,\CC)$--valued
series $\Psi_{\lambda}(x)$ along a counterclockwise loop based at $x=b$
near the origin $x=0$ that goes around the point 
$x=a$ (see again figure \ref{fig:cont2}, or
\ref{fig:cont3}). Before we apply the monodromy formula of 
theorem \ref{thm:mon}, we make the remark that the
only relation in $\LL$ is (\ref{onlyrel}), and, as noted before,
there are only two possible regular triangulations of the set $\cA.$ 
The transition between the triangulation
corresponding to the weighted projective space and the one
corresponding to the  Landau--Ginzburg point is defined by 
\begin{equation}
\begin{split}
I_- & = \{ 1, \ldots, k \}, \ h_j=-d_j, \ 1 \leq j \leq k,\\
I_+ & = \{ k+1, \ldots, k+n+1 \}, h_j=q_{j-k}, \ k+1 \leq j \leq k+n+1.
\end{split}
\end{equation}
We apply theorem \ref{thm:mon} to express the monodromy at the
point $x=a$ along a counterclockwise
loop based at a point near the maximal unipotent point $x=0.$
By analytic continuation along this loop, the $H^*(W,\CC)$--valued series 
$\Phi_\lambda (e^{i \pi \epsilon} x)$ becomes 
\begin{equation}\label{eq:mon11}
\Phi_\lambda(e^{i \pi \epsilon} x) -
\int_{C_\xi} e^{2 \pi i \epsilon' (\lambda- \xi)}
\frac{\prod_{j=1}^k (1- e^{-2 \pi i d_j \lambda})} 
{\prod_{j=1}^{n+1} (1-e^{-2 \pi i q_j\xi})} \;
\Phi_\xi (e^{i \pi \epsilon} x) \; d\xi,
\end{equation}
where the integer $\epsilon'$ depends on $q$ and it is given by 
$\epsilon'= [(q+2)/2]$ (see (\ref{def:eps})), which means that
$2\epsilon'= \epsilon + q.$ The contour $C_\xi$ encloses
the points $\xi \in  \bigcup_{j=1}^{n+1} A_{q_j},$
with $A_{q_j}$ defined by (\ref{def:A}) (compare also to (\ref{def:A+})).

In fact, what we really need is the corresponding 
monodromy formula for $\Psi_\lambda(x).$ Combining (\ref{eq:phipsi})
and (\ref{eq:mon11}) it can be seen that the analytic continuation
of $\Psi_\lambda (x)$ around the prescribed loop is given by 
\begin{equation}\label{eq:mon12}
\Psi_\lambda (x)- \int_{C_\xi}
\frac{\prod_{j=1}^k (1- e^{-2 \pi i d_j \xi})} 
{\prod_{j=1}^{n+1} (1-e^{-2 \pi i q_j\xi})} \;
 \Psi_\xi (x) \; d\xi.
\end{equation}
This is the point where we employ condition (\ref{cond2}) to
simplify he contour $C_\xi.$ Indeed, condition
(\ref{cond2}) implies that all the poles of the integrand,
except $\xi=0,$ are just apparent singularities. From now 
on we assume that the contour $C_\xi$ encloses {\em only}
the origin $\xi=0.$

As promised, set $\alpha= 2 \pi i \lambda$ (and $\beta= 2 \pi i \xi$),
and define $I_\alpha(x):=\Psi_\lambda (x).$ The monodromy formula
(\ref{eq:mon12}) becomes 
\begin{equation}\label{eq:mon13}
I_\alpha(x) - \frac{1}{2\pi i} 
\int_{C_\beta} \frac{\prod_{j=1}^k (1- e^{-d_j \beta})} 
{\prod_{j=1}^{n+1} (1-e^{-q_j \beta})}
I_\beta (x) \; d\beta,
\end{equation}
with the contour $C_\beta$ chosen accordingly to enclose the origin 
$\beta=0.$ But the contour integral with respect to $\beta$ is
just the coefficient of $\beta^{n-k}$ in the expression
\begin{equation}\label{eq:mon14}
\beta^{n+1-k} \; \frac{\prod_{j=1}^k (1- e^{-d_j \beta})} 
{\prod_{j=1}^{n+1} (1-e^{-q_j \beta})}
I_\beta (x)= 
\big( \prod_{j=1}^k \frac{1- e^{-d_j \beta}}{\beta}\big)
\; \big( \prod_{j=1}^{n+1} \frac{\beta}{1-e^{-q_j \beta}} \big)
I_\beta (x).
\end{equation}
Alternatively, we can view the parameter $\beta$ (or $\alpha$), as the 
class of the divisor corresponding to $v_{k+n+1}$ in 
$X=\PP(q_1,\ldots,q_n,1).$ Recall that $X$ is defined by the 
fan in $\ZZ^n$ with the generators $v_j$ given by (\ref{def:wpr}).
The vectors $v_j,$ $k+1 \leq j \leq k+n,$ form the standard
basis in $\ZZ^n,$ so the cone $\sigma$ generated by them determines a smooth 
point in X. If $D_j$ is the divisor in $X$ corresponding to $v_{k+j},$ 
it is easy to see that $[D_j]= q_j \beta.$ Standard intersection theory 
in toric varieties (see, for example, remark 10.9 in \cite{Dani}) shows 
that in $A^*(X)_\QQ$ we have 
\begin{equation}\label{eq:beta}
\big(\prod_{j=1}^n q_j \big) \; \beta^n= [D_{1}] 
\cdot \ldots \cdot [D_{n}]= [\sigma],
\end{equation}
where $[\sigma]$ is the class in $A^*(X)_\QQ$ given by the maximal
dimensional cone $\sigma,$ i.e. the class of the corresponding smooth 
point. The class in $A^*(X)_\QQ$ of the smooth complete intersection 
$W$ is $\prod_{j=1}^k d_j \beta,$ and this implies that 
\begin{equation}\label{eq:mon15}
\int_W \beta^{n-k} = \frac{\prod_{j=1}^k d_j}{\prod_{j=1}^n q_j}.
\end{equation}
This shows that the smoothness of $W$
implies that ${\prod_{j=1}^n q_j} \; \big\vert \; {\prod_{j=1}^k d_j},$
which is also a necessary (but, unfortunately, not sufficient) 
condition for our assumption (\ref{cond2}) to hold.

By combining (\ref{eq:mon14}) and (\ref{eq:mon15}), we conclude that
\begin{equation}
\begin{split}
& I_\alpha(x) - \frac{1}{2\pi i} 
\int_{C_\beta} \frac{\prod_{j=1}^k (1- e^{-d_j \beta})}
{\prod_{j=1}^{n+1} (1-e^{-q_j \beta})}\; 
I_\beta (x) \; d\beta =\\
&=I_\alpha(x) - \int_W 
\big( \prod_{j=1}^k \frac{1- e^{-d_j \beta}}{d_j \beta}\big)
\; \big(\prod_{j=1}^{n+1} \frac{q_j \beta}{1-e^{-q_j \beta}} \big)\;
I_\beta (x).  
\end{split}
\end{equation}
According to the last part of proposition \ref{toricoh}, the Todd class
of the smooth complete intersection $W$ is given by 
\begin{equation}
\Td_W= \big( \prod_{j=1}^k \frac{1- e^{-d_j \beta}}{d_j \beta}\big)
\; \big(\prod_{j=1}^{n+1} \frac{q_j \beta}{1-e^{-q_j \beta}} \big).
\end{equation}
In other words, the action corresponding to the
loop under discussion is
\begin{equation}
I_\alpha (x) \mapsto I_\alpha(x) - \big(\int_W \Td_W \cdot I_\alpha(x)\big) 
\cdot 1_W.
\end{equation}
This means that we can identify indeed the monodromy action
on the cohomology of the mirror Calabi--Yau manifold $M$ 
corresponding to a counterclockwise loop around the point
$x=a=\prod_{j=1}^k d_j^{\;-d_j}\prod_{j=1}^{n} q_j^{q_j} $ 
with the action on the cohomology of the complete intersection
$W$ induced by the ideal sheaf of the diagonal $\Ia_{\Delta}$ in
$W \times W$ (formulae (\ref{eq:toddcplex}) and (\ref{eq:todd})). 

As explicitly shown by Kontsevich in his Rutgers lecture \cite{Kont2},
in this case, we have a complete dictionary between loops and kernels, 
as follows:

\iffigs
\begin{figure}
  \centerline{\epsfxsize=7cm\epsfbox{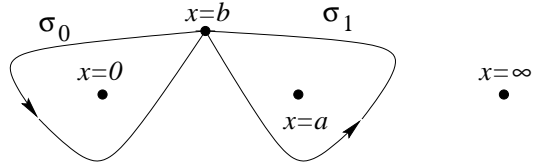}}
  \caption{The loops $\sigma_0$ and $\sigma_1$ (based at $x=b$)
generate the fundamental group of 
$\PP^1 \setminus \{0,a, \infty \},$ with $a=\prod_{j=1}^k d_j^{\;-d_j}
\prod_{j=1}^{n} q_j^{q_j}.$}
  \label{fig:cont3}
\end{figure}
\fi

\begin{thm}\label{thm:1par}
\
\begin{itemize}
\item The kernel $p_2^*(\La) \otimes j_* {\cal O}_{\Delta},$ 
with the line bundle $\La$ determined by the restriction to $W$
of a degree one divisor in $\PP(q_0, \ldots, q_n,1),$ 
corresponds to the loop $\sigma_0$ in figure \ref{fig:cont3}. 
\item The kernel $\Ia_\Delta$ given by the ideal sheaf of the 
diagonal in $W \times W$ (or, equivalently, defined by the complex
(\ref{eq:toddcplex})) corresponds to the loop $\sigma_1$ in figure
\ref{fig:cont3}.
\item The kernel (\ref{eq:cplex}) defined by the complex 
\begin{equation}
\ldots \to 0 \to p_2^*({\cal L}) \otimes p_1^* ({\cal L^{\vee}})
\to j_* {\cal O}_{\Delta} \to 0 \to \ldots ,
\end{equation}
with $j_* {\cal O}_{\Delta}$ at the $0$-th place and
$\La$ as above,
corresponds to the loop $\sigma_0 \cdot \sigma_1 \cdot \sigma_0^{-1}$
(this is a direct consequence of the previous two facts).
\end{itemize}
\end{thm}

\subsection{The Two--Parameter Case}\label{chap:2par}

We discuss now the mirror correspondence between a class of 
Calabi--Yau smooth complete intersections $W$
with two dimensional K\"ahler moduli space and the 
mirror Calabi--Yau manifolds $M$ with two dimensional complex moduli space.
The class generalizes the K3 fibrations Calabi--Yau
hypersurfaces in $\PP(2,2,2,1,1)$ and $\PP(6,2,2,1,1)$
considered, for example, in \cite{2par}, \cite{DavePles}, \cite{KachruVafa}. 
The behavior of the periods given by the GKZ system in two--parameter 
examples was also analyzed in \cite{GMV}, \cite{GY}. 
The Calabi--Yau complete intersections that we study are also
fibrations over $\PP^1$ with Calabi--Yau fibers. In some sense,
this is the simplest class of two--parameter examples that can be
constructed using weighted projective spaces, because our 
assumptions imply that we only need to resolve an $A_1$
type singularity.

Let $\widetilde{W}=\widetilde{E}_1 \; \cap \ldots \; \cap \widetilde{E}_k \; $ 
be a generic $(n+1-k)$-dimensional Calabi--Yau complete intersection in 
the weighted projective space $\PP(2q_1, \ldots,2q_n,1,1)$
with $\widetilde{E}_j$ hypersurfaces of homogeneous degree $2d_j,$
$q:=q_1+ \ldots +q_n +1=d_1+ \ldots+ d_k.$
As before, lemma 3.5.6. in \cite{CK} implies that this weighted
projective space is Fano if and only if $q_j \vert q$ for all $j.$
In particular, this implies that the g.c.d. of $q_1, \ldots, q_{n}$ is 1.

The singular locus of the weighted projective space 
$\PP(2q_1, \ldots,2q_n,1,1)$ may have various strata corresponding
to the possible non--trivial common factors of the integers $q_j.$ The
assumption we make is that $\widetilde{W}$ intersects {\it only} 
the largest singular stratum (singularities of type $A_1$) and 
{\it none} of the lower dimensional ones,
\begin{equation} 
\hbox{Sing}(\widetilde{W})= \widetilde{W} \; \cap \; \{(x_1, \ldots, x_{n+2}):
x_{n+1}= x_{n+2}= 0 \}
\end{equation}
and the intersection is transverse. 


The variety $\PP(2q_1, \ldots,2q_n,1,1)$ is described torically
by the fan in $\ZZ^{n+1}$ which is the interior point fan of a reflexive 
polytope and its large cones are generated by any distinct $n+1$ vectors 
from the set of vectors 
\begin{equation}
\begin{split}
v_i &= (0, \ldots, \overset{i-k}1, \ldots , 0), \ 
k+1 \leq i \leq k+n+1, \\
v_{k+n+2} &= (-2q_1, \ldots, -2q_n,-1).
\end{split}
\end{equation}
To resolve the singularities of the complete intersection 
$\widetilde{W},$
we consider the blow--up $X \to \PP(2q_1, \ldots,2q_n,1,1)$ 
torically defined by introducing the new vector
\begin{equation}
v_{k+n+3}:=\dfrac{v_{k+n+1}+ v_{k+n+2}}{2}=
(-q_1, \ldots, -q_n,0).
\end{equation}
Let $W$ be the smooth blow--up of $\widetilde{W}$ which is a
complete intersection in $X$ of the proper transforms $E_j$
of $\widetilde{E}_j,$ $1 \leq j \leq k.$ The set $\cA \subset \ZZ^{k+n+1}$
is given by $\cA= \{ \ov_1, \ldots , \ov_{k+n+3} \},$
where the vectors $\ov_1, \ldots, \ov_k$ have the first $k$
components equal to the standard basis vectors in $\ZZ^k$ 
and the rest of $n+1$  components equal to zero.
As before, the hypersurfaces $E_1, \ldots, E_k$ correspond to a
partition of the $n+3$ vectors generating the fan defining $X$
in subsets $\{v_{k+1}, \ldots, v_{k+d'_1} \},\ldots,$ 
$\{ v_{k+d'_{k-1}+1}, \ldots, v_{k+n+3}\}.$ The vector $\ov_i$ is obtained 
by adjoining the standard basis vector $e_j$ in $\ZZ^k$ to the vector $v_i$ 
if the divisor in $X$ 
corresponding to $v_i$ contributes to the hypersurface $E_j,$
$1\leq j \leq k.$ Hence, we have
\begin{equation}
\begin{split}
\ov_i & :=(e_i, 0,\ldots,0 ), \ 1 \leq i \leq k,\\
\ov_i & :=(e_j, \bar{v}_{i-k}), \ k+1 \leq i \leq k+n+3, \ 
\hbox{if $v_i$ contributes to $E_j$.}
\end{split}
\end{equation}
The lattice of relations $\LL$ among the elements of $\cA$ is
generated by the relations
\begin{equation}\label{2parL}
\begin{split}
q_1 \ov_{k+1}+ \ldots + q_n \ov_{k+n} + \ov_{k+n+3}
& = d_1 \ov_1+ \ldots + d_k \ov_k, \\
\ov_{k+n+1}+ \ov_{k+n+2} & = 2 \; \ov_{k+n+3}.
\end{split}
\end{equation}
In particular, according to our discussion in section \ref{chap:triang},
the vectors of the Gale transform are given by the columns 
of the $2 \times (k+n+3)$ matrix
\begin{equation}
\left(
\begin{array}{rrrrrrrrr}
-d_1&\ldots&-d_k&q_1&\ldots&q_n&0&0&1\\
0&\ldots&0&0&\ldots&0&1&1&-2
\end{array}
\right)
\end{equation}

Hence,  as a complete fan in $\LL_\RR^\vee \cong \RR^2,$ the secondary 
fan has 4 big cones denoted by
$\Ca_1, \Ca_2, \Ca_3, \Ca_4,$ that are generated 
by the pairs of vectors 
$\{(1,0),(0,1)\},$ $\{(0,1), (-1,0) \},$ $\{(-1,0),(1,-2)\},$ 
and $\{(1,-2), (1,0) \},$ respectively (see figure \ref{fig:dave}).
Denote by $\Ta_1,$ $\Ta_2,$ $ \Ta_3,$ $\Ta_4$ the corresponding 
regular triangulations of the set $\cA.$ The triangulation $\Ta_1$ 
corresponds to the smooth
complete intersection $W \subset X$ (the smooth phase), while
the triangulation $\Ta_3$ corresponds to the Landau--Ginzburg phase.
Proposition \ref{prop:cone}. provides the maximal simplices 
of the four regular triangulations as listed below.
\begin{equation}
\begin{split}
\Ta_1^{n+1}=& \big\{ \{1,\ldots,k+n+3 \} \setminus \{i,j\}, 
\; i \in \{k+1,\ldots,k+n,k+n+3\},\\ 
&  j \in \{k+n+1, k+n+2\} \big\},\\
\Ta_2^{n+1}=& \big\{ \{1,\ldots,k+n+3 \} \setminus \{i,j\},
\; i \in \{1,\ldots, k \}, \\
& j \in \{k+n+1, k+n+2\}  \big\}, \\
\Ta_3^{n+1}=&  \big\{ \{1,\ldots,k+n+2 \} \setminus \{i\},
\; i \in \{1,\ldots, k \} \big\}, \\
\Ta_4^{n+1}=& \big\{ \{1,\ldots,k+n+2 \} \setminus \{i\},
\; i \in \{k+1,\ldots, k+n+2 \} \big\}. 
\label{2par:triang}
\end{split}
\end{equation}
The four phase transitions correspond to relations in $\LL$
as follows 
\begin{equation}\label{2par:trans}
\begin{split}
\Ta_1 \to \Ta_2 \ \hbox{given by}& \ 
q_1 \ov_{k+1}+ \ldots + q_n \ov_{k+n} + \ov_{k+n+3}
= d_1 \ov_1 + \ldots + d_k \ov_k,\\
\Ta_2 \to \Ta_3 \ \hbox{given by}& \ 
\ov_{k+n+1}+ \ov_{k+n+2}  = 2 \;\ov_{k+n+3}, \\
\Ta_1 \to \Ta_4 \ \hbox{given by}& \ 
\ov_{k+n+1}+ \ov_{k+n+2} = 2 \;\ov_{k+n+3}, \\
\Ta_4 \to \Ta_3 \ \hbox{given by}&  \
2q_1 \ov_{k+1}+ \ldots + 2q_n \ov_{k+n} + \ov_{k+n+1}+\ov_{k+n+2}
=\\
&= d_1 \ov_1 + \ldots + d_k \ov_k.
\end{split}
\end{equation}
If $D_j, 1\leq j \leq k+n+3,$ are the toric divisors in $X$ corresponding
to the vectors $v_j,$ than we can choose a basis $\mu, \nu$ of 
${\rm Pic}(X)$ such that
\begin{equation}
\begin{split}
[D_j]=-d_j \mu, \ 1 \leq j \leq k,\ [D_j]=q_j \mu, \ k+1 \leq j \leq k+n,\\
[D_{k+n+1}]= [D_{k+n+2}]= \nu, \ [D_{k+n+3}]=\mu- 2\nu.
\label{2par:linrel}
\end{split}
\end{equation}
The Stanley--Reisner ring of the triangulation
$\Ra_{\Ta_1} \otimes \QQ \cong H^*(X,\QQ)$ is seen to be given by
\begin{equation}
\Ra_{\Ta_1} \otimes \QQ= \QQ[\mu,\lambda] / \Ia_{mon}, \ 
\Ia_{mon}=(\mu^n(\mu-2\lambda), \lambda^2).
\end{equation}
Standard intersection theory in $A^*_\QQ$ gives that
\begin{equation}\label{2par:int}
\big(\prod_{j=1}^n q_j \big) \; \mu^n \lambda= [D_{k+1}] 
\cdot \ldots \cdot [D_{k+n+1}]= [\sigma],
\end{equation}
where $[\sigma]$ is the class in $A^*(X)_\QQ$ given by the maximal
dimensional cone $\sigma$ generated by the vectors $v_1,\ldots,$
$v_{n+1},$ i.e. the class of the corresponding smooth 
point. The class in $A^*(X)_\QQ$ of the smooth complete intersection 
$W$ is $\prod_{j=1}^k d_j \mu,$ and this implies that 
\begin{equation}\label{2par:integral}
\int_W \mu^{n+1-k} = 2 \int_W \mu^{n-k} \nu=
2\frac{\prod_{j=1}^k d_j}{\prod_{j=1}^n q_j}.
\end{equation}
As in the one parameter case, we see that the smoothness of $W$
implies that 
\begin{equation}
{\prod_{j=1}^n q_j} \; \big\vert \; {\prod_{j=1}^k d_j}
\end{equation}
(this follows from assumption (\ref{cond2}), too). We also
have that the ring $\Ha$ (used extensively in the previous sections)
corresponding to  the $(n+1-k)$-dimensional smooth complete intersection
$W$ is given in this case by 
\begin{equation}
\Ha \cong H^*(W,\ZZ) \cong \ZZ[\mu,\lambda] / \Ia_{mon}, \ 
\Ia_{mon}=(\mu^{n-k}(\mu-2\lambda), \lambda^2).
\end{equation}

The complex moduli space of the mirror Calabi--Yau 
variety $M$ has been investigated extensively in the literature in the 
context of studying the mirrors of Calabi--Yau hypersurfaces
in $\PP(2,2,2,1,1)$ and $\PP(6,2,2,1,1).$ Our situation
is more general, but the structure of the discriminant
(the locus of singular complex structures) is identical
up to numerical factors. We will follow very closely the
exposition in \cite{Dave} and we use an identical
pictorial representation of the situation in figure \ref{fig:dave}. 
As we mentioned before, we only consider the simplified 
polynomial moduli space of \cite{AGM2} denoted by $\Ma$ 
whose toric compactification is given by the secondary fan
obtained from the study of the K\"ahler moduli space of
$W$. The variations of the complex structure on $M$ are
described by the series $\Psi_{\mu,\nu}(x,y)$ given by
\begin{equation}\label{2par:psi}
\Psi_{\mu,\nu}(x,y):=\sum_{m,n \geq 0} R_{\mu,\nu}(m,n)
x^{m+\mu} y^{n+\nu},
\end{equation}
with the coefficients $R_{\mu,\nu}(m,n)$ equal to
\begin{equation}
\frac{\prod_{j=1}^k \Gamma(d_j(m+\mu)+1)}
{(\prod_{j=1}^n \Gamma(q_j(m+\mu)+1)) \Gamma^2(n+\nu+1) 
\Gamma((m+\mu)-2(n+\nu)+1)}. 
\end{equation}

\iffigs
\begin{figure}
  \centerline{\epsfxsize=7cm\epsfbox{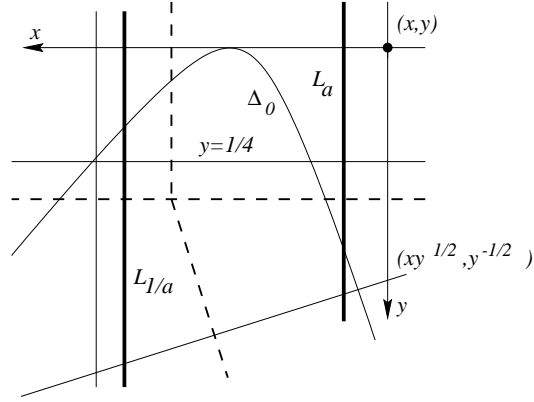}}
  \caption{The complex moduli space $\Ma$ of  $M,$ with the toric 
structure defined by the secondary fan of $W$
(marked with dashed lines).} 
  \label{fig:dave}
\end{figure}
\fi

This series is convergent around the point $(x,y)=(0,0),$
which, in the toric variety $\Ma,$ is the toric fixed point 
corresponding to the cone $\Ca_1.$ Our main task
will be to compute the monodromy of the series $\Psi_{\mu,\nu}(x,y)$
around the so--called principal component of the discriminant
locus, labeled $\Delta_0$ in figure \ref{fig:dave}. To describe the
equation of $\Delta_0$ we use Horn uniformization
(see, for example, page 295 in \cite{GKZ2}, with the observation
that the formulae shown there contain an error). Namely, if we
define the functions
\begin{equation}
F(m,n) := \frac{R_{\mu,\nu}(m,n)}{R_{\mu,\nu}(m+1,n)}, \  
G(m,n) := \frac{R_{\mu,\nu}(m,n)}{R_{\mu,\nu}(m,n+1)} 
\end{equation}
then the principal component of the discriminant $\Delta_0$ 
is parametrized by
\begin{equation}
\begin{split}
x =f(u,v)\quad & ,\quad  \ y  =g(u,v),\ \hbox{where} \\
f(u,v) :=\lim_{t \to \infty}
F(tu,tv)\ &, \  g(u,v) :=\lim_{t \to \infty} G(tu,tv).
\end{split}
\end{equation}

A direct calculation shows that in our context, $\Delta_0$ admits the
rational parameterization
\begin{equation}
x=\prod_{j=1}^k d_j^{\;-d_j}\prod_{j=1}^{n} 
q_j^{q_j} \ \frac{u-2v}{u}, \ y=\frac{v^2}{(u-2v)^2},
\end{equation}
so
\begin{equation}
\Delta_0= \Big\{ \; y= \frac{1}{4} \Big( 1 -  
\prod_{j=1}^k d_j^{\;-d_j}\prod_{j=1}^{n} q_j^{q_j}\; \frac{1}{x} \Big)^2 
\;\Big\}.
\end{equation}
The other components of the discriminant locus are identical with
those from \cite{Dave} and are described by $\{x=0\},$ $\{y=0\},$
$\{x^{-1}=0\}$ and $\{y=1/4\}.$ All of these, as well as the orbifold
locus $\{y^{-1/2}=0\},$ are pictured in figure \ref{fig:dave}.

We now explain the strategy for computing the monodromy around
the principal component $\Delta_0.$ The type of argument to be
used here is typical and goes back to the work of Zariski 
\cite{Zariski} on the fundamental group of the complement of
an algebraic curve in $\CC^2.$ Our first task is to express
the loop around the component $\Delta_0$ as a product
of loops that are allowed by our analytic continuation formulae. 

Fix two small positive real numbers $a$ and
$\theta,$ and consider the complex lines $L_a$ and $L_{1/a}$
(see figure \ref{fig:dave}) defined by
\begin{equation}
L_a:= \{ \; x=a e^{i\theta} \}, \ L_{1/a}:= \{ \; x=a^{-1} e^{i\theta} \},
\end{equation}
The points marked in figure \ref{fig:La} are the points
of intersection of the two complex lines with different
components of the discriminant locus. The loops represented
in the figures are based at the points $(x,y)=(a\; e^{i\theta},b) 
\in L_a$ and $(x,y)=(1/a \; e^{i\theta},b) \in L_{1/a},$ respectively.

\iffigs
\begin{figure}
  \centerline{\epsfxsize=13cm\epsfbox{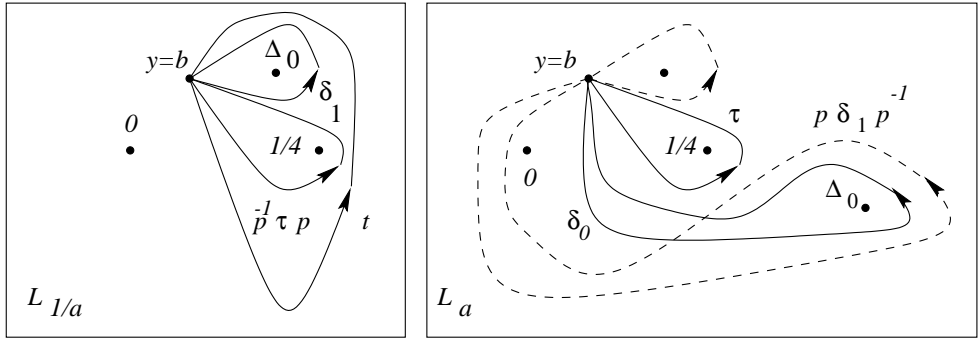}}
  \caption{The loop $\delta_1 \subset L_{1/a}$ around the 
principal component $\Delta_0$ is mapped into the loop
$p \delta_1 p^{-1} \simeq v\cdot \delta_0 \cdot v^{-1} \subset L_a$ under
the continuous transformation from the complex line $L_{1/a}$
to the complex line $L_a$ ($v$ is a counterclockwise
loop in $L_a$ around the origin $y=0$ and $p$ is path in $\Ma$ that
starts at $(a\; e^{i\theta},b)$ and ends at $(1/a\; e^{i\theta},b)$
-- $v$ and $p$ are not shown in the figure).}
  \label{fig:La}
\end{figure}
\fi

Using the $H^*(W,\CC)$--valued series $\Psi_{\mu,\nu}(x,y),$ 
we will show that the monodromy action of the loop $\delta_0$
induces an action on $H^*(W,\CC)$ that 
matches the action on
$H^*(W,\CC)$ induced by the automorphism (\ref{eq:toddcplex}) 
of $\Der(W)$ given by the ideal sheaf of the diagonal in $W \times W.$ 
This means that we have to show that the monodromy action 
can be written as the transformation $Td$ defined by 
\begin{equation}
\label{2par:todd}
I_{\mu,\nu}(x,y) \overset{Td}\longmapsto I_{\mu,\nu}(x,y)- 
(\int_W I_{\mu,\nu}(x,y) \cdot {\rm Todd}_W) \cdot {\rm 1}_W,
\end{equation}
where, as in the one parameter case, the series $I_{\mu,\nu}(x,y)$
is an appropriate version of the series $\Psi$ that
incorporates the solutions to the Picard--Fuchs equations of $M.$

To this end, note that in the fundamental
group of the complement of the discriminant locus,
we have that $t \simeq p^{-1} \cdot \tau \cdot p \cdot \delta_1,$ 
which gives
\begin{equation}\label{eq:path}
\delta_0 \simeq v^{-1} \cdot \tau^{-1} \cdot p \cdot t \cdot p^{-1} \cdot v.
\end{equation}
The loops $\tau$ and $t$ are indicated in figure \ref{fig:La},
and $p$ is the path $p(s)= s\; e^{i\theta},$ $a \leq s \leq 1/a,$
joining the points $(a\; e^{i\theta},b)$ and 
$(1/a\; e^{i\theta},b).$ The nice fact about the relation
(\ref{eq:path}) is that the monodromy action of each path involved on 
the right--hand side is directly expressible using our 
analytic continuation results from section \ref{chap:mon}. 

We will interchangeably use the names of the paths to denote the paths
themselves as well as the
corresponding actions on $H^*(W,\CC)$ (their induced group actions are 
compatible). We prove in fact the following.
\begin{proposition}
The action on $H^*(W,\CC)$ induced by the loop 
$p \cdot t \cdot p^{-1} \cdot v$ is equal to the action induced by
$\tau \cdot v \cdot Td.$
\end{proposition}
\begin{proof} According to the last part of proposition \ref{toricoh}.
we have 
\begin{equation}\label{eq:2partodd}
\Td_W= \big( \prod_{j=1}^k \frac{1- e^{-d_j \mu}}{d_j \mu}\big)
\; \big(\prod_{j=1}^{n} \frac{q_j \mu}{1-e^{-q_j \mu}} \big)
\big( \frac{\nu}{1-e^{- \nu}} \big)^2
\big(\frac{\mu-2\nu}{1-e^{-(\mu-2 \nu)}}\big).
\end{equation}
We will need the following lemma.
\begin{lemma}\label{lemma:2par}
For any cohomology class $\alpha(\mu, \nu) \in H^*(W,\CC)$ 
\begin{equation}
\int_W \alpha \cdot \Td_W = \int_{C_{\mu}} \int_{C_{\nu}}
T(\mu,\nu) \; \alpha(\mu,\nu) \; d\nu d\mu,
\end{equation}
with the function $T(\mu,\nu)$ defined by 
\begin{equation}
T(\mu,\nu):= \frac{\prod_{j=1}^k (1- e^{-2\pi id_j \mu})}
{\prod_{j=1}^{n} (1-e^{-2 \pi i q_j \mu}) (1-e^{-2 \pi i \nu})^2 
(1-e^{-2 \pi i (\mu - 2 \nu)})},
\end{equation}
and the contours $C_\mu,$ $C_\nu$ give the residues at 
$\mu=0$ and $\nu=0$ respectively.
\end{lemma}
\begin{proof}(of the lemma) For any series $\phi(\mu,\nu)$ 
around the origin $(0,0),$ we denote by $c_{m,n}(\phi)$ the coefficient 
of the monomial $\mu^m \nu^n.$ Using the relations (\ref{2par:integral})
in $H^*(W,\CC),$ we can write that
\begin{equation}
\begin{split}
&\int_W \alpha \cdot \Td_W= \frac{\prod_{j=1}^k d_j}{\prod_{j=1}^n q_j}\;
(c_{n-k,1}+2c_{n+1-k,0})(\alpha \cdot \Td_W)=\\
&=(c_{n-k,1}+2c_{n+1-k,0}) \Big( T(\frac{\mu}{2\pi i},\frac{\nu}{2\pi i})
\alpha(\frac{\mu}{2\pi i},\frac{\nu}{2\pi i})
\mu^{n-k}(\mu -2 \nu) \Big)
\label{eq:2par10}
\end{split}
\end{equation}
We can write that

\begin{equation}
\begin{split}
&c_{n-k,1}\Big( T(\frac{\mu}{2\pi i},\frac{\nu}{2\pi i})\;
\alpha(\frac{\mu}{2\pi i},\frac{\nu}{2\pi i})\;
\mu^{n-k}\nu^2(\mu -2 \nu) \Big)=\\
&= \int_{C_{\mu}} \int_{C_{\nu}}
\frac{\mu^{n-k}\nu^2(\mu -2 \nu)}
{\mu^{n-k+1}\nu^2} T(\mu,\nu) \alpha(\mu,\nu)\;
d\nu d\mu=\\
&= \int_{C_{\mu}} \int_{C_{\nu}}
 T(\mu,\nu) \; \alpha(\mu,\nu) \; d\nu d\mu \; - 
2 \int_{C_{\mu}} \int_{C_{\nu}}
\frac{\nu}{\mu}
\; T(\mu,\nu) \; \alpha(\mu,\nu) \; d\nu d\mu,
\end{split}
\end{equation}
and 
\begin{equation}\label{2par:form11}
\begin{split}
&2 c_{n-k+1,0}\Big( T(\frac{\mu}{2\pi i},\frac{\nu}{2\pi i})\;
\alpha(\frac{\mu}{2\pi i},\frac{\nu}{2\pi i})\;
\mu^{n-k}\nu^2(\mu -2 \nu) \Big)=\\
&=2 \int_{C_{\mu}} \int_{C_{\nu}}
\frac{\mu^{n-k}\nu^2(\mu -2 \nu)}
{\mu^{n-k+2}\nu} T(\mu,\nu) \alpha(\mu,\nu) \;
d\nu d\mu=\\
&= 
2 \int_{C_{\mu}} \int_{C_{\nu}}
\frac{\nu}{\mu}
\; T(\mu,\nu) \; \alpha(\mu,\nu)\big) \; d\nu d\mu-
4 \int_{C_{\mu}} \int_{C_{\nu}}
\frac{\nu^2}{\mu^2}
\; T(\mu,\nu) \; \alpha(\mu,\nu) \; d\nu d\mu.
\end{split}
\end{equation}
The lemma follows after pointing out that the last term in 
(\ref{2par:form11}) is in fact zero, since the factor $\nu^2$
cancels the pole at $\nu=0$ of the integrand.
\end{proof}
We proceed with the proof of the proposition. We first investigate
the action of $v \cdot \tau \cdot Td$ on $H^*(W,\CC).$ 
The loop $\tau$ is exactly the loop involved
in the statement of theorem \ref{thm:mon} corresponding to the 
transition $\Ca_1 \to \Ca_4 \to \Ca_1.$ The subtle point here
is that the curve $L_a$ is {\it not} the rational curve 
connecting the toric fixed points $p_1$ and $p_4$ in $\Ma.$ Indeed
as marked in the figure \ref{fig:dave}, if the coordinates around the 
point $p_1$ are $(x,y),$ then the coordinates around the 
point $p_4$ are $(xy^{1/2}, y^{-1/2}).$ As a consequence
of this fact, the analytic continuation of the series $\Psi$
from a neighborhood of $p_1$ to a neighborhood of $p_4$ will
be a series, which, when viewed in the curve $L_a,$ will not converge
in the neighborhood of the point $y=\infty.$ Hence, the loop obtained 
as in theorem \ref{thm:mon} will only go around the point $y=1/4$ and 
not around the point of intersection of $L_a$ with $\Delta_0.$

Using theorem \ref{thm:mon}, as well as (\ref{2par:trans}) and 
lemma \ref{lemma:2par}, we can write explicitly the effect
of the transformation $\tau \cdot v \cdot Td$ on the series 
$\Psi_{\mu,\nu}(x,y)$ (given by (\ref{2par:psi}))

\begin{equation}
\begin{split}
\Psi_{\mu,\nu}(x,y) &\overset{\tau}\longmapsto 
\Psi_{\mu,\nu}(x,y)
 - \int_{C_{\tnu}} e^{2 \pi i (\nu -\tnu)}
\frac{1-e^{2\pi i (\mu-2\nu)}} {(1-e^{-2\pi i \tnu})^2} \;
\Psi_{\mu,\tnu}(x,y)\; d\tnu \\
&\overset{v \cdot Td}\longmapsto
\Big( \delta_{\nu,\tnu} - \int_{C_{\tnu}} e^{2 \pi i (\nu -\tnu)}
\frac{1-e^{2\pi i (\mu-2\nu)}}{(1-e^{-2\pi i \tnu})^2} \Big)\cdot
e^{2 \pi i \tnu}\cdot\\
&\cdot \Big( \psi_{\mu,\tnu}(x,y) - 
\int_{C_{\rmu}} \int_{C_{\rnu}}
 \Psi_{\rmu,\rnu}(x,y) T(\rmu,\rnu)\; d\rnu d\rmu \Big) d\tnu,
\end{split}
\end{equation}
where $\delta_{\nu,\tnu}$ is the Dirac delta function. But
\begin{equation}
\int_{C_{\tnu}} \frac{1}{(1-e^{-2\pi i \tnu})^2}\; d\tnu =1,
\end{equation}
which means that the action of the transformation $\tau \cdot v \cdot Td$
on the $H^*(W,\CC)$--valued series $\Psi_{\mu,\nu}(x,y)$ is
given by
\begin{equation}\label{eq:2par1mon}
\begin{split}
e^{2 \pi i \nu}
&\Big(\Psi_{\mu,\nu}(x,y) - e^{2 \pi i (\mu - 2\nu)} 
\int_{C_{\rmu}} \int_{C_{\rnu}}
\Psi_{\rmu,\rnu}(x,y) T(\rmu,\rnu) d\rnu d\rmu-  \\
&- (1-e^{2\pi i (\mu-2\nu)})
\int_{C_{\tnu}} 
\frac{1} {(1-e^{-2\pi i \tnu})^2} \;
\Psi_{\mu,\tnu}(x,y)\; d\tnu \Big) .
\end{split}
\end{equation}

We now turn our attention to the computation of the action 
on $H^*(W,\CC)$ induced by the loop $p \cdot t \cdot p^{-1} \cdot v.$
As in the one--parameter case we have to work with the 
series $\Phi^{\Ca_1}_{\mu,\nu}(x,y)$ with the 
property (compare (\ref{eq:phipsi})
\begin{equation}\label{2par:phipsi}
\Phi^{\Ca_1}_{\mu,\nu}(e^{i \pi \epsilon} x,y)= 
e^{i \pi (\epsilon + q) \mu}
\prod_{j=1}^k \frac{e^{-2\pi i d_j \mu}-1}
{2 \pi i } \; \Psi_{\mu,\nu} (x,y),
\end{equation}
with $\epsilon$ given by 
\begin{equation}\epsilon:=
\begin{cases}
1& \hbox{if $q$ is odd},\\
2& \hbox{if $q$ is even},
\end{cases}\label{def:2parepsilon}
\end{equation}
The results in section \ref{chap:mon} refer to the $\Phi$--series, and not
directly to the $\Psi$--series. We point out that this fact had
no relevance for the computation of the transformation $\tau \cdot v
\cdot Td$ simply because the variable $\mu$ was not essentially 
involved in that computation.

The expressions for the transformations associated to the 
paths $p$ and $p^{-1}$ are obtained by applying corolary \ref{cor:contform}
to the transitions $\Ca_1 \to \Ca_2$ and $\Ca_2 \to \Ca_1,$
respectively, while the action induced by the loop $t$ is 
given by applying theorem \ref{thm:mon} to the transition 
$\Ca_2 \to \Ca_3 \to \Ca_2.$ In the following group of formulae, we
suppress the variables $x$ and $y.$
\begin{equation}\label{2par:form99}
\begin{split}
\Phi^{\Ca_1}_{\mu,\nu} \overset{p}\longmapsto & \int_{C_{\tmu}}
\frac{e^{\pi i (\epsilon +q)(\mu - \tmu)}}{e^{2\pi i (\tmu -\mu)}-1}
\frac{\prod_{j=1}^k (1-e^{-2\pi i d_j \mu})}
{\prod_{j=1}^k (1-e^{-2\pi i d_j \tmu})} \Phi^{\Ca_2}_{\tmu,\nu}\; d\tmu\\
\Phi^{\Ca_2}_{\tmu,\nu} \overset{t}\longmapsto &
\Phi^{\Ca_2}_{\tmu,\nu} - \int_{C_{\rnu}}
e^{2 \pi i (\nu - \rnu)} \frac{1- e^{2\pi i (\tmu-2\nu)}}
{(1- e^{-2 \pi i \rnu})^2} \Phi^{\Ca_2}_{\tmu,\rnu}\; d\rnu \\
\Phi^{\Ca_2}_{\tmu,\rnu} \overset{p^{-1}\cdot v} \longmapsto &
\int_{C_{\rmu}}
\frac{e^{\pi i (\epsilon +q)(\tmu-\rmu)}}
{e^{2\pi i (\rmu -\tmu)-1 }}
\frac{\prod_{j=1}^n (1-e^{-2\pi i q_j \tmu}) (1-e^{-2\pi i (\tmu -2\rnu)})}
{\prod_{j=1}^n (1-e^{-2\pi i q_j \rmu}) (1-e^{-2\pi i (\rmu -2\rnu)})}\\
& \cdot e^{2\pi i \rnu} \; \Phi^{\Ca_1}_{\rmu,\rnu} \; d\rmu
\end{split}
\end{equation}

Since we are interested in studying the composed transformation
on the series $\Psi,$ we use (\ref{2par:phipsi}) and write 
that the action of the loop $p \cdot t \cdot p^{-1} \cdot v$ on
the series $\Psi_{\mu,\nu}(x,y)$ is given by
\begin{equation}\label{2par:form100}
\begin{split}
\Psi_{\mu,\nu}(x,y)& \mapsto  \int_{C_{\tmu}} 
\frac{1}{\prod_{j=1}^k (1-e^{-2\pi i d_j \tmu})}
\frac{1}{e^{2\pi i (\tmu -\mu)}-1} \cdot \\
& \cdot \Big( \delta_{\nu,\rnu} - \int_{C_{\rnu}}
e^{2 \pi i (\nu - \rnu)} \frac{1- e^{2\pi i (\tmu-2\nu)}}
{(1- e^{-2 \pi i \rnu})^2} \Big) \cdot \\
& \cdot \int_{C_{\rmu}}
\frac{\prod_{j=1}^n (1-e^{-2\pi i q_j \tmu}) (1-e^{-2\pi i (\tmu -2\rnu)})
e^{2 \pi i \rnu}} {e^{2\pi i (\rmu -\tmu)}-1} \cdot \\
& \cdot \frac{\prod_{j=1}^k (1-e^{-2\pi i d_j \rmu})}
{\prod_{j=1}^n (1-e^{-2\pi i q_j \rmu}) (1-e^{-2\pi i (\rmu -2\rnu)})}\;
\Psi_{\rmu,\rnu}(x,y) \;
d\rmu d\rnu d\tmu.
\end{split}
\end{equation}
A first remark is that using the delta function is equivalent
to ignoring the loop $t.$ So the delta function gives 
a transformation corresponding to the analytic continuation
along a contractible loop $p \cdot p^{-1}$ followed 
by the loop $v,$  
\begin{equation}\label{tran1}
\Psi_{\mu,\nu}(x,y) \mapsto e^{2 \pi i \nu} \Psi_{\mu,\nu}(x,y).
\end{equation}  
The contour $C_\tmu$ is chosen according to corolary \ref{cor:contform}
such that it encloses all the values $0,1/d_j, \ldots, (d_j-1)/d_j$
for $1 \leq j \leq k,$ but not the value $\tmu=\mu$ (we dealt
with the subtleties associated with this choice in the proof of
theorem \ref{thm:mon}).  
We first deform the contour such that it encloses
the value $\tmu=\mu,$ and then 
operate the change of variable $t=e^{2 \pi i \tmu}.$
The computation can be finished if we can estimate 
the sum of the residues of the newly obtained function at 
$t=e^{2 \pi i \mu}, 0, \infty.$

As a first step, it is useful to study the residue 
$\tmu=\mu$ by using the $\Phi$--functions instead of $\Psi$--functions,
that is by using the formulae (\ref{2par:form99}). The part
of the transformation that we are interested in is given by 
\begin{equation}
\begin{split}
\Phi^{\Ca_1}_{\mu,\nu}(x,y) \mapsto& - e^{2\pi i \nu} 
(1-e^{2\pi i (\mu-2\nu)})
\int_{C_{\rnu}} \frac{1} {(1-e^{-2\pi i \rnu})^2} \cdot \\
& \cdot \Rz_{\tmu=\mu} \; (p\cdot p^{-1}\cdot \Phi^{\Ca_1}_{\mu, \rnu}(x,y))
 \; d\rnu,
\end{split}
\end{equation}
where $p\cdot p^{-1}$ denotes the action of the contractible loop 
with the same name. But this action is trivial, 
and it is given exactly by the the residue at $\tmu=\mu$ in the 
sequence of changes of variables $(\mu, \rnu) \to (\tmu,\rnu)
\to (\rmu, \rnu)$ (this was the detail that made the proof 
of theorem \ref{thm:mon} work). We conclude that the partial
transformation obtained from the residue at $\tmu=\mu$ is
\begin{equation}\label{tran2}
\Psi_{\mu,\nu}(x,y) \mapsto - e^{2\pi i \nu} 
(1-e^{2\pi i (\mu-2\nu)}) \int_{C_{\rnu}} \frac{1} {(1-e^{-2\pi i \rnu})^2}
\Psi_{\mu, \rnu}(x,y) \; d\rnu,
\end{equation}
with the contour $C_{\rnu}$ chosen to only enclose the pole $\rnu=0.$
 
Finally, we can analyze the last piece of the transformation. 
According to corollary \ref{cor:contform}, the contour $C_\rmu$ needs
to be chosen such that it encloses only the poles
$\rmu=0,1,\ldots,(q_j-1)/q_j,$ for $1 \leq j \leq n,$ and
the value $\rmu=2\rnu.$ But our assumption (\ref{cond2})
guarantees that these points are not poles of the integrand 
(they are canceled out by the product in the numerator),
with the only possible exceptions $\rmu=0$ and $\rmu=2\rnu.$
So we can assume that the contour $C_\rmu$ encloses only
the residues $\rmu=0,2\rnu.$ Since the contour $C_\rnu$
can be chosen to be a small circle with radius $\xi$
around the origin, we can write
\begin{equation}
\begin{split}
&\int_{|\rnu|=\xi}\Big( \int_{|\rmu|=\xi' < \xi}
+ \int_{|\rmu - 2\rnu|=\xi' <2 \xi} f(\rmu,\rnu)d\rmu \Big) d\rnu=\\
=& \int_{|\rmu|=\xi'' > 2\xi} \int_{|\rnu|=\xi} f(\rmu,\rnu) d\rnu d\rmu.
\end{split}
\end{equation}
This means that in formula (\ref{2par:form100}) we can change the 
order of integration. We need to calculate the last
part of the transformation, that is 
\begin{equation}\label{2par:form123}
\begin{split}
&\Psi_{\mu,\nu}(x,y) \mapsto - \!\int_{C_{\rmu}} \!\int_{C_{\rnu}}
\frac{\prod_{j=1}^k (1-e^{-2\pi i d_j \rmu})e^{2\pi i \nu}
\Psi_{\rmu,\rnu}(x,y)}
{\prod_{j=1}^n (1-e^{-2\pi i q_j \rmu}) (1- e^{-2 \pi i \rnu})^2
(1-e^{-2\pi i (\rmu -2\rnu)})} \\
&\cdot \int_{C_{\tmu}} 
\frac{\prod_{j=1}^n (1-e^{-2\pi i q_j \tmu}) (1-e^{-2\pi i (\tmu -2\nu)})
(1-e^{-2\pi i (\tmu -2\rnu)})}
{\prod_{j=1}^k (1-e^{-2\pi i d_j \tmu}) (e^{2\pi i (\tmu -\mu)}-1)
(e^{2\pi i (\rmu -\tmu)}-1)}\; d\tmu d\rnu d\rmu.
\end{split}
\end{equation}
The change of variable $t=e^{2\pi i \tmu}$ gives 
\begin{equation}
\begin{split}
&\int_{C_{\tmu}} 
\frac{\prod_{j=1}^n (1-e^{-2\pi i q_j \tmu}) (1-e^{2\pi i (\tmu -2\nu)})
(1-e^{-2\pi i (\tmu -2\rnu)})}
{\prod_{j=1}^k (1-e^{-2\pi i d_j \tmu}) (e^{2\pi i (\tmu -\mu)}-1)
(e^{2\pi i (\rmu -\tmu)}-1)}\; d\tmu=\\
&=\int_{C_t}
\frac{\prod_{j=1}^n (1- t^{-q_j}) (1-t e^{-2 \pi i(2 \nu)})
(1-t^{-1} e^{4 \pi i \rnu})}
{\prod_{j=1}^k (1-t^{-d_j}) (t e^{- 2 \pi i \mu} -1)
(t^{-1} e^{ 2 \pi i \rmu} -1)}\; \frac{dt}{t}.
\end{split}
\end{equation}
Since $d_1+\ldots+d_k=q_1+\ldots+q_k+1,$ direct power counting shows
that $t=0$ is not a pole. By computing the value of the residue
at $t=\infty$ the value of the integral is obtained to be 
$e^{2 \pi i (\mu - 2 \nu)},$ and, after using lemma \ref{lemma:2par},
the transformation (\ref{2par:form123}) can be written as
\begin{equation}\label{tran3}
\Psi_{\mu,\nu}(x,y) \mapsto - e^{2\pi i \nu} e^{2 \pi i (\mu - 2 \nu)}
\int_{C_{\rmu}} \int_{C_{\rnu}} \Psi_{\rmu,\rnu}(x,y)
T(\rmu, \rnu) \; d\rnu d\rmu.
\end{equation}

By adding up the partial transformations (\ref{tran1}), (\ref{tran2}), 
and (\ref{tran3}), we conclude that the action 
of the loop $p \cdot t \cdot p^{-1} \cdot v$ on 
the $H^*(W,\CC)$--valued series $\Psi_{\mu,\nu}(x,y)$
is given indeed by (check with (\ref{eq:2par1mon}))
\begin{equation}
\begin{split}
& e^{2 \pi i \nu} \Psi_{\mu,\nu}(x,y)
- e^{2\pi i \nu} 
(1-e^{2\pi i (\mu-2\nu)}) \int_{C_{\rnu}} \frac{1} {(1-e^{-2\pi i \rnu})^2}
\Psi_{\mu, \rnu}(x,y) \; d\rnu - \\
&- e^{2 \pi \nu} e^{2 \pi i (\mu - 2 \nu)}
\int_{C_{\rmu}} \int_{C_{\rnu}}\Psi_{\rmu,\rnu}(x,y)
 T(\rmu, \rnu) \; d\rnu d\rmu,
\end{split}
\end{equation}
which ends the proof of the proposition.
\end{proof}
We summarize the results of this section below.

\begin{thm} \label{thm:2par}
\
\begin{itemize}
\item The kernel $p_2^*(\La_\mu) \otimes j_* {\cal O}_{\Delta}$ 
where the line bundle $\La_\mu$ has $c_1(\La_\mu)= \mu,$
corresponds to the loop $u$ which goes once around the
component $x=0$ of the discriminant locus.
\item The kernel $p_2^*(\La_\nu) \otimes j_* {\cal O}_{\Delta}$ 
where the line bundle $\La_\nu$ has $c_1(\La_\nu)= \nu,$
corresponds to the loop $v$ which goes once around the
component $y=0$ of the discriminant locus.
\item The kernel $\Ia_\Delta$ given by the ideal sheaf of the 
diagonal in $W \times W$ (or, equivalently defined by the complex
(\ref{eq:toddcplex})) corresponds to the loop $\delta_0$ which
goes once around the ``closest'' part (to the point $(x,y)=0$) 
of the principal component $\Delta_0$ of the discriminant locus
(see figure \ref{fig:dave}).
\item The kernel (\ref{eq:cplex}) defined by the complex 
\begin{equation}
\ldots \to 0 \to p_2^*(\La_\nu) \otimes p_1^* ({\cal L_\nu^{\vee}}) 
\to j_* {\cal O}_{\Delta} \to 0 \to \ldots,
\end{equation}
with $j_* {\cal O}_{\Delta}$ at the $0$-th place and
$\La_\nu$ as above, corresponds to the loop 
$v \cdot \delta_0 \cdot v^{-1}$ which is homotopically 
equivalent to a loop that goes once 
around the ``farthest'' part (from the point $(x,y)=(0,0)$)
of the principal component $\Delta_0$ of the discriminant locus
(this is the loop denoted by $p \cdot \delta_1 \cdot p^{-1}$ in 
the explanation to figure \ref{fig:dave}).
\end{itemize}
\end{thm}

To complete the dictionary between kernels and loops for this 
situation we still need to add to the above list the kernel corresponding
to the loop around the component $y=1/4.$ In the next section,
we will deal with this issue in a much more general context.

\subsection{The General Case}\label{chap:generalcase}

In this section we prove that the action in cohomology of the kernel 
$\Ea^\bullet(F)$ introduced in the definition \ref{def:bigdef}
corresponding to an edge 
$F$ in the secondary polytope does indeed match the monodromy action 
given by theorem \ref{thm:mon}. The edge $F$
is assumed to have one vertex corresponding to a triangulation $\Ta_0$
giving a smooth phase of the Calabi--Yau complete intersection $W,$ and the 
other vertex corresponding to a triangulation $\Ta_1.$
The edge corresponds to two adjacent big cones in the secondary fan, and 
defines a closed toric orbit which is a projective curve $\PP^1$ embedded 
in the moduli space of complex structures of the mirror Calabi--Yau $M$
whose toric compactification is given by the secondary fan. 

We return to the general toric context presented in the beginning
of this chapter, so 
$W=E_1 \cap \ldots \cap E_k \subset X_{\Sigma_0}$ is an 
$(n-k)$-dimensional 
Calabi--Yau complete intersection determined by a nef--partition
in the $n$-dimensional toric Fano variety $X_{\Sigma_0}.$ 
The elements of $\Sigma_0(1)$ are denoted by $\{ v_{k+1}, \ldots, v_N \},$
and the corresponding $\cA$--set in $\ZZ^{n+k}$ is given by 
$\ovl{\Sigma}_0(1)= 
\{ \ov_1, \ldots, \ov_k, \ov_{k+1}, \ldots, \ov_N \}.$
The zero subscript indicates that the fan $\overline{\Sigma}_0$ is obtained 
from the triangulation $\Ta_0.$ 
Assume that the edge $F$ in the secondary polytope is
associated to a circuit $I$ corresponding to a relation of the form
\begin{equation}\label{eq:circuit}
\sum_{j \in I_+} q_j \ov_j= \sum_{j \in I_-'} q_j \ov_j +
\sum_{j \in I_-''} d_j \ov_j,
\end{equation}
with $(I_+ \cup I_-') \subset \{k+1, \ldots, N\},$ 
$I_-'' \subset \{1, \ldots , k\},$ $q_j, d_j > 0.$ The circuit $I$
defines the modification $s_I(\Ta_0)= \Ta_1$ which replaces the 
simplices of type Conv($I \setminus j),$ $j \in I_-' \cup I_-'',$ in 
$\Ta_0$ with simplices of type Conv($I \setminus j),$ $j \in I_+,$ in
$\Ta_1.$ According to proposition \ref{prop:edge} the polyhedral 
subdivision $\Sa(F)= \Sa(\Ta_0, \Ta_1)$ contains all the 
common simplices of $\Ta_0$ and $\Ta_1,$ as well as polyhedra
of type Conv($I \cup J$) for all separating subsets 
$J$ (definition \ref{def:sep}). The triangulations $\Ta_0, \Ta_1$ 
as well as the polyhedral subdivision $\Sa(F),$ determine different fan
structures denoted by $\ovl{\Sigma}_0, \ovl{\Sigma}_1$ and $\ovl{\Sigma}_F,$
respectively. These fans are all supported on the same reflexive 
Gorenstein cone generated by the $\cA$--set, and they define different 
non-compact toric varieties $\EE_{\Ta_0}, \EE_{\Ta_1},$ and $\EE_{\Sa(F)},$ 
respectively, subject to the diagram
\begin{equation}
\begin{split}
\xymatrix{
\EE_{\Ta_0}\ar[rd]^{f}&             & \EE_{\Ta_1}\ar[ld]_{g} \\
                   & \EE_{\Sa(F)} &                   }
\end{split}
\end{equation} 

We would like to understand the exceptional loci of the map 
$f$ and its restrictions to $X_{\Sigma_0}$ and $W.$ The map
$f$ comes from a map of fans (see section \ref{chap:toricgen})
which is simply the identity at the lattice level. The construction
of the fan $\ovl{\Sigma}_F$ described above shows that in passing 
from $\ovl{\Sigma}_0$ to $\ovl{\Sigma}_F,$ the cone $\ovl{\sigma}_0$ 
in $\ovl{\Sigma}_0,$ with 
\begin{equation}\label{eq:exc1}
\ovl{\sigma}_0=\sum_{j \in I_-' \cup I_-''} \RR_{\geq 0} \ov_j,
\end{equation}
is mapped as a strict subcone of the cone $\ovl{\sigma}_F$ in 
$\ovl{\Sigma}_F,$ where
\begin{equation}\label{eq:exc2}
\ovl{\sigma}_F=\sum_{j \in I} \RR_{\geq 0} \ov_j.
\end{equation}
Since $\ovl{\sigma}_F$ is the minimal cone in $\ovl{\Sigma}_F$ 
that contains $\ovl{\sigma}_0,$
it follows that the closed toric orbits 
$V(\ovl{\sigma}_0) \to V(\ovl{\sigma}_F)$ are mapped onto and they are the
loci where the map $f : \EE_{\Ta_0} \to \EE_{\Ta_1}$ is not an isomorphism. 

The compact toric variety $X_{\Sigma_0},$ embedded
as the zero section in $\EE_{\Ta_0},$ is mapped onto
a compact toric subvariety $X_{\Sigma_F}$ in $\EE_{\Sa(F)}.$ The 
fan $\Sigma_F$ is obtained by performing a modification 
(as described above) to the fan $\Sigma_0$ along the circuit 
(\ref{eq:circuit1}) given below.
We collect these facts and their direct consequences 
in the following proposition.
\begin{proposition}\label{prop:exclocus}
\
\begin{enumerate}
\item The loci where the map $f : \EE_{\Ta_0} \to \EE_{\Sa(F)}$
fails to be an isomorphism are the closed toric orbits 
$V(\ovl{\sigma}_0) \subset \EE_{\Ta_0}$ and 
$V(\ovl{\sigma}_F) \subset \EE_{\Sa(F)},$ with 
$\ovl{\sigma}_0$ and $\ovl{\sigma}_F$ cones in $\ovl{\Sigma}_0$ and 
$\ovl{\Sigma}_F$ defined by (\ref{eq:exc1}) and (\ref{eq:exc2}), 
respectively. We have that 
\begin{equation}
V(\ovl{\sigma}_0)=(\bigcap_{j \in I_-'}
\ovl{D}_j) \cap (\bigcap_{j \in I_-''} \ovl{E}_j),
\end{equation}
where $\ovl{D}_j \subset \EE_{\Ta_0}$ are the toric divisors 
determined by the vectors $\ov_j, j \in j \in I_-'.$
The restriction of $f$ maps $V(\ovl{\sigma}_0)$ onto
$V(\ovl{\sigma}_F).$
\item The exceptional locus $E_0 \subset X_{\Sigma_0}$ 
of the restriction $f_{\vert  X_{\Sigma_0}}$ is irreducible and 
torically described by the closed orbit $V(\tau),$ where $\tau$ is the cone 
in $\Sigma_0$ generated by the vectors $v_j, j \in I_-'.$  We have that
\begin{equation}\label{eq:e0}
E_0=V(\tau)=\bigcap_{j \in I_-'} D_j,
\end{equation}
where $D_j \subset X_{\Sigma_0}$ are the toric divisors 
determined by the vectors $v_j, j \in j \in I_-'$ (if 
$I_-'=\emptyset,$ then $E_0= X_{\Sigma_0}$). The relation 
(\ref{eq:circuit}) gives rise to a relation between the elements of 
$\Sigma_0(1)$ which reads
\begin{equation}\label{eq:circuit1}
\sum_{j \in I_+} q_j v_j= \sum_{j \in I_-'} q_j v_j. 
\end{equation}
\item The toric structure of $E_0$ is determined by a fan included 
in the lattice $(\Nrm/ \Nrm_\tau)_{\RR},$ ($\Nrm_\tau$ is the 
sublattice in $\Nrm$ generated by $\Nrm \cap \tau$)
consisting of all the cones in $\Sigma_0$ that contain $\tau$ as 
a subcone. More explicitly, the maximal dimensional (simplicial) cones of
the fan defining $E_0$ are generated by $v_i'$ 
(the images of the vectors $v_i$ under the projection 
$\Nrm \to \Nrm/ \Nrm_\tau$) with
$i \in ((I_+ \setminus j) \cup J) \setminus \{1, \ldots , k\},$
for all $j , j \in I_+,$ and $J$ a separating set for $\Ta_0$ and
$\Ta_1.$ In particular, the relation (\ref{eq:circuit1}) 
becomes
\begin{equation} 
\sum_{j \in I_+} q_j v_j' = 0.
\end{equation}
By replacing all the cones
of the form Conv$(I_+ \setminus j)$ with the cone Conv$(I_+),$
the fan defining $E_0$ becomes the fan defining $Z_0=f(E_0).$ 
Any fiber of $f: E_0 \to Z_0$ is a weighted projective space
of dimension equal to $\dim E_0 - \dim Z_0.$

\end{enumerate}
\end{proposition}
 
In fact, the last property shows that the lattice 
$\Nrm/ \Nrm_\tau$ admits a decomposition $\Nrm/ \Nrm_\tau= 
\Nrm_1 \oplus \Nrm',$ where $\Nrm'$ is the sublattice generated 
in $\Nrm/ \Nrm_\tau$ by the vectors $v_j', j \in I_+.$ 
In order to construct the fibered product $E_0 \times_{Z_0} E_0,$
consider the lattice $\Nrm_1 \oplus \Nrm' \oplus \Nrm'',$ with 
$\Nrm''$ isomorphic to $\Nrm'.$ For each vector $v_j', j \in I_+,$ denote by 
$v_j''$ its image in $\Nrm''$ under the isomorphism $\Nrm' \cong \Nrm''.$ 
We have that
\begin{equation}
\sum_{j \in I_+} q_j v_j'' = 0.
\end{equation}
Given a ray $v=(v_1, v')$ in the fan defining $E_0,$ the 
corresponding ray in the fan defining $E_0 \times_{Z_0} E_0$
is given by $\tv=(v_1,v',v'').$ As a direct consequence of this 
construction we obtain the following.

\begin{proposition} \label{prop:fanprod}
The maximal cones of the (simplicial)
fan in $(\Nrm_1 \oplus \Nrm' \oplus \Nrm'')_\RR$ defining the toric variety 
$E_0 \times_{Z_0} E_0$ have the form 
\begin{equation}\label{eqn:fiberedfan}
\RR_{\geq 0} \tv_i + \sum_{j \in I_+ \setminus j'} \RR_{\geq 0} v_j' +
\sum_{j \in I_+ \setminus j''} \RR_{\geq 0} v_j'',
\end{equation}
for some $j',j'' \in I_+,$ and some vectors 
$v_i$ and $v_j', j \in I_+ \setminus j',$ which generate a maximal cone 
in the fan defining $E_0.$ 
The maps of fans given by the diagram 
\begin{equation}
\begin{split}
\xymatrix{
 &\Nrm_1 \oplus \Nrm' \oplus \Nrm'' \ar[ld]\ar[rd] & \\
\Nrm_1 \oplus \Nrm' \ar[rd]& & \Nrm_1 \oplus \Nrm'' \ar[ld] \\
 & \Nrm_1 & 
}
\end{split}
\end{equation}
correspond to the canonical maps in the diagram  
\begin{equation}
\begin{split}
\xymatrix{
 &E_0 \times_{Z_0} E_0 \ar[ld]\ar[rd] & \\
 E_0 \ar[rd] & & E_0 \ar[ld] \\
 & Z_0 &
}
\end{split}
\end{equation}
\end{proposition}

\begin{remark} \label{rem:mori}
The toric morphism $f : X_{\Sigma_0} \to X_{\Sigma_F}$ 
is an elementary contraction in the sense of Mori theory
as discussed by M. Reid \cite{reid} in the toric context. 
This interpretation will shed some light on our construction,
so we detail it here. 

First, recall that the closed K\"ahler cone of $X_{\Sigma_0}$ consisting 
of positive $(1,1)$-classes in $H^2(X_{\Sigma_0}, \RR)$ 
is isomorphic to the big cone in the secondary 
fan corresponding to the triangulation $\Ta_0$ (see theorem 4.5
in \cite{batqcoh}). Since $X_{\Sigma_0}$ is simplicial and complete,
the Mori cone of $X_{\Sigma_0}$ is the cone of effective $1$-cycles in 
$A_1(X_{\Sigma_0}) \otimes \RR \cong H_2(X_{\Sigma_0},\RR),$ and it is 
dual to the K\"ahler cone. The Mori cone can be then naturally 
identified with the ``big'' cone (defined by (\ref{eq:defncone}))
at the vertex in the secondary polytope corresponding to the 
triangulation $\Ta_0.$ In particular, the edges of this cone
give rise to extremal rays. For a concise review of M. Reid's
work \cite{reid} on Mori theory for toric varieties, see also section 
3.3.2 in \cite{CK}. 

We only review here the facts that will be used in this work. According
to Reid's results, for each edge $F$ as above, 
the class of the extremal ray in $H_2(X_{\Sigma_0}, \RR)$
is given by the closed toric orbit corresponding to some internal
$(n-1)$-dimensional cone in the $n$-dimensional simplicial fan 
$\Sigma_0.$ If $v_{j_1}, \ldots, v_{j_{n-1}}$ are the generators 
of such a cone $\sigma,$ there exist two other generators 
$v_{j_n}$ and $v_{j_{n+1}}$ such that $v_{j_1}, \ldots, v_{j_{n-1}}, v_{j_n}$ 
and $v_{j_1}, \ldots, v_{j_{n-1}}, v_{j_{n+1}}$ generate the two 
corresponding adjacent cones of dimension $n.$ It follows that there is a 
relation 
\begin{equation}\label{eq:relm1}
\sum_{i=1}^{n+1} \lambda_{j_i} v_{j_i} =0,
\end{equation}
with $\lambda_{j_n}, \lambda_{j_{n+1}} > 0.$ The cone $\sigma$ 
determines a curve $C_{\sigma} \subset X_{\Sigma_0}$ which, 
according to proposition \ref{katzcox}, determines a relation 
\begin{equation}\label{relm2}
\sum_{j=k+1}^N (D_j \cdot C_{\sigma}) v_j=0,
\end{equation}
where $D_j$ is the toric divisor corresponding to $v_j,$ and, in line 
with our previous notation, $v_{k+1}, \ldots, v_N$ are the elements
of $\Sigma_0(1).$
 
As stated in lemma 3.3.2 in \cite{CK}, and as a direct consequence of
Reid's work, the two relations are in fact
identical up to a multiplication by an integral constant. More
important for us, they are also identified with the relation 
(\ref{eq:circuit1}) defining the circuit $I_+ \cup I_-'$ (the 
restriction of the circuit $I_+ \cup I_-' \cup I_-''$ to the
elements of $\Sigma_0(1)$). The vectors $v_{j_n}$ and $v_{j_{n+1}}$
have the property $\{ j_n, j_{n+1} \} \subset I_+.$

\begin{proposition} (Cor. 2.5 and 2.6 and Prop. 2.7 in \cite{reid})
\label{prop:contraction}
\
\begin{enumerate}
\item The elementary contraction corresponding to the class of the 
curve $C_\sigma$ is the map $f: X_{\Sigma_0} \to X_{\Sigma_F}.$ The 
restriction $f: E_0 \to Z_0$ to the loci where $f$ is not an 
isomorphism, is a flat morphism, all of whose fibers 
are weighted projective spaces of dimension $\dim(E_0) - \dim(Z_0).$
\item 
\begin{equation}\label{eq:int}
\begin{split}
D_j &\cdot C_{\sigma} > 0, \quad \hbox{if} \; j \in I_+, \\
D_j &\cdot C_{\sigma} < 0, \quad \hbox{if} \; j \in I_-', \\ 
D_j &\cdot C_{\sigma} = 0, \quad \hbox{if} \; j \notin I_+ \cup I_-'.
\end{split}
\end{equation}
\end{enumerate}
\end{proposition}

\end{remark}

As described in section \ref{chap:functors}, the definition \ref{def:bigdef} 
of the complex $\Ea^\bullet(F)$ involves in an 
essential way the subvariety $Y= E \times_{Z} E,$ where 
$E= W \cap E_0$ and $Z=f(E).$ Using the 
results of the previous remark, we can give an explicit 
description of the subvariety $Y$ as a complete intersection in the 
toric variety $E_0 \times_{Z_0} E_0.$

\begin{proposition} \label{prop:yconstr}
If $E_1, \ldots, E_k$ are the Cartier divisors in
$X_{\Sigma_0}$ given by the nef--partition defining the Calabi--Yau
manifold $W=E_1 \cap \ldots \cap E_k,$ ($I_-'' \subset \{1 ,\ldots, k \}$),
then 
\begin{equation}
Y=         \bigcap_{j \in I_-''} f_1^*(E_j \cap E_0) \cap
           \bigcap_{j \in I_-''} f_2^*(E_j \cap E_0) \cap
           \bigcap_{j \notin I_-''} E_j \cap (E_0 \times_{Z_0} E_0),
\end{equation} 
where $f_1$ and $f_2$ are the toric maps
\begin{equation}
\begin{split}
\xymatrix{
 &E_0 \times_{Z_0} E_0  \ar[ld]_{f_2} \ar[rd]^{f_1} & \\
 E_0 & & E_0
}
\end{split}
\end{equation}
and, for $j \notin I_-'',$ we have that
\begin{equation}\label{eq:pullequal}
E_j \cap (E_0 \times_{Z_0} E_0)=f_1^*(E_j \cap E_0) =f_2^*(E_j \cap E_0).
\end{equation}
\end{proposition}

\begin{proof} The Cartier divisors (spanned by global sections)
$E_1, \ldots, E_k$ come from a nef--partition of the vectors in $\Sigma_0(1).$
As a consequence of the second part of proposition \ref{prop:contraction},
for any divisor $E_j,$ with $j \in \{1, \ldots ,k \} \setminus I_-'',$
we have that $E_j \cdot C_\sigma =0.$ 
This shows that any curve that is contracted under 
the map $f : X_{\Sigma_0} \to X_{\Sigma_F}$ does not intersect the 
generic divisors $E_j$ for $j \in \{1, \ldots ,k \} \setminus I_-''.$
Proposition \ref{prop:fanprod} implies then that $E_0 \cap f_1^*(D_i)= 
E_0 \cap f_2^*(D_i),$ for any $i \in \{k+1, \ldots , N \} \setminus 
(I_+ \cup I_-').$ Hence, for any $j \notin I_-'',$ we have that
$f_1^*(E_j \cap E_0)= f_2^*(E_j \cap E_0).$ Note also that,
for $j \in I_-'',$ $E_j \cdot C_\sigma > 0,$ and 
the result follows. 

\end{proof}

We are now in the position to analyze the action of the complex
$\Ea^\bullet(F)$ on $H^*(W,\CC)$ in order to make the connection with the 
monodromy computations. 
One of the key ingredients will be the 
computation of the Gysin map $H^*(Y,\CC) \to H^*(E,\CC)$ induced by the 
canonical projection $Y=E \times_{Z} E \to E.$ We need the following 
technical (but crucial) result.

\begin{proposition} \label{prop:toricfibr}
Consider a toric fibration $p: X_{\Sigma} \to 
X_{\Sigma_1}$ induced by the map of fans $\Nrm_1 \oplus \Nrm_2 \to \Nrm_1.$
Assume that the typical fiber of $p$ is a weighted projective space
$\PP(h_1, h_2, \ldots, h_{l+1})$
corresponding to a set of vectors $v_1, \ldots, v_{l+1}$ in $\Sigma(1),$
with $\{ v_1, \ldots, v_{l+1} \} \subset \{0\} \oplus \Nrm_2,$ 
which satisfy the linear relation
\begin{equation}\label{eq:sum0rel}
h_1 v_1 + h_2 v_2 + \ldots + h_{l+1} v_{l+1} =0.
\end{equation}

If $\mu \in H^2(X_\Sigma, \ZZ)$ denotes the generator in cohomology
determined by this relation, then, as in section \ref{chap:mon}, 
we assume that $[D_j]= h_j \mu + \mu_j', 
1 \leq j \leq l,$ where $D_j$ is the toric divisor in $X_{\Sigma}$
induced by the vector $v_j.$ Then, for any 
$\gamma \in H^*(X_{\Sigma},\CC),$ we have that the Gysin map
$p_*: H^*(X_\Sigma,\CC) \to H^{*-l}(X_{\Sigma_1},\CC)$ is given by 
\begin{equation}\label{eq:intgamma}
p_* (\gamma)= \frac{1}{2\pi i}
\int_{C_\xi} \frac{1}{\prod_{j=1}^{l+1} (h_j \xi + \mu_j')}
\, \gamma(\xi) \,d\xi,
\end{equation}
where it is understood that $\gamma$ is written in terms of the same 
generators in cohomology used to express the classes $\mu, \mu_j',$ 
and $C_\xi$ is a contour enclosing all the poles of the integrand.
\end{proposition}

\begin{proof} The fan defining the weighted projective space 
$\PP(h_1, h_2, \ldots, h_{l+1}),$ has the property that
the multiplicity index in $\Nrm_2$ of any maximal cone generated by the 
vectors $v_j, j \in \{1, 2, \ldots, l+1 \} \setminus j_0,$ is equal to 
$h_{j_0}.$ 

Since the vectors $v_1, v_2, \ldots, v_{l+1}$ do not
belong to any maximal cone of the fan $\Sigma_0,$ we have that
\begin{equation}\label{eq:zeroint}
[D_1] \cdot [D_2] \ldots [D_{l+1}] =0.
\end{equation}

For each $i, 1 \leq i \leq l+1,$ the closed toric orbits $V_i,$ 
given by the cone with generators
$v_j, j \in \{1, \ldots, l+1 \} \setminus i,$ 
is mapped under the morphism $p$ to the 
closed toric orbit in $X_{\Sigma_1}$ represented by $X_{\Sigma_1}$ itself. 
Remark 10.9 in \cite{Dani} shows that, in $A_*(X_{\Sigma_0})_\QQ,$ we have 
that
\begin{equation}
h_i \prod_{j=1,j \not= i}^k [D_j]= [V_i].
\end{equation}
We obtain that, for any $i, 1 \leq i \leq l+1,$
\begin{equation}\label{eq:gysmap}
p_*\bigl( h_i \prod_{j=1, j \not= i}^{l+1} [D_j] \bigr)= 1 \in 
H^0(X_{\Sigma_1}, \CC).
\end{equation}
According to lemma 10.7.1 in \cite{Dani}, any class $\gamma \in 
H^*(X_\Sigma, \CC)$ can be written as a sum of ``monomials'' without 
repetitions of classes corresponding to toric divisors. To prove the 
result, we only have to consider the case when $\gamma$ is 
such a monomial. Assume that
\begin{equation}\label{eq:exgamma}
\gamma=\prod_{j \in J} [D_j] \cdot \prod_{j \in J'} [D_j],
\end{equation}
where $J \subset \{1,2 \ldots, l+1\},$ and 
the divisors $D_j, j \in J',$ correspond to vectors in
$\Sigma(1)$ other than $v_1, v_2,\ldots, v_{l+1}.$ For any
$j \notin \{1,2, \ldots, l+1 \}$, there is no danger 
of confusion if we also denote by $[D_j]$ the class of the toric divisor in 
$X_{\Sigma_1}$ whose pull-back is the class $[D_j]$ in $X_\Sigma.$
Hence 
\begin{equation}
p_*(\gamma)= \prod_{j \in J'} [D_j] \cdot 
p_* \bigl(\prod_{j \in J} [D_j]\bigr).
\end{equation}
If $\vert J \vert < l,$ then $p_* \bigl(\prod_{j \in J} [D_j]\bigr) =0,$
and relation (\ref{eq:zeroint}) shows that the only case when the 
result is non-zero is the case when $J= \{1,2, \ldots, l+1\} \setminus i,$ for 
some $i.$ In that case, we have 
\begin{equation}
p_*(\gamma)=\frac{1}{h_i} \prod_{j \in J'} [D_j].
\end{equation}   

Consider now a system of generators of $H^*(X_\Sigma, \ZZ)$ corresponding
to a basis of the space of linear relations among the vectors in $\Sigma(1),$
chosen such that $\mu$ denotes the generator corresponding to the relation 
(\ref{eq:sum0rel}). We can write that
\begin{equation}
[D_j]=h_j \mu + \mu_j', \; 1 \leq j \leq l+1,
\end{equation}
while the expressions for the classes $[D_j], j \notin \{1,2, \ldots , l+1 \},$
do not involve $\mu.$ If the monomial $\gamma$ involves all the classes
$[D_1], [D_2], \ldots, [D_{l+1}],$ then the integrand in the equation
(\ref{eq:intgamma}) does not have any residue, so the value of the integral is
indeed zero. Also, if the subset $J$ in (\ref{eq:exgamma}) has strictly
less than $l$ elements, the sum of the residues at finite poles calculated
by the right hand side of the formula (\ref{eq:intgamma}) is again zero, simply
because the residue at $\infty$ is zero. The remaining case is again 
$J= \{1,2, \ldots, l+1\} \setminus i,$ for some $i.$ But then 
\begin{equation}
\begin{split}
&\frac{1}{2\pi i}
\int_{C_\xi} \frac{1}{\prod_{j=1}^{l+1} (h_j \xi + \mu_j')}
\, \gamma(\xi) \,d\xi  \\
&= \frac{1}{2\pi i} \int_{C_\xi} \frac{1}{h_i \xi + \mu_i'}\,d\xi 
\; \prod_{j \in J'} [D_j] \\
&= \frac{1}{h_i} \prod_{j \in J'} [D_j],
\end{split}
\end{equation}
which ends the proof of the proposition.

\end{proof}

By applying this proposition to the toric map $E_0 \times_{Z_0} E_0
\to E_0,$ we can compute the Gysin map induced by the morphism
$Y \to E$ using the following immediate lemma.

\begin{lemma}\label{lemma:map}
Given the diagram 
\begin{equation}
\begin{split}
\xymatrix@C+10pt{
Y\; \ar@{^{(}->}[r]^-{i} \ar[d]_{q_2} & 
E_0 \times_{Z_0} E_0 \ar[d]_{f_2}\\
E\; \ar@{^{(}->}[r]^{j} & E_0
}
\end{split}
\end{equation}
with $i$ and $j$ the canonical inclusions, and 
$\gamma \in H^*(E_0 \times_{Z_0} E_0, \CC),$ then
\begin{equation}
(q_2)_* (i^*(\gamma)) = j^*(\gamma'),
\end{equation}
for some $\gamma' \in H^*(E_0, \CC)$ with the property
\begin{equation}
\gamma' \cdot [E]= (f_2)_* (\gamma \cdot [Y]).
\end{equation}
\end{lemma}
 
\begin{proof} Indeed, we have that
\begin{equation}
j_*(q_2)_* (i^*(\gamma))= (f_2)_*  i_* (i^*(\gamma)) 
= (f_2)_* (\gamma \cdot [Y]), 
\end{equation}
so  
\begin{equation}
\gamma' \cdot [E]= j_* (j^*(\gamma')= (f_2)_* (\gamma \cdot [Y])
\end{equation}

\end{proof}

We can now put things together and compute the action of the complex
$\Phi_{\Ea^\bullet (F)}$ on $H^*(W,\CC).$ As before, we assume that we have
chosen a system of generators for the cohomology of the toric variety
$X_{\Sigma_0}$ such that
the class corresponding to the relation given by the edge $F$ is
denoted by $\mu.$. In particular, if the relation among the elements
of the corresponding set $\cA=\{ \ov_1, \ldots ,\ov_k, \ov_{k+1}, 
\ldots, \ov_N \}$  is written as
\begin{equation}\label{eq:circuitt}
\sum_{j \in I_+} q_j \ov_j= \sum_{j \in I_-'} q_j \ov_j +
\sum_{j \in I_-''} d_j \ov_j,
\end{equation}
with $I_+ \cup I_- \subset \{k+1, \ldots, N \},$ and 
$I_-'' \subset \{1, \ldots, k \}.$ The classes of the toric divisors in 
$X_{\Sigma_0}$ look like this:
\begin{equation}
\begin{split}
[D_j]&= q_j \mu + \mu_j', \; j \in I_+,\\
[D_j]&= -q_j \mu + \mu_j', \; j \in I_-',\\
[E_j]&= d_j \mu - \mu_j', \; j \in I_-'',
\end{split}
\end{equation}
and $W= E_1 \cap E_2 \ldots \cap E_k.$
In what follows, for the ease of notation, we will sometimes
omit to write the obvious pull-back maps of the classes in 
$H^*(X_{\Sigma_0}, \CC)$ to the various sub-varieties involved.

\begin{proposition}\label{prop:edgeact}
The action of $\Phi_{\Ea^\bullet(F)}$ on $H^*(W,\CC)$ induced by 
the complex $\Ea^\bullet(F)$ (introduced in definition \ref{def:bigdef}) 
is given by 
\begin{equation}\label{eq:minuss}
\gamma \mapsto \gamma - \prod_{j \in I_-'} (1- e^{-q_j \mu + \mu_j')})
\int_{C_\xi} \frac{\prod_{j \in I_-''} (1- e^{-(d_j \xi - \mu_j')})}
{\prod_{j \in I_+} (1- e^{-(q_j \xi + \mu_j')})} \;\gamma (\xi) \; d\xi,
\end{equation}
where $C_\xi$ is a contour enclosing all the poles $\xi$ such that
$q_j \xi + \mu_j' = 0,$ for some $j, j \in I_+.$

\end{proposition}

\begin{proof}
According to the proposition \ref{prop:action},
the automorphism $\Phi_{\Ea^\bullet(F)}$ acts on $H^*(W,\CC)$ by
\begin{equation}\label{eq:wrw}
\gamma \mapsto \gamma - \prod_{j \in I_-'} (1-e^{[D_j]})
q_{2_*} \Bigl( (j_1 \circ q_1)^* (\gamma) \cdot  
\Td_Y \cdot q_2^* \bigl( (\Td_E)^{-1} \bigr)
\end{equation}
with the canonical maps (see diagram (\ref{diagr:big})) 
$q_i : Y \to E, j_i : E \to W, i=1,2.$ 

Proposition \ref{prop:yconstr} shows that 
\begin{equation}
\begin{split}
\Td_Y &= \Td_{E_0 \times_{Z_0} E_0} \cdot 
\prod_{j \in J_-''} (\Td_{q_1^*(E_j)} \cdot \Td_{q_2^*(E_j)})^{-1} \cdot
\prod_{j \notin J_-''} (\Td_{E_j})^{-1},\\
\Td_E &= \Td_{E_0} \cdot \prod_{j=1}^k (\Td_{E_j})^{-1},
\end{split}
\end{equation}
where, for $j \notin I_-'', 1 \leq j \leq k,$ 
we denote by $E_j$ the pull-back hypersurfaces $q_1^*(E_j)= q_2^*(E_j)$
(see (\ref{eq:pullequal})). It follows that 
\begin{equation}
\Td_Y \cdot q_2^* \bigl((\Td_E)^{-1} \bigr)= 
\prod_{j \in J_+} \Td_{q_1^*(D_j)} \cdot 
\prod_{j \in J_-''} (\Td_{q_1^*(E_j)})^{-1}.
\end{equation}

According to the construction of proposition (\ref{prop:yconstr}), the subvarieties 
$Y$ and $E$ are complete intersections in $E_0 \times_{Z_0} E_0$ and 
$E_0,$ respectively. Hence, in 
$H^*(E_0 \times_{Z_0} E_0, \CC)$ we can write 
\begin{equation}
[Y] \cdot f_2^*([E])^{-1}= \prod_{j \in I_-''} f_1^*([E_j]).
\end{equation}
Lemma \ref{lemma:map} gives then that, for a class $\alpha \in 
H^*(E_0 \times_{Z_0} E_0, \CC),$ we have 
\begin{equation}\label{eq:qwww}
\begin{split}
&(q_2)_* \Bigl( \alpha \cdot  \Td_Y \cdot q_2^* \bigl( (\Td_E)^{-1} \bigr)\Bigr)
=\\ &j^* (f_2)_* \Bigl( \alpha \cdot 
\prod_{j \in J_+} \Td_{f_1^*(D_j)} \cdot 
\prod_{j \in J_-''} (1- e^{-[E_j]})\Bigr),
\end{split}
\end{equation}
where $j$ is the inclusion $j : E \to E_0.$ 

Finally,
proposition \ref{prop:toricfibr} applied to the toric fibration 
$f_2 : E_0 \times_{Z_0} E_0 \to E_0$ shows that 
$(q_2)_* \Bigl( (j_1 \circ q_1)^* (\gamma) \cdot  
\Td_Y \cdot q_2^* \bigl( (\Td_E)^{-1} \bigr)$
can be written in terms 
of the chosen generators for $H^*(X_{\Sigma_0}, \CC)$ as 
\begin{equation}
\frac{1}{2 \pi i} \int_{C_\xi} 
\frac{\prod_{j \in I_-''} (1- e^{-(d_j \xi - \mu_j')})}
{\prod_{j \in I_+} (1- e^{-(q_j \xi + \mu_j')})} \;\gamma (\xi) \; d\xi,
\end{equation}
where $C_\xi$ is a contour enclosing all the poles $\xi$ such that
$q_j \xi + \mu_j' = 0,$ for some $j, j \in I_+.$

The result follows after combining this last formula with (\ref{eq:wrw}),
since for $j \in I_-',$ $[D_j]= -q_j \mu + \mu_j'.$ 

\end{proof}

\begin{thm} \label{thm:genthm} For any edge $F$ of the 
secondary polytope with one vertex corresponding to the 
smooth phase $W,$ the action of $\Phi_{\Ea^\bullet(F)}$ on $H^*(W,\CC)$
induced by the complex $\Ea^\bullet(F)$ (given by definition 
\ref{def:bigdef}) is identified under mirror symmetry
with the monodromy action
on $H^*(M,\CC)$ of a loop homotopic to a loop included in
the rational curve determined by the edge $F$ in the 
moduli space of complex structures on the mirror Calabi--Yau manifold $M.$
\end{thm}

\begin{proof} The proof reduces to a calculation very similar to 
the one performed in the beginning of section 
\ref{chap:1par} and it is designed to match the formulae of the proposition
\ref{prop:edgeact} and theorem \ref{thm:mon}.

Recall (\ref{phicone}) that the solutions to the GKZ system are given by 
the series 
\begin{equation}\label{eq:phi2}
\Phi^{\Ca}_{\lambda}(z):=
\sum_{l\in\Ca^{\vee}\cap \LL }\;
\prod_{j=1}^k \frac{z_j^{\lambda_j+l_j-1}} {\Gamma (\lambda_j+l_j)} 
\prod_{j=k+1}^N \frac{z_j^{\lambda_j+l_j}} {\Gamma (\lambda_j+l_j+1)},
\end{equation}
where the parameter $\lambda=(\lambda_1,\ldots,\lambda_N)\in\CC^N$
satisfies (\ref{phicond}), 
and $\Ca^\vee \subset \LL_\RR$ is the dual of the cone 
$\Ca \subset \LL^{\vee}_\RR,$ given by 
\begin{equation}
\Ca^\vee :=\{ l :
\langle w,l \rangle \geq 0, \ \hbox{for all} \ w \in \Ca \}.
\end{equation}
As explained in section \ref{chap:mon}, the choice of generators 
for $H^*(W,\CC)$ implies that the product 
\begin{equation}
\prod_{j=1}^k \frac{1} {\Gamma (\lambda_j+l_j)} 
\prod_{j=k+1}^N \frac{1} {\Gamma (\lambda_j+l_j+1)}, \; 
l \in \Ca^{\vee}\cap \LL
\end{equation}
is written as 
\begin{equation}
\prod_{j=1}^k \frac{1} {\Gamma (-d_j (m+\mu) + l_j'+\mu_j'+1)} 
\prod_{j=k+1}^N \frac{1} {\Gamma (q_j (m+\mu) +l_j'+\mu_j'+1)}
\end{equation}
for some integers $m, l_j'.$

Consider now the series 
\begin{equation}
\Psi(z):= \sum_{l\in\Ca^{\vee}\cap \LL }
\frac{\prod_{j=1}^k \Gamma (-\lambda_j-l_j+1)}
{\prod_{j=k+1}^N \Gamma (\lambda_j+l_j+1)} \;
\prod_{j=1}^k 
z_j^{\lambda_j+l_j-1} \prod_{j=k+1}^N z_j^{\lambda_j+l_j}.
\end{equation}
In terms of explicit generators, each coefficient 
\begin{equation}
\frac{\prod_{j=1}^k \Gamma (-\lambda_j-l_j+1)}
{\prod_{j=k+1}^N \Gamma (\lambda_j+l_j+1)}
\end{equation}
is written as
\begin{equation}
\frac{\prod_{j=1}^k \Gamma (d_j (m+\mu) - l_j'- \mu_j')}
{\prod_{j=k+1}^N \Gamma(q_j (m+\mu) +l_j'+\mu_j'+1)}.
\end{equation}
A very similar calculation to (\ref{eq:phipsi}) shows that, possibly after
some change of the coordinates $z_j$ to some new coordinates
$\tilde{z}_j$ differing only by complex phase changes, the series $\Psi$ 
and $\Phi$ are related by 
\begin{equation} 
\Phi^{\Ca}_{\lambda}(\tilde{z})= \prod_{j=1}^k 
\frac{1-e^{-2 \pi i (d_j \mu - \mu_j')}}{2 \pi i} \Psi(z).
\end{equation}
Proposition \ref{stienstra} implies that in fact, the
$\Psi$-series gives all the periods corresponding to the
mirror Calabi--Yau $M$ as the coefficients of the 
monomials in the series given by the generators of $H^*(W,\CC).$
The monodromy formula for the series $\Psi$ is obtained from
the monodromy formula for the series $\Phi$ given by 
theorem \ref{thm:mon}. Up to a conjugation in the fundamental group
(or, equivalently a change in the coordinates $z_j$), the calculation from 
section \ref{chap:1par} mentioned before shows that the analytic 
continuation of $\Phi(z)$ around the prescribed loop is 
\begin{equation}
\Psi(z) - \frac{1}{2 \pi i} \prod_{j \in I_-'} 
(1- e^{2 \pi i(-q_j \mu + \mu_j')})
\int_{C_\xi} \frac{\prod_{j \in I_-''} (1- e^{-2 \pi i (d_j \xi - \mu_j')})}
{\prod_{j \in I_+} (1- e^{-2 \pi i(q_j \xi + \mu_j')})} \;\Psi(\xi) \; d\xi.
\end{equation}
A change of variable which replaces any generator of $H^*(W,\CC)$
of the form $2\pi i \alpha$ by the 
generator $\alpha$ ends the proof of the theorem.

\end{proof}


\end{document}